\documentclass[authoryear,square,preprint,12pt,3p,times]{elsarticle}

\usepackage{graphicx}
\usepackage{wrapfig}

\usepackage{amssymb}

\usepackage{url}

\newcommand {\iii} {3$^{\rm rd}$}

\newcommand {\ix} {9$^{\rm th}$}

\newcommand {\xiii} {13$^{\rm th}$}
\newcommand {\xv} {15$^{\rm th}$}
\newcommand {\xvi} {16$^{\rm th}$}
\newcommand {\xvii} {17$^{\rm th}$}
\newcommand {\xviii} {18$^{\rm th}$}
\newcommand {\xix} {19$^{\rm th}$}
\newcommand {\xx} {20$^{\rm th}$}
\newcommand {\xxi} {21$^{\rm st}$}

\newcommand {\first} {1^{\rm st}}
\newcommand {\second} {2^{\rm nd}}
\newcommand {\third} {3^{\rm rd}}
\newcommand {\fourth} {4^{\rm th}}
\newcommand {\fifth} {5^{\rm th}}
\newcommand {\sixth} {6^{\rm th}}

\newcommand {\tenth} {10^{\rm th}}
\newcommand {\ith} {i^{\rm th}}
\newcommand {\jth} {j^{\rm th}}
\newcommand {\kth} {k^{\rm th}}

\newcommand {\Aryabhatiya} {\={A}ryabha\d{t}\={\i}ya}
\newcommand {\Aryabhata} {\={A}ryabha\d{t}a}

\newcommand {\minus} {\;-\;}
\newcommand {\plus} {\;+\;}
\newcommand {\equals} {\;=\;}
\newcommand {\back}{\hspace*{-0.5em}}

\newcommand {\next} {\\ \hspace*{1em} \hfill}
\newcommand {\same} {\hfill}

\newcommand {\citeibid} [2] [] {[ibid., #1]}

\journal{arXiv,org. \hfill  \begin {minipage} [t]  {3in} \rm Please read and cite the corrected, published article that can be accessed through the DOI link on the archive page for this document.\end {minipage}}

\begin{document}

\begin{frontmatter}



\title {How Ordinary Elimination Became Gaussian Elimination}


\author{Joseph F. Grcar}
\ead{jfgrcar@na-net.ornl.gov}
\address{6059 Castlebrook Drive, Castro Valley, CA 94552-1645}

\begin {abstract}
Newton, in notes that he would rather not have seen published, described a process for solving simultaneous equations that later authors applied specifically to linear equations. This method --- that Euler did not recommend, that Legendre called ``ordinary,'' and that Gauss called ``common'' --- is now named after Gauss: ``Gaussian'' elimination. Gauss's name became associated with elimination through the adoption, by professional computers, of a specialized notation that Gauss devised for his own least squares calculations. The notation allowed elimination to be viewed as a sequence of arithmetic operations that were repeatedly optimized for hand computing and eventually were described by matrices. 

\bigskip

\noindent
In Aufzeichnungen, die Newton lieber nicht der Ver\"offentlichung preisgegeben h\"atte, beschreibt er den Prozess f\"ur die L\"osung von simultanen Gleichungen, den sp\"atere Autoren speziell f\"ur lineare Gleichungen anwandten. Diese Methode --- welche Euler nicht empfahl, welche Legendre ``ordinaire'' nannte, und welche Gau{\ss} ``gew\"ohnlich'' nannte --- wird nun nach Gau{\ss} benannt: Gau{\ss}sches Eliminationsverfahren. Die Verbindung des Gau{\ss}schen Namens mit Elimination wurde dadurch hervorgebracht, dass professionelle Rechner eine Notation \"ubernahmen, die Gau{\ss} speziell f\"ur seine eigenen Berechnungen der kleinsten Quadrate ersonnen hatte, welche zulie{\ss}, das Elimination als eine Sequenz von arithmetischen Rechenoperationen betrachtet wurde, die wiederholt f\"ur Handrechnungen optimisiert wurden und schlie{\ss}lich durch Matrizen beschrieben wurden. 

\bigskip
\end {abstract}

\begin {keyword}
algebra before 1800 \sep Gaussian elimination \sep human computers \sep least squares method \sep mathematics education




\MSC 01-08 \sep 62J05 \sep 65-03 \sep 65F05 \sep 97-03


\end {keyword}

\end{frontmatter}



\renewcommand {\baselinestretch} {1.0}


\section {Introduction}
\label {sec:part1}

\subsection {Overview}

The familiar method for solving simultaneous linear equations, Gaussian elimination, originated independently in ancient China and early modern Europe. The details of the Chinese development may be lost because of a tradition of anonymity in the primary text. In contrast, the more recent European origin is traceable through its appearance in comparatively few significant documents produced over a period of some five hundred years. This paper explains the European history of Gaussian elimination up to the invention of electronic computers. 

The development has three phases: first came the ``schoolbook'' rule, second were methods that professional hand computers used primarily for least squares work, and third was the adoption of matrix notation which reconciled the schoolbook and professional methods. Section \ref {sec:part1} summarizes the history and connects it to the present by documenting where in education and technology Gaussian elimination is found today. 

Section \ref {sec:part2} examines phase one, the period before Gauss during which Gaussian elimination was invented and became a staple of algebra textbooks. Beginning with a discussion of the Chinese origin of the technique, it continues on to the European origin in the late Renaissance and an examination of the contribution of Newton and Rolle. 

Section \ref {sec:part3} introduces phase two, the period during which the need to solve simultaneous linear equations, and solve them repeatedly, arose for a clear social purpose: statistical inference. After Legendre and Gauss invented the method of least squares, Gauss started a tradition of seeking to improve the methods of calculation. 

Section \ref {sec:part4} continues phase two into the period after Gauss and describes the efforts of professional hand computers to ease their labor. This section focuses on the work of three men: Myrick Doolittle, a computer working at the United States Coast and Geodetic Survey; Andr\'e Cholesky, a French military geodesist; and Prescott Crout, a professor of mathematics at the Mathematics Institute of Technology.

Section \ref {sec:part5} treats phase three, the matrix explications of Gaussian elimination that arose for didactic purposes in the first half of the \xx\ century. It was understood from the matrix descriptions that the schoolbook methods of phase one, and the professional methods of phase two, were all essentially related. This phase begins with the Cracovian algebra of the astronomer, geodesist, and mathematician Tadeusz Banachiewicz. The latter was an important influence on the geodesist Henry Jensen whose contribution is also discussed, along with that of the statistician Paul Sumner Dwyer and the mathematician Ewald Bodewig. The section concludes with a short discussion of the work of John von Neumann and Herman Goldstine. 

\subsection {Gaussian Elimination Today}
\label {sec:today}

Both elementary and advanced textbooks discuss Gaussian elimination. For example, the precalculus algebra textbook of \citet [743-747, sec.\ 10.2] {Cohen2006} shows how ``elementary operations'' on equations produce an ``equivalent system'' in ``upper-triangular form'' that can be solved by ``back-substitution'' (Figure \ref {fig:today}, equation (\ref {eqn:Cohen})). The distinguishing features are as follows.
\begin {enumerate}
\item The equations and variables may be rearranged so the leading equation contains whichever variable is chosen to be the leading variable for immediate elimination.\footnote {A large coefficient is desired for the leading variable in the leading equation.}
\item The leading equation is used (in various ways) to remove the leading variable from each of those following.
\item These steps are then applied recursively to the following (modified) equations.
\item Once a single equation remains, then the back-substitution is made. 
\end {enumerate}
The form of the algorithm employed in equation (\ref {eqn:Cohen}) is viewed as canonical: the leading equation remains unchanged while variables are removed by subtracting appropriate multiples of it from each equation below. This paper follows current usage by referring to any algorithm that is essentially equivalent to equation (\ref {eqn:Cohen}) as ``Gaussian elimination'' whatever its period or source; the phrase should not be interpreted as an attribution.

\begin {figure}
{%
\footnotesize%
Schoolbook method:}
\begin {equation}
\label {eqn:Cohen}
\newcommand{\myspace}{\hspace*{0.5em}}
\renewcommand {\arraystretch} {1.125}
\setlength {\arraycolsep} {0.125em}
\left.
\begin {array} {r c r}
x + 2y + z & =& 3\\
x + y + 2z& =& 9\\
2x + y + z& =&16\\
\end{array}
\myspace \right\} \Rightarrow \myspace
\left.
\begin {array} {r c r}
x + 2y + z &= &3\\
- y + z& = &6\\
-3y -z& =& 10\\
\end{array}
\myspace \right\} \Rightarrow \myspace
\left.
\begin {array} {r c r}
x + 2y + z &= &3\\
- y + z& = &6\\
-4z& =& -8\\
\end{array}
\myspace \right\} \Rightarrow \myspace
\begin {array} {r c r}
z& =& 2\\
y& = & -4\\
x& =& 9
\end{array}
\end{equation}
{%
\footnotesize%
``Classic'' professional method:}
\vspace*{-1ex}
\begin {equation}
\newcommand {\depth} {\vrule depth0.5ex height0pt width0pt}
\label {eqn:matrix}
\left.
{
\setlength {\arraycolsep} {0.25em}
\begin {array} {c}
\underbrace {
\scriptsize
\hspace*{0.5em}
\left[ \begin {array} {c c c} 1&& \\ & 1& \\ & 3& 1 \end {array} \right]^{-1} \,
\left[ \begin {array} {c c c} 1&& \\ 1& 1& \\ 2&& 1 \end {array} \right]^{-1} \,
\left[ \begin {array} {c c c} 1& 2& 1\\ 1& 1& 2\\ 2& 1& 1 \end {array} \right] =
\left[ \begin {array} {c c c} 1& 2& 1\\ & -1& 1\\ && \back-4\end {array} \right]_{\depth}
\hspace*{0.5em}
}
\\[3.0ex]
\Downarrow
\\[0.5ex]
\scriptsize
\left[ \begin {array} {c c c} 1& 2& 1\\ 1& 1& 2\\ 2& 1& 1 \end {array} \right]
=
\left[ \begin {array} {c c c} 1&& \\ 1& 1& \\ 2& 3& 1 \end {array} \right] \,
\left[ \begin {array} {c c c} 1& 2& 1\\ & \back-1& 1\\ && \back-4\end {array} \right]
\end {array}
}
\vrule depth8.5ex width0pt \right\} \Rightarrow \quad
{
\setlength {\arraycolsep} {0.25em}
\begin {array} {c}
\underbrace {
\scriptsize
\hspace*{0.5em}
\left[ \begin {array} {c c c} 1& 2& 1\\ 1& 1& 2\\ 2& 1& 1 \end {array} \right]
\left[ \begin {array} {c} x\\ y\\ z \end {array} \right] =
\left[ \begin {array} {c} 3\\ 9\\ 6 \end {array} \right]_{\depth}
\hspace*{0.5em}
}
\\[3.0ex]
\Downarrow
\\[0.5ex]
\scriptsize
\left[ \begin {array} {c c c} 1&& \\ 1& 1& \\ 2& 3& 1 \end {array} \right] \,
\left[ \begin {array} {c c c} 1& 2& 1\\ & \back-1& 1\\ && \back-4\end {array} \right]
\left[ \begin {array} {c} x\\ y\\ z \end {array} \right] =
\left[ \begin {array} {c} 3\\ 9\\ 6 \end {array} \right]
\end {array}
}
\vspace*{-3ex}
\end {equation}
\caption {Gaussian elimination in precalculus algebra and as presently understood by computing specialists.}
\label {fig:today}
\end {figure}

\nocite {Waugh1945}

In contrast, the technical literature views Gaussian elimination as a method for factoring matrices. Elementary operations ``reduce'' the coefficient matrix of equation (\ref {eqn:Cohen}) to an ``upper-triangular'' matrix thereby accomplishing a ``triangular factorization,'' or ``decomposition,'' from which the equations can be solved by ``forward substitution'' followed by backward substitution (Figure \ref {fig:today}, equation (\ref {eqn:matrix}), but note that the arithmetical operations are identical to those in equation (\ref {eqn:Cohen})). The advanced textbook of \citet [107] {Petersen2004} explains that this algorithm is the standard test for the speed of computers in scientific work. Its widespread use in so large a field as scientific computing results in many algorithmic variants that are collectively called, simply, Gaussian elimination. The variations are designated by acronyms, adjectives, and eponyms. At this level of differentiation the version of equation (\ref {eqn:Cohen}) is named either ``classic'' elimination or ``Doolittle's method.'' Nevertheless, advanced or specialized texts always begin by identifying exactly this algorithm as Gaussian elimination.\footnote {For example see \citet [43--44] {Duff1986}, \citet [94] {Golub1996}, \citet [158--160] {Higham2002}, \citet [23--25] {Petersen2004}, and \citet [148--154] {Stewart1998}.}

Today, textbooks at all levels encourage the belief that Gauss introduced the method of equation (\ref {eqn:Cohen}) and that his work was somehow remarkable when compared to prior art. For example, \citet [743] {Cohen2006} claim that ``Gauss used this technique'' to analyze the orbit of Pallas \citep {Gauss1810-HM} although ``the essentials'' had already appeared in ancient Chinese texts. \citet [p.\ 187] {Higham2002} acknowledges that ``a variant'' was used by the Chinese, but he finds ``the first published appearance of Guassian elimination'' in another paper of \citet {Gauss1809-HM}. Some authors do allude to a possible European origin before Gauss. \citet [3] {Farebrother1988} makes the remark that ``Gauss's formalization'' of the ``traditional schoolbook method'' appeared in his Pallas work. These inconsistent attributions beg the question: what exactly did Gauss contribute?  

\section {Before Gauss}
\label {sec:part2}

\subsection {Gaussian Elimination in Antiquity or Not}
\label {sec:antiquity}

Periods before Gauss are surveyed here to see where and when Gaussian elimination may have developed. It will be seen that consideration of any simultaneous linear equations is comparatively rare in the primary sources. Moreover, almost all early sources consist of terse ``word problems'' rather than discourses on solution methods.


The earliest written mathematics is believed to be ``Babylonian'' cuneiform tablets from the Tigris and Euphrates valley. \citet [343--344] {Robson2008} counts around a thousand published tablets that record mathematics taught in scribal schools. Her inventory shows that most of these tablets contain arithmetical or metrological tables. Tablets VAT 8389 and 8391 from the Old Babylonian period, 2000--1600 BC, contain what are believed to be among the earliest problems that can be interpreted as systems of linear equations.\footnote {This dating is the widest possible based on \citet [6, 321] {Robson1999}. \citet [440--441] {Friberg2007amazing} explains the acquisition of VAT 8389 and 8391 by the Vorderasiatisches Museum in Berlin and their lost provenance. \citet [II, plates 24, 26] {Neugebauer1973-HM} provides photographs of the tablets.} 

The actual text of problem 1 on VAT 8389 and a literal restatement in modern English are given by \citet [77--82] {Hoyrup2002}. For brevity, equation (\ref {eqn:Babylon}) of Figure \ref {fig:ancient} expresses the solution symbolically; see also \citet [334-335] {Friberg2007remarkable} and \citet [16] {Katz1998}. The frequent admonition by the anonymous scribe to ``keep it in your head'' \citep [see] [79] {Hoyrup2002} is poignant evidence of the computing technology available to him. \citet [82] {Hoyrup2002} comments that the solution method is pedagogically superior to the method of double false position found later. \citet [334] {Friberg2007remarkable} reports that relatively few tablets pose linear problems. Indeed, the tablets summarized by \citet [3] {Bashmakova2000-HM} suggest, that when ``Babylonian'' problems can be interpreted as simultaneous equations, usually one condition is nonlinear.


\begin {figure}
\newcommand {\beginequation} {\vspace*{-1.0ex}}
\newcommand {\size} {\footnotesize \setlength {\baselineskip}{12pt}}
\newcommand {\halfem} {\hspace*{0.5em}}
{%
\size%
\textbf {2000--1600 BC:} Problem 1 on tablet VAT 8389: the total area of two fields is 1800 sar, the rent for one is 2 sil\`a of grain per 3 sar, for the other is 1 sil\`a per 2 sar, the total rent on the first exceeds that on the other by 500 sil\`a.\par
}
\beginequation
\begin {equation}
{%
\size%
\label {eqn:Babylon}
\setlength {\arraycolsep} {0.125em}
\begin {array} {c c c c c c c c c c}
x& +& y& \equals& 1800\\[0.5ex]
{2 \over 3} \cdot x& -& {1 \over 2} \cdot y& =& 500\\
\end {array}
\halfem \Rightarrow \halfem
\begin {array} {c c c c l}
(900 + s)& +& (900 - s)& \equals& 1800\\[0.5ex]
{2 \over 3} \cdot (900 +s)& -& {1 \over 2} \cdot (900 - s)& =& 150 + {7 \over 6} \cdot s 
\halfem \Rightarrow \halfem
s = {6 \over 7} \cdot \mbox {350, hence $x$, $y$}\\
\end {array}
}
\end {equation}
{%
\size%
\textbf {1550 BC:} Problem 40 of the Rhind papyrus: divide 100 into 5 parts so the common difference is the same and the three largest are 7 times the other two \citep [170--172] {Gillings1972}. Solution: an arithmetic progression is understood from the phrase ``common difference.'' The scribe sums $\{$23, 17$1\over2$, 12, 6$1\over2$, 1$\}$, which he somehow knows beforehand satisfies the $7:1$ condition, and upon obtaining 60, he then scales the terms by $5\over3$ to enforce the 100 condition.\par
}
\beginequation
\begin {equation}
{%
\size%
\label {eqn:Rhind}
\newcommand {\depth} {\vrule depth0.5ex height0pt width0pt}
\underbrace {\depth \, 23 + 17.5 + 12 \,}_{\mbox {$52.5 = 7 \times 7.5$}} + \underbrace {\depth\,6.5 + 1 \,}_{\mbox {7.5}} 
= \mbox {60} 
\quad \Rightarrow \quad \{
\mbox {$38 {1 \over 3}$}, \;
\mbox {$29 {1 \over 6}$}, \;
\mbox {$20$}, \;
\mbox {$10 {5 \over 6}$}, \;
\mbox {${5 \over 3}$}
\}
}
\end {equation}
{%
\size%
\textbf {Before 200 AD:} Problem 1 in chapter 8 of the \textit {Nine Chapters\/} \citep [391 and 399--403] {Kangshen1999}: from 3 top-grade rice paddies, 2 medium-grade, and 1 low-grade, the combined yield is 39 dou of grain, and so on for two more collections of paddies. What is the yield of a paddy of each grade? Solution:\par
}
\beginequation
\begin {equation}
\label {eqn:China}
{%
\size%
\newcommand{\myspace}{\hspace*{0.5em}}
\renewcommand {\arraystretch} {1.125}
\setlength {\arraycolsep} {0.125em}
\begin {array} {r c r}
3x + 2y + z&= & 39\\
2x + 3y + z& = & 34\\
x + 2y + 3z& =& 26\\
\end{array}
\myspace \Leftrightarrow \myspace
\begin {array} {r r r}
1& 2& 3\\
2& 3& 2\\
3& 1& 1\\
26& 34& 39
\end{array}
\myspace \Rightarrow \myspace
\begin {array} {r r r}
& & 3\\
4& 5& 2\\
8& 1& 1\\
39& 24& 39
\end{array}
\myspace \Rightarrow \myspace
\begin {array} {r r r}
& & 3\\
& 5& 2\\
36& 1& 1\\
99& 24& 39
\end{array}
\myspace \Rightarrow \myspace
\renewcommand {\arraystretch} {1.5}
\begin {array} {r c r}
z& =& {11 \over 4}\\
y& = & {17 \over 4}\\
x& =& {37 \over 4}\\[1.0ex]
\end{array}
}
\end{equation}
{%
\size%
\textbf {250 AD:} Problem 19 in book 1 of the \textit {Arithmetica\/} of Diophantus \citep [136] {Heath1910}: find four numbers such that the sum of any three exceeds the fourth by a given amount. Solution (using symbols for clarity): let $n_i$ be the numbers, $2h$ the sum of all four, and $e_i$ the excesses. Then $2h - n_i = n_i + e_i$ or $n_i = h - e_i/2$. Summing the latter for all four numbers gives $2h = 4h - (e_1 + e_2 + e_3 + e_4) / 2$. Thus $h$ and then the $n_i$ can be evaluated. As the problem would be expressed today:\par
}
\beginequation
\begin {equation}
{%
\size%
\label {eqn:Diophantus}
\setlength {\arraycolsep} {0.125em}
\begin {array} {c c c c c c c c c c}
-& n_1& +& n_2& +& n_3& +& n_4& \equals& e_1\\
& n_1& -& n_2& +& n_3& +& n_4& =& e_2\\
& n_1& +& n_2& -& n_3& +& n_4& =& e_3\\
& n_1& +& n_2& +& n_3& -& n_4& =& e_4\\
\end {array}
\quad \Rightarrow \quad
\mbox {\size $\displaystyle n_i \equals {e_1 + e_2 + \cdots + e_n \over 4} - {e_i \over 2}$}
}
\end {equation}
{%
\size%
\textbf {499 AD:} Problem 29 in chapter 2 of the \textit {\Aryabhatiya\/} of \Aryabhata\ \citep [40] {Clark1930} \citep [134] {Plofker2009}: to find several numbers when the results of subtracting each from their sum are known. Solution: sum the known differences and divide by the quantity of terms less one. The result is the sum of all the numbers, from which they can be determined.\par
}
\beginequation
\begin {equation}
\label {eqn:Aryabhata}
{%
\size%
\setlength {\arraycolsep} {0.125em}
\begin {array} {c c c c c c c c c c c}
0& +& n_2& +& n_3& +& \cdots& +& n_n& \equals& d_1\\
n_1& +& 0& +& n_3& +& \cdots& +& n_n& \equals& d_2\\
n_1& +& n_2& +& 0& +& \cdots& +& n_n& \equals& d_3\\
\vdots&& \vdots&& \vdots&& \ddots&& \vdots&& \vdots\\
n_1& +& n_2& +& n_3& +& \cdots& +& 0& \equals& d_n\\
\end {array}
\quad \Rightarrow \quad
\mbox {\size $\displaystyle n_i \equals {d_1 + d_2 + \cdots + d_n \over n-1} - d_i$}
}
\end {equation}
\caption {Representatives of all the known ancient problems, and their solutions, that might be interpreted in symbolic algebra as simultaneous linear equations.}
\label {fig:ancient}
\end {figure}

One of the main sources for Egyptian mathematics is the Rhind papyrus which dates from around 1550 BC.\footnote {The Rhind papyrus originated in Thebes and is now in the British Museum. For dating, provenance, and a view of the papyrus see the museum website: \url {http://www.britishmuseum.org/explore/highlights/highlight_objects/aes/r/rhind_mathematical_papyrus.aspx}.} Problem 40 in the Rhind papyrus (Figure \ref {fig:ancient} equation (\ref {eqn:Rhind})) is sometimes mentioned as an early example of simultaneous linear equations \citep [303] {Dedron1959}. \citet [170--172] {Gillings1972}, however, argues that the solution is based on knowledge of arithmetic progressions, the latter being a distinguishing feature of mathematics in ancient Egypt. 

By far the most impressive treatment of simultaneous linear equations known from antiquity --- and the only lengthy presentation prior to the late \xvii\ century --- is chapter 8 of the \textit {Jiuzhang Suanshu\/} (\textit {Nine Chapters of the Mathematical Art\/}), a problem ``book'' anonymously and collectively written in China. An early signed commentary, including comments on chapter 8, was written by Liu Hui who is known to have lived in the \iii\ century AD. The method for solving linear problems is therefore at least as old, and Hui asserts it is much older, although no original text survives. According to \citet [128--131] {Martzloff1997-HM}, chapters 1--5 are known from a \xiii\ century copy, and chapters 6--9 are reconstructed from \xviii\ century quotations of a lost \xv\ century encyclopedia. 

Mathematicians in ancient China represented numbers by counting rods. They organized elaborate calculations by placing the rods inside squares arranged in a rectangle.\footnote {The Romans also arranged complicated calculations in tabular form apparently using pebbles in place of rods. Neither Chinese nor Roman artifacts survive for these ``table'' calculations. Both were supplanted by abaci.} In chapter 8 of the \textit {Nine Chapters\/}, each column of squares corresponds to a modern linear equation, so in order to obtain the coefficient tableau of modern mathematics, the ancient rectangle must be rotated counterclockwise by 90 degrees. Problem 1 in chapter 8 is frequently displayed as representative of the ancient calculation (Figure \ref {fig:ancient} equation (\ref {eqn:China})). Chapter 8 includes the solution of 18 different systems of equations all treated in the same systematic way. Unlike equation (\ref {eqn:Cohen}), when two columns of the table are combined to remove a number, each column is multiplied by the leading number in the other. Subtracting the right column from the left column then removes the leading number from the left while preserving integer coefficients.\footnote {For more information on chapter 8 see \citet {LayYong1989}.}

The main extant source for Greco-Roman algebra is the \textit {Arithmetica\/} by Diophantus. An approximate date of composition is usually given as around 250 AD, but \citet [3, 9] {Schappacher2005} believes Diophantus could have lived as late as 400 AD, and in any case, the sources of the extant texts cannot be traced any earlier than the \ix\ century.\footnote {Of 13 supposed books in the \textit {Arithmetica\/}, 1--6 are known from several copies in the original Greek \citep [5] {Heath1910}, while 4--7 are known in Arabic translation \citep {Sesiano1982}, and the rest are lost. Three manuscripts are at the Biblioteca Apostolica Vaticana in Rome: Vat.gr.191, ff.\ 360-390v, Vat.gr.200, ff.\ 1-193, and Vat.gr.304, ff.\ 77-117v. Four manuscripts are at the Biblioth\`eque nationale de France in Paris: BnF grec 2378, 2379, 2380, and 2485.}


\citet [59] {Heath1910} quotes Diophantus for a rule to solve problems with one unknown and one condition, which need not be linear: ``If a problem leads to an equation in which any terms are equal to the same terms but have different coefficients, we must take like from like.'' The use of the word ``equation'' is of course an anachronism. Many problems in the \textit {Arithmetica\/} involve more than one unknown. Heath writes that Diophantus ``shows great adroitness'' in expressing all the unknowns in terms of a single newly introduced value, ``the ultimate result being that it is only necessary to solve a simple equation with one unknown quantity.'' The process is illustrated by the verbal solution given by Diophantus for equation (\ref {eqn:Diophantus}) in Figure \ref {fig:ancient}. This problem is among the most elaborate of the few in the \textit {Arithmetica\/} that are purely linear. As Heath suggests, the problem involves a special linear system whose solution is a repetitive formula and is thus amenable to special reasoning that would not suffice for the more general problems of the \textit {Nine Chapters\/}.

According to \citet [122, 317--318] {Plofker2009}, the earliest surviving work of ancient Hindu mathematics is the \textit {\Aryabhatiya\/} of \Aryabhata\ who is known to have been born in 476 AD. \Aryabhata's linear problems are reminiscent of the \textit {Arithmetica\/} but are stated more generally for any quantity of unknowns (Figure \ref {fig:ancient}, equation (\ref {eqn:Aryabhata})). The translation by \citet {Clark1930} has \Aryabhata\ explaining how to calculate the solution but not how to find it. The repetitive conditions of the problem allow it to be reduced to a single unknown, as in the case of the problems of Diophantus.


In the medieval period, the equivalent of single polynomial equations, but notably not simultaneous equations, were solved by several Arabic-speaking mathematicians. Examples are found in the work of the encyclopedist Al-Khwarizmi of Baghdad and the writings of Leonardo of Pisa (Fibonacci) who travelled in the western Arabic world \citep [78] {Hogendijk1994}. 


Symbolic algebra developed during the European Renaissance in the arithmetization of geometry and in the theory of equations. By the end of the \xvi\ century, an audience had developed for textbooks that taught arithmetic, how to express ``questions'' (word problems in modern parlance) as symbolic equations, and methods for solving them. To obtain a comprehensive picture of algebra in the late Renaissance, \citet {Kloyda1938} surveyed what may be all the algebra texts printed between 1550 and 1660. Only four of the 107 texts she examined discuss simultaneous linear equations, which suggests such equations were of little interest during this period. 

\begin {figure}[b]
\newcommand {\beginequation} {\vspace*{-1.0ex}}
\newcommand {\size} {\footnotesize \setlength {\baselineskip}{12pt}}
{%
\size%
Problem and solution of \citet [111--112] {Peletier1554-HM}, with back substitution from the $\tenth$, $\fifth$, and $\first$ equations:\par
}
\beginequation
\begin {equation}
{%
\size%
\label {eqn:Peletier}
\setlength {\arraycolsep} {0.125em}
\newcommand{\myspace}{\hspace*{0.5em}}
\begin {array} {r r r c l}
&1.& 2 R + A + B& =& 64\\
&2.& R + 3A + B& =& 84\\
&3.& R + A + 4B& = & 124\\
\noalign {\medskip}
2. + 3. \Rightarrow& 4.& 2R + 4A + 5B& = & 206\\
4. - 1. \Rightarrow& 5.& 3A+4B& =& 146\\
1. + 2. \Rightarrow& 6.& 3R + 4A + 2B& =&148\\
1. + 3. \Rightarrow& 7.& 3R + 2A + 5B& =& 188\\
6. + 7. \Rightarrow& 8.& 6R + 6A + 7B& = & 336\\
6 \times 3. \Rightarrow& 9.& \myspace 6R + 6 A + 24 B& =& 744\\
9. - 8. \Rightarrow& \myspace 10.& 17B& =& 408
\end {array}
}
\end {equation}
{%
\size%
Problem and solution of \citet [190--191] {Buteo1560-HM}, with back substitution from the $\sixth$, $\fifth$, and $\third$ equations:\par
}
\beginequation
\begin {equation}
{%
\size
\label {eqn:Buteo}
\setlength {\arraycolsep} {0.125em}
\newcommand{\myspace}{\hspace*{0.5em}}
\begin {array} {r r r c l}
&1.& 3 A + B + C& =& 42\\
&2.& A + 4B + C& =& 32\\
&3.& \myspace A + B + 5C& = & 40\\
\noalign {\medskip}
3 \times 2. - 1. \Rightarrow& 4.& 11 B + 2C& = & 54\\
3 \times 3. - 1. \Rightarrow& 5.& 2 B + 14 C& =& 78\\
11 \times 5. - 2 \times 4. \Rightarrow & \myspace 6.& 150C& =&750
\end {array}
}
\end {equation}
\caption {In the compilation by \citet {Kloyda1938}, these are two of the apparently only four Renaissance examples of solving linear systems by symbolic means, restated in modern notation. Peletier wrote an infix ``p.''\ (piu) for modern $+$, also ``m.''\ (meno) for $-$, and for equality he wrote the word. Buteo wrote ``.''\ or ``,'' for $+$, and ``['' for $=$. The problems appear to have been chosen to avoid negatives.}
\label {fig:renaissance}
\end {figure}

The four authors found by Kloyda who did solve linear systems display a rapid exploration of the possibilities of symbolic algebra. \citet [107--112] {Peletier1554-HM} illustrates the power of symbolic reasoning by solving a problem in two different ways. About three men and their money, Peletier attributed the problem to \citet [71--73] {Cardano1545} where it can be found solved by a mixture of verbal reasoning and equations in two variables. Peletier reviewed that argument, and then passed to a purely symbolic solution beginning from three variables and three equations that most directly express the conditions of the problem (Figure \ref {fig:renaissance}, equation (\ref {eqn:Peletier})). For no apparent reason he took seven steps where three are needed. \citet {Buteo1560-HM} solved a similar problem by Gaussian elimination using three steps (Figure \ref {fig:renaissance}, equation (\ref {eqn:Buteo})). The method was the same double-multiply elimination used in the \textit {Nine Chapters\/}. The solution is from the $\sixth$, $\fifth$, and $\third$ equations whereas the \textit {Nine Chapters\/} and the canonical elimination of equations (\ref {eqn:Cohen}) would retain the $\sixth$, $\fourth$, and $\first$ equations. \citet [193--196] {Buteo1560-HM} also solved a problem with four equations and unknowns that is repeated by \citet [82 reverse to 84 obverse] {Gosselin1577-HM}. Eighty years later, \citet [85--86] {Rahn1659-HM} solved a system of three linear equations. The four authors found by Kloyda provided examples of Gaussian elimination, but they did not explicitly state the algorithm.

The few problems discussed here and shown in Figures \ref {fig:ancient} and \ref {fig:renaissance} are representative of the systems of linear equations solved from ancient times up until the latter years of the \xvii\ century. During this period, simultaneous linear equations were seldom considered, and when they were, Gaussian elimination was absent except in the \textit {Nine Chapters\/} and in the work of Peletier, Buteo, Gosselin, and Rahn. Since the latter are unlikely to have had been in receipt of communications from China, Gaussian elimination apparently developed independently in Europe as a natural exercise in symbolic algebra. A statement of the method and uses for it were yet to come.

\subsection {``this bee omitted by all that have writ introductions''}
\label {sec:Newton}

Isaac Newton's work on algebra extended roughly from his appointment to the Lucasian professorship in 1669 until he began composing ``the Principia'' in 1684.\footnote {Newton's mathematical papers have been published by \citet {Whiteside1964-1967, Whiteside1968-1982}. Especially relevant to simultaneous equations in the latter set are volume 2 for the years 1667--1670 and volume 5 for 1683--1684.} In 1669--1670 Newton wrote amendments for a Latin version of a Dutch algebra text, authored by \citet {Kinckhuysen1661-HM}, which John Collins planned to publish in England. Collins abandoned the project when other books appeared.\footnote {A cursory inspection of Google books reveals that between 1650 and 1750 at least 40 algebra textbooks were published in England.} Newton himself lectured on algebra at Cambridge for eleven years beginning with the 1673--1674 academic term.\footnote {\citet [82] {Clark2006} describes the academic politics peculiar to Cambridge and Oxford that resulted in professorial lectures, including Newton's, being given ``to the walls.''} During that time, he wrote and repeatedly revised an incomplete manuscript for his own algebra treatise that was to be named \textit {Arithmetic{\ae}\ Universalis\/}. His last algebra manuscript was prepared in 1684 when, for unknown reasons, he suddenly honored the requirements of the Lucasian professorship by depositing with Cambridge University his lectures for the algebra course. The bulk of the notes were transcribed by his secretary from his previous algebra manuscripts. After Newton left academic life, his lectures were published in their original Latin (1707, 1722) and in translation (1720, 1728) under the intended title of his aborted treatise, \textit {Universal Arithmetic\/}. Newton had no claim to material that the university had paid him to prepare, nevertheless as explained by \citet [V, 8--11] {Whiteside1968-1982}, he strongly objected to its publication in case it should be misinterpreted as representing his latest research.\footnote {The titles pages of \citet {Newton1707, Newton1720-HM} pointedly identify no author.} The second English edition, incorporating changes from Newton, appeared the year after his death.

Thanks to Whiteside's impressive scholarship in preparing Newton's papers, a portion of the textbook that is relevant to Gaussian elimination can be traced directly to Newton. In the much earlier commentary on Kinckhuysen, Newton scribbled a margin note to Collins saying he intended to rectify a lacuna in contemporary textbooks: an explanation for solving collections of equations.
\begin {quotation}
\noindent Though this bee omitted by all that have writ introductions to this Art, yet I judge it very propper \& necessary to make an introduction compleate. 
\next --- Isaac Newton, margin note, circa 1669--1670
\next quoted by \citet [II, 400, n.\ 62] {Whiteside1968-1982}
\end {quotation}
Newton proposed to insert a new chapter that first explained the overall strategy of solving simultaneous equations and then listed the tactics by which it might be accomplished.\footnote {The new chapter for Kinckhuysen's textbook can be found in \citep [II, 400--411, parallel Latin and English texts] {Whiteside1968-1982}. The material appears also in the transcribed lecture notes \citep [V, 122--129, parallel texts] {Whiteside1968-1982}. That text was copied into the Latin and English editions of the unauthorized textbook.}
\begin {quotation}
\noindent \textbf {Of the Transformation of two or more {\AE}quations into one, in order to exterminate the unknown Quantities.}

\medskip
When, in the Solution of any Problem, there are more {\AE}quations than one to comprehend the State of the Question, in each of which there are several unknown Quantities; those {\AE}quations (two by two, if there are more than two) are to be so connected, that one of the unknown Quantities may be made to vanish at each of the Operations, and so produce a new {\AE}quation\dots. And you are to know, that by each {\AE}quation one unknown Quantity may be taken away, and consequently, when there are as many {\AE}quations and unknown Quantities, all at length may be reduc'd into one, in which there shall be only one Quantity unknown. 
\next --- \citet [60--61] {Newton1720-HM} and prior years
\end {quotation}
On the pages following this rule, Newton offered several methods for removing a variable from two equations, including ``equating'' and ``substituting.'' 
\begin {quotation}
\noindent \textbf {The Extermination of an unknown Quantity by an Equality of its Values.}

\medskip
When the quantity to be exterminated is only of one Dimension in both {\AE}quations, both its Values are to be sought by the Rules already deliver'd, and the one made equal to the other. 

Thus, putting $a+x=b+y$ and $2x+y = 3b$, that $y$ may be exterminated, the first Equation will give $a+x-b=y$, and the second will give $3b-2x=y$. Therefore $a+x-b=3b-2x$, \dots

\medskip
\noindent \textbf {The Extermination of an unknown Quantity, by substituting its Value for it.}

\medskip
When, at least, in one of the {\AE}quations the Quantity to be exterminated is only of one Dimension, its Value is to be sought in that {\AE}quation, and then to be substituted in its Room in the other {\AE}quation. \dots 
\same --- \citet [61--62] {Newton1720-HM} and prior years
\end {quotation}
The context and accompanying examples make clear that Newton meant a general approach for solving simultaneous \textit {nonlinear\/} equations. Newton considered only a few linear systems. In addition to the illustration quoted above, the textbook has a problem about the composition of alloys that is solved by extermination through equality of values \citep [75--77] {Newton1720-HM}.\footnote {Newton's unpublished treatise has three linear problems about the cost of some goods (wheat, wine, and silk), about money divided among beggars, which also appears in the textbook, and the alloy problem \citep [V, 567--573, probs.\ 3, 4, 7] {Whiteside1968-1982}.}

Newton's comparatively accessible algebra textbook became, as characterized by \citet [V, 54--55, n.\ 1] {Whiteside1968-1982}, ``the most widely read and influential of his writings.'' \citet [143--144] {Macomber1923-HM} notes ``before the death of Newton there came to be a demand for suitable text books of algebra for the public schools; and during the \xviii\ century, a number of texts appeared, all more closely resembling the algebra of Newton than those of earlier writers.'' Further, \citet [132] {Macomber1923-HM} finds that Newton's rule for solving simultaneous equations is ``the earliest appearance of this method on record.''\footnote {Macomber supported her claim by her own survey of algebra textbooks from the 16--\xviii\ centuries, which was corroborated by \citet {Kloyda1938}. As noted, Newton himself in 1669--1670 believed he had supplied the missing lesson, although his text remained unpublished until 1707. Macomber knew of Rolle's work of 1690, but she evidently believed that Newton's approach more resembled what was taught in early \xx\ century schoolbooks.} She refers to Newton's statement of Gaussian elimination, but not by that name because it had not yet been named for Gauss, even as late as 1923.


During the 37 years while Newton's work lay unpublished, it was overtaken by another development. Michel \citet [42-43] {Rolle1690-HM} explained how to solve simultaneous, specifically linear, equations by ``la Methode'' of precise rules for Gaussian elimination. He arranged the calculation in two columns to highlight similarities in the two different phases of the work: algebra on the left, where he substituted formulas into variables, and arithmetic on the right, where he substituted numbers into variables (Figure \ref {fig:Rolle}). The left column begins with the original equations, and records successively smaller systems by using the first equation in each system to remove the same variable from the other equations in that system. The right column begins with the first row of each reduced system, and chooses substitutions using the same sequence of rows as on the left. Rolle worked entirely in terms of substitution. 

\begin {figure}[h!]
\vspace*{-1.0ex}
\begin {displaymath}
{%
\footnotesize%
\setlength {\arraycolsep} {0.125em}
\newcommand{\myspace}{\hspace*{0.5em}}
\begin {array} {r r r r c l}
&& \multicolumn {4} {r} {\mbox {\textit {Columne de direction\/}}}\\
\noalign {\medskip}
&&& \multicolumn {3} {r} {\mbox {Premiere Classe}}\\
&& 1)& \hspace*{3em} x + y + z& =& 6\\
&& 2)& x + y + v& =& 7\\
&& 3)& x + z + v& = & 8\\
&& 4)& y + z + v& =& 9\\
\noalign {\medskip}
&&& \multicolumn {3} {r} {\mbox {Seconde Classe}}\\
\mbox {$x$ from 1\ sub.\ into 2}& \Rightarrow& 5)& v - z& =& 1\\
\mbox {$x$ from 1\ sub.\ into 3}& \Rightarrow& 6)& v - y& =& 2\\
\mbox {$x$ from 1\ sub.\ into 4}& \Rightarrow& 7)& y + z + v& =& 9\\
\noalign {\medskip}
&&& \multicolumn {3} {r} {\mbox {Troisi\'eme Classe}}\\
\mbox {$v$ from 5\ sub.\ into 6}& \Rightarrow& 8)& z - y& =& 1\\
\mbox {$v$ from 5\ sub.\ into 7}& \Rightarrow& 9)& y + 2z& =& 8\\
\noalign {\medskip}
&&& \multicolumn {3} {r} {\mbox {Quartri\'eme Classe}}\\
\mbox {$z$ from 8\ sub.\ into 9}& \Rightarrow& 10)&3y& =& 6\\
\end {array}
\quad \vrule \quad
\begin {array} {r c l l l l}
\multicolumn {3} {l} {\mbox {\textit {Columne de retour\/}}}\\
\noalign {\medskip}
\multicolumn {3} {l} {\mbox {Premiere Classe}}\\
\hspace*{0.5em} y& =& 2& (1& \Leftarrow& \mbox {$y$ from left column 10}\\
z & =& 1 + y& (2& \Leftarrow& \mbox {$z$ from left column 8}\\
v& = & 1 + z& (3& \Leftarrow& \mbox {$v$ from left column 5}\\
x& =& 6 - y - z& (4& \Leftarrow& \mbox {$x$ from left column 1}\\
\noalign {\medskip}
\multicolumn {3} {l} {\mbox {Seconde Classe}}\\
z& =& 3& (5& \Leftarrow& \mbox {$y$ from 1 sub.\ into 2}\\
v& =& 1 + z& (6& \Leftarrow& \mbox {$y$ from 1 sub.\ into 3}\\
x& =& 4 - z& (7& \Leftarrow& \mbox {$y$ from 1 sub.\ into 4}\\
\noalign {\medskip}
\multicolumn {3} {l} {\mbox {Troisi\'eme Classe}}\\
v& =& 4& (8& \Leftarrow& \mbox {$z$ from 5 sub.\ into 6}\\
x& =& 1& (9& \Leftarrow& \mbox {$z$ from 5 sub.\ into 7}\\
\noalign {\medskip}
\multicolumn {3} {l} {\mbox {Quartri\'eme Classe\hspace*{0.5em}}}\\
x& =& 1& (10& \Leftarrow& \mbox {$v$ from 8 sub.\ into 9}\\
\end {array}
}
\end {displaymath}
\caption {Gaussian elimination as performed by Michel \citet [47] {Rolle1690-HM}. Rolle did not use ``='' for equality.}
\label {fig:Rolle}
\end {figure}

It is difficult to see Rolle's influence in subsequent algebra textbooks. His emphasis on substitution may survive in the ``method of substitution'' as Gaussian elimination is sometimes called, and the name for his second column, \textit {retour\/}, may be present in \textit {back-\/} or \textit {backward\/} substitution, but these are speculations. Authors of influential French textbooks near the turn of the century, such as \citet {Bezout1788-HM}, \citet {Lacroix2nd1800-HM}, and \citet {Bourdon1828-HM}, appear to owe more to Newton than to Rolle.

\subsection {Gaussian Elimination in Textbooks After Newton}
\label {sec:afterNewton}

The influence of Newton's textbook can be traced in three signature features. First is a series of lessons for removing one variable from a pair of equations. The practice of devoting separate lessons to each elimination process became a common pedantry. Newton and Rolle taught substitution, and Newton also had ``extermination of an unknown quantity by an equality of its values'' which became ``comparison.'' Later authors added the rule used in the \textit {Nine Chapters\/} that forms linear combinations of two equations, called ``addition and/or subtraction.'' Second is a separate lesson for more than two simultaneous equations. For Newton it was a recursive rule for any number and all kinds of equations. Some later authors limited the rule to linear equations, repeated it for each type of elimination, or built on the two-equation case with separate rules for three and four equations. Third is Newton's use of the word \textit {exterminate\/}. \citet [II, 401, n.\ 63] {Whiteside1968-1982} reports that Newton first wrote ``elimino'' and then replaced it in almost all instances by ``extermino.'' The usage of this terminology grew from the \xviii\ to the \xix\ centuries and then abruptly vanished.

Among the earliest authors influenced by Newton was the banker Nathaniel Hammond \citep [44] {Macomber1923-HM}. He served as chief accountant for the Bank of England from 1760 to 1768 \citep {Roberts1995}, and his successful algebra textbook went to four editions between 1742 and 1772. Hammond got down to business by emphasizing his clear instructions for solving word problems. 
\begin {quotation}
As the principal Difficulty in this Science, is acquiring the Knowledge of solving of Questions, I have given a great Variety of these in respect to Numbers and Geometry, and their solutions I chose to give in the most particular, distinct, and plain Manner; and for which the Reader will find full and explicit Directions.
\next --- \citet [vii] {Hammond1742}
\end {quotation}
\citet [142, 219--220, and 296--297] {Hammond1742} may have been the first to replace Newton's recursive rule by a progression of rules for two, three, and four equations. The last of these is:
\begin {quotation}
\noindent \textbf {The Method of resolving Questions, which contain four Equations, and four unknown Quantities.}

\medskip
\noindent 72. When the Question contains four Equations, and there are four unknown Quantities in each Equation; find the Value of one of the unknown Quantities in one of the given Equations, and for that unknown Quantity in the other three Equations write the Value of it, which then reduces the Question to three Equations, and three unknown Quantities. 

Then find the Value of one of these three unknown Quantities in one of these three Equations, and for that unknown Quantity in the other two Equations write the Value of it, which reduces the Question to two Equations, and two unknown Quantities. 

Then find the Value of one of the unknown Quantities in each of these two Equations, and make these Equations equal to one another, when we shall have an equation with only one unknown Quantity, which being reduced, will answer the Question. \dots

And in the same Manner may any other Question in the like Circumstances be answered.
\same --- \citet [296--298] {Hammond1742}
\end {quotation}
The method for removing variables from two equations was equality-of-values, and from three or four equations it was substitution. \citet [142] {Hammond1742}, like Newton, used the word \textit {exterminating\/}.

Thomas Simpson's algebra text is another early and very prominent example of Newton's influence. This popular textbook appeared in ten editions from 1745 to 1826. The definitive second edition states Newton's recursive rule separately for each type of elimination, and uses Newton's favored word:
\begin {quotation}
\noindent \textbf {Of the Extermination of unknown quantities, or the reduction of two, or more equations, to a single one.} 

\dots. This, in most cases, may be performed various ways, but the following are the most general.

$1^\circ$. Observe which, of all your unknown quantities, is the least involved, and let the value of that quantity be found in each equation \dots\ looking upon all the rest as known; let the values thus found be put equal to each other \dots; whence new equations will arise, out of which the quantity will be totally excluded; with which new equations the operation may be repeated, and the unknown quantities exterminated, one by one, till, as last, you come to an equation containing only one unknown quantity.

$2^\circ$. Or, let the value of the unknown quantity, which you would first exterminate, be found in that equation where in it is the least involved, considering all the other quantities as known, and \dots\ be substituted \dots\ in the other equations; and with the new equations thus arising repeat the operation, till you have only one unknown quantity in one equation.

$3^\circ$. Or lastly, let the given equations be multiply'd or divided by such numbers or quantities, whether known or unknown, that the term which involves the highest power of the unknown quantity to be exterminated, may be the same in each equation; and then, by adding, or subtracting the equations \dots\ that term will vanish, and a new equation emerge, wherein the number of dimensions (if not the number of unknown quantities) will be diminished. 
\same --- \citet [63--64] {Simpson1755-HM}
\end {quotation}
Simpson followed Newton in not limiting the discussion to linear equations. His hesitancy about the outcome of case $3^\circ$ for polynomial equations suggests he may have been the originator of the addition and/or subtraction lesson. He includes twelve examples, of which eight are linear and solved by adding or subtracting pairs of equations to remove a variable common to both.

The innovations of Newton and Rolle are underscored by the failure of some authors to formulate an explicit rule for solving simultaneous equations. For example, a contemporary of Hammond and Simpson, the sightless Lucasian professor Nicholas \citet [164--166] {Saunderson1761-HM}, solved the three-equation problem of Peletier and Cardano using Gaussian elimination by the method of addition and/or subtraction. Saunderson explained that equations from which variables were removed should be grouped into \textit {ranks\/}, which were equivalent to Rolle's \textit {classes\/}, but Saunderson did not develop this terminology into a general prescription. 

Leonhard Euler, in his algebra textbook which was much admired for its concise style, provided a clear example of a missing rule. The text was originally written in German but appeared first in Russian in 1768 and eventually in many other languages.\footnote {\citet [3--4, 6] {Heefer2005} enumerates the editions of Euler's textbook. He traces most of Euler's exercises to an algebra textbook by Christoff Rudolff of 1525 which was reprinted by Michael Stifel in 1553.} Euler began with the compelling testimonial that the book was dictated for the instruction of his secretary, who mastered the subject from the text without additional instruction \citep [xxiii] {Euler1822-HM}. He included a chapter specifically for simultaneous linear equations. To find the values for two unknowns in two equations, he repeated the equality-of-values method.
\begin {quotation}
The most natural method of proceeding \dots\ is, to determine, from both equations, the value of one of the unknown quantities, as for example $x$, and to consider the equality of these two values; for then we shall have an equation, in which the
unknown quantity $y$ will be found by itself\dots. Then, knowing $y$, we shall only have to substitute its value in one of the quantities that express $x$.
\next --- \citet [part 2, sec.\ 1, chap.\ 4, art.\ 45] {Euler1771-HM} translated in \citet [206] {Euler1822-HM} 
\end {quotation}
Euler continued with equality-of-values for three equations. However, he cautioned against adopting a rote approach, and therefore did not state a general algorithm. 
\begin {quotation}
If there were more than three unknown quantities to determine, and as many equations to resolve, we should proceed in the same manner; but the calculation would often prove very tedious. It is proper, therefore, to remark, that, in each particular case, means may always be discovered of greatly facilitating the solution.
\next --- \citet [part 2, sec.\ 1, chap.\ 4, art.\ 53] {Euler1771-HM} translated in \citet [211] {Euler1822-HM} 
\end {quotation}



Euler's ambivalence toward Gaussian elimination became a dichotomy between research and instruction in the work of \'Etienne B\'ezout. The preface to B\'ezout's 1779 masterpiece on the theory of equations reviewed efforts to solve simultaneous polynomial equations. B\'ezout noted that Euler and Gabriel Cramer had studied how to ``reach the final equation'' in the reduction process that Newton and Rolle, both unmentioned, first made explicit. He then announced a way to ``improve the elimination method for first-order equations'' that takes advantage of zero coefficients \citep [xxiii and 138--144] {Bezout2006-HM}. Although B\'ezout referred to this process as elimination, in modern terminology he evaluated determinant solution formulas by a novel method for sparse systems, those with many zero coefficients.


B\'ezout did not let the determinants in his research intrude on his successful textbooks which he prepared for military academies. In these he taught Gaussian elimination. A posthumous and final version of his lessons began with ``\'egalez ces deux valuers'' to reduce two equations to one, and then continued to larger systems \citep [53--67, art.\ 65--73] {Bezout1788-HM}. He solved three equations by equality-of-values, and unlike Euler, B\'ezout did not hesitate to recommend the method for a greater number of equations. 
\begin {quotation}
\noindent \textbf {Equations of the first degree, with three or more unknowns}

\medskip
\noindent 71. Once what we have said is well understood, it is easy to see how we should behave when the number of unknowns and equations is greater.

We always assume that we have as many equations as unknowns. If there are three, \textit {we take in each the value of a single unknown, as if everything else were known. We then will match the first value to the second and first to third, or we will match the first to the second, and the second to third. We will, by this method, obtain two equations with two unknowns only, and which are treated by the preceding rule (66).\/}

If all the unknowns are not in every equation, the calculation would be easier, but it would always be done in a similar manner.

72. We see then that if there were a greater number of equations, the general rule would be: \textit {Take in each equation, the value of one unknown; equate one of these values to the others, and you have one equation and one unknown less. Treat these new equations as you have done for the first, and you will again have one equation and one unknown less. Continue until you finally have one equation that has no more than one unknown.\/} 
\same --- \citet [57 and 59, original emphasis] {Bezout1788-HM}\footnote {\textbf {Des \'Equations du premier degr\'e, \`a trois \& \'a un plus grand nombre d'inconnues.} \par 71. Ce que nous venons de dire \'etant une fois bien con\c{c}u, il est facile de voir comment on doit se conduire, lorsque le nombre des inconnues \& des \'equations est plus consid\'erable.\par Nous supposerons toujours qu'on ait autant d'\'equations que d'inconnues. Si l'on en a trois, \textit {on prenda dans chacune la valeur d'une m\^eme inconnue, comme si tout le rest \'etoit connu. On \'egalera ensuite la premi\`ere valeur \`a la seconde, \& la premi\`ere \`a la troisi\`eme; ou bien l'on \'egalera la premi\'ere \`a la seconde, \& la seconde \`a la troisi\`eme. On aura, par ce pro\'ced\'e, deux \'equations \`a deux inconnues seulement, \& on les traitera par la r\`egle pr\'ec\'edene (66).\/}\par Si toutes les inconnues n'entroient pas \`a la sois dans chaque \'equation, le calcul seroit plus simple, mais se seroit toujours d'une mani\`ere analogue.\par 72. On voit par-l\`a que s'il y avoit un plus grand nombre d'\'equations, la r\`egle g\'en\'erale seroit\dots. \textit {Prenez, dans chaque \'equation, la valeur d'une m\^eme inconnue; \'egalez l'une de ces valeurs \`a chacune des autres, \& vous aurez une \'equation \& une inconnue de moins. Traitez ces nouvelles \'equations comme vous venez de saire pour les premi\`eres, \& vous aurez encore une \'equation \& une inconnue de moins. Continuez ainsi jusqu'\`a ce qu'enfin vous parveniez \`a n'avoir plus qu'une inconnue.\/}}
\end {quotation}

As the \xviii\ century ended, Sylvestre Lacroix wrote a textbook which demonstrated that the time when authors could neglect Gaussian elimination was passing. \citet [1--2] {Domingues2008} characterizes Lacroix as a minor mathematician but an astute textbook author who sought to present the best approaches in an original, uniform style.\footnote {\citet [471] {MonthlyReview1801-HM} recommended Lacroix's \textit {Elemens d'alg\`ebre\/} for ``a great variety of valuable methods; and the excellence of the selection forms, in our opinion, the chief merit \dots.''} Lacroix's textbook included lessons on simultaneous linear equations that underwent considerable revision between the $\second$ and $\fifth$ editions of 1800 and 1804. In the $\second$ edition Lacroix discussed one and two unknowns in the same number of equations, and then passed to a derivation of explicit, determinant-like formulas for the unknowns in systems of two, three, and even four equations \citep [79--104, art.\ 75--91] {Lacroix2nd1800-HM}. In the $\fifth$ edition, between these two rather different discussions, he inserted text similar to Hammond's for up to four simultaneous linear equations. In the section title Lacroix made it plain the method could be applied to any number of equations.
\begin {quotation}
\noindent \textbf {Of the resolution of any given number of Equations of the First Degree, containing an equal number of Unknown Quantities.}


\medskip
\noindent 78. \dots\ if these unknown quantities are only of the first degree, [then] according to the method adopted in the preceding articles, \textit {we take in one of the equations the value of one of the unknown quantities, as if all the rest were known, and substitute this value in all the other equations, which will then contain only the other unknown quantities\/}. 

This operation, by which we exterminate one of the unknown quantities, is called \textit {elimination\/}. In this way, if we have three equations with three unknown quantities, we deduce from them two equations with two unknown quantities, which are to be treated as above; and having obtained the values of the two last unknown quantities, we substitute them in the expression for the value of the first unknown quantity. 

If we have four equations with four unknown quantities, we deduce from them, in the first place, three equations with three unknown quantities, which are to be treated in the manner just described; having found the values of the three unknown quantities, we substitute them in the expression for the value of the first, and so on.
\next --- \citet [114, art.\ 78, original emphasis] {Lacroix5th1804-HM}
\next translated in \citet [84, art.\ 78] {Lacroix1818-HM}
\end {quotation}
The American translator, John Farrar, inserted Newton's ``exterminate'' into the English text above, whereas the original is as follows.
\begin {quotation}
Cette op\'eration, par laquelle on chasse une des inconnues, se nomme \textit {\'elimination\/}.
\next --- \citet [114, art.\ 78, original emphasis] {Lacroix5th1804-HM}
\end {quotation}
With this passage Lacroix initiated the characterization of Gaussian elimination as \textit {elimination\/}. Dozens of \xix\ century American textbooks repeated ``[this] is called elimination'' once Farrer's 1818 translation of Lacroix appeared.\footnote {A search in Google books for the phrase ``is called elimination'' finds many examples in the period 1820-1900.} 

Following Lacroix, Louis Pierre Marie Bourdon was the next author of an influential algebra text that became grist for American translators. The $\fifth$ edition of \citet [79--80] {Bourdon1828-HM} identified three ways to eliminate a variable from two equations: Newton's equality-of-values or comparison, Newton and Rolle's substitution, and Simpson's addition and/or subtraction, the latter being the most highly recommended. Bourdon taught that simultaneous equations should be solved by using elimination to replace $m$ equations with $m-1$ equations, then $m-2$, and so on. Bourdon's text was translated by Edward \citet {Ross1831-HM} and abridged by Charles \citet {Davies1835-HM} to become \textit {Davies' Bourdon\/} which had ``extensive circulation in all parts of the United States'' \citep [372] {Cajori1890-HM}.

\section {Gaussian Elimination at the Time of Gauss}
\label {sec:part3}

\subsection {``Without any desire to do things which are useful''}
\label {sec:Neumann}

As the preceding survey of European texts shows, Newton brought order to the treatment of simultaneous equations through his presentation of a method that evolved into schoolbook Gaussian elimination. That is not to say simultaneous linear equations were more than pedagogic exercises through the period considered. Neither the scientist, Newton, nor the banker, Hammond, offered any compelling examples from their respective domains of expertise. 

None of the secondary literature argues that solving algebraic equations was needed, and several authors intimate to the contrary. \citet [71--72] {Neugebauer1969} points out that ancient economies required only arithmetic to function. The exercises in chapter 8 of the \textit {Nine Chapters\/} are plainly contrived. Although mathematics lessons can supply valuable information about daily life in the past \citep [416] {Libbrecht1973}, \citet [25] {Katz1997b} sees an ``artificial'' quality to algebra problems in Babylonian, Greek, Arabic, and European Renaissance texts. \citet {Hogendijk1994} explains that Islamic civilization needed arithmetic for commerce and advanced mathematics for astronomy, but the mathematical studies that in hindsight were the most sophisticated, such as algebra, were undertaken for their own sake. 

``By and large it is uniformly true in mathematics that there is a time lapse between a mathematical discovery and the moment when it is useful'' John \citet {vonNeumann1954-HM} opined in his measured prose, and in the meantime,``the whole system seems to function without any direction, without any reference to usefulness, and without any desire to do things which are useful''\footnote {Von Neumann was not an apologist for mathematics in the sense of G.\ H.\ Hardy. \citet [89] {Dieudonne1981} called von Neumann ``the last of the great mathematicians'' who understood both mathematics and its uses.} It is not often that a mathematical discovery is immediately useful, as was the case with the method of least squares. And it was the method of least squares which finally created a recurring need to solve simultaneous linear equations. When a use for schoolbook Gaussian elimination finally arose, Carl Friedrich Gauss invented something better for the problem at hand.

\subsection {``M\'ethode des moindres quarr\'es''}
\label {sec:squares}

The genesis of the method of least squares lay in a scientific question which in the \xviii\ century had not been completely resolved: how to make accurate predictions from measurements. In the middle of the century, differences between astronomical observations and orbital formulas, derived from Newton's principles, at times had even cast doubt on the inverse square law of gravitation.\footnote {See \citet [29] {Gillispie1997-HM}, \citet [142] {Goldstine1977}, and \citet [30] {Stigler1986} for elaboration.} Both Euler and Pierre-Simon Laplace had speculated that the law of gravity might need modification for astronomical distances. \citet [17, 28] {Stigler1986} explains that fitting the orbital formulas to more observations than were absolutely necessary had been inconceivable. Euler ``distrusted the combination of equations, taking the mathematician's view that errors actually increase with aggregation rather than taking the statistician's view that random errors tend to cancel one another'' \citep [28] {Stigler1986}. A new paradigm was adopted once Tobias Mayer successfully applied ad hoc fitting methods to predict the lunar orbit. Laplace in particular then derived fitted orbits that both vindicated Newton and established the mathematics that he, Laplace, employed, i.e. ``analysis.''\footnote {\citet [sect.\ 2] {Hawkins1977b} summarizes the status of analysis at the time of Lagrange and Laplace.} These successes raised the question, what was the best fitting method to use.\footnote {\citet {Farebrother1999} surveys the several fitting methods available at the end of the \xviii\ century including those used by Mayer and Laplace.}

Two inventors are recognized for the method of least squares. Gauss claimed that he knew in 1794 or 1795\footnote {\citet [art.\ 186] {Gauss1809-HM} wrote 1795 but \citet [241] {Plackett1972} quotes Gauss recalling 1794. Gauss may not have remembered the exact date because \citet [469] {Dunnington2004-HM} reports that Gauss began his famous diary of mathematical discoveries only in 1796.} ``the principle that the sum of the squares of the differences between the observed and computed quantities must be a minimum.''  Thus when the dwarf planet Ceres was sighted and lost in 1801, he quickly found its orbit by procedures that included least squares methods. Gauss was reticent about the calculations in the hasty announcement of his results, although he did send an explanation in 1802 to his friend Wilhelm Olbers.\footnote {Heinrich Wilhelm Matth\"aus Olbers was a physician and an accomplished amateur astronomer.} The explanation was inexplicably returned to Gauss only in 1805, yet in the interim Gauss had been publishing orbits for Ceres and other celestial bodies \citep [53 and 420--421] {Dunnington2004-HM}.

Meanwhile, in an appendix to a long paper, Adrien-Marie Legendre posed the general problem of finding the most accurate parameterization furnished by a given set of observations.\footnote {\citet [13 and 15] {Stigler1986} recommends Legendre's appendix as among the most elegant introductions of a significant mathematical concept.} \citet [72] {Legendre1805-HM} observed that the problem often involves many systems of equations each of the form,
\begin {displaymath}
E = a + bx + cy + fx + \mbox {etc.,}
\end {displaymath}
(his notation) where $a$, $b$, $c$, $f$, \dots\ are numbers that vary among the equations, and $x$, $y$, $z$, \dots\ are parameters common to all the equations. As in linear models today, the numbers in each equation represent the data of one observation; the unknowns are the model parameters. Legendre viewed the problem as finding values for the parameters to make $E$ small. If the number of equations exceeds that of the unknowns (so the equations cannot be solved exactly), then he suggested minimizing the sum of the squares of the $E$'s. He called this overall process the \textit {m\'ethode des moindres quarr\'es\/} (modern \textit {carr\'es\/}). The solution was found by differentiating the sum of squares to derive the ``equations of the minimum'' (modern normal equations). Since these simultaneous equations were linear and equinumerous with the variables, \citet [73] {Legendre1805-HM} said they could be solved by ``ordinary methods'' (par les m\'ethodes ordinaires).

Over a decade passed between Gauss's own discovery and its publication. Although Gauss intended to publish immediately after Legendre, the manuscript was delayed by the Napoleonic wars. In \textit {Theoria Motus\/} (1809), Gauss finally explained his process for orbital calculations, and he returned to the conceptual problem of fifty years earlier: how to justify values calculated from erroneous data. Lacking justification, he intimated that the only reason to minimize squares was convenience \citep [220--221, art.\ 186] {Gauss1809-HM}. Rather than begin by minimizing the discrepancy in the equations as Legendre had done, instead Gauss formulated ``the expectation or probability that all these values will result together from observation'' \citep [210, art.\ 175] {Gauss1809-HM}. Assuming the errors in the unknowns follow a Gaussian or normal distribution (later terminology), Gauss showed that it was maximizing the expectation which implied that the sum of squares should be minimized. He then echoed Legendre's sentiment about the resulting calculation to find the parameters. 
\begin {quotation}
\noindent We have, therefore, as many linear equations as there are unknown quantities to be determined, from which the values of the latter will be obtained by common elimination. (\dots\ per eliminationem vulgarem elicientur.) 
\next --- \citet [214, art.\ 180] {Gauss1809-HM}, this translation \citet [261] {Gauss1857-HM}
\end {quotation}


There then followed, in the published work of Gauss and Laplace, a scientific dialogue that explored the nature of probability and estimation.\footnote {Three short monographs with historical emphasis but each on slightly different aspects of this work --- statistical fitting procedures, parametric statistical inference, and classical analysis of variance --- have recently been written by \citet {Farebrother1999}, \citet {Hald2007}, and \citet {Clarke2008}, respectively. Histories with a wider scope include \citet [chap.\ 25] {Gillispie1997-HM}, \citet [chaps.\ 4.10] {Goldstine1977}, and \citet [chap.\ 4] {Stigler1986}.} During this period, \citet {Gauss1823-HM} gave a second and unqualified justification for the method of least squares, now known as the Gauss-Markov theorem for the minimum variance linear unbiased estimator. \citet [98, 105--109]{Hald2007} suggests that few contemporary readers, if any, understood all that Gauss and Laplace wrote. However, he considers that Gauss's first proof for the method of least squares, coupled with the sufficiency in many instances of the assumption of normally distributed errors, allowed the likes of \citet {Hagen1867-HM}, \citet {Chauvenet1868-HM}, and \citet {Merriman1884-HM} to maintain and extend a statistically respectable methodology of estimation from the time of Gauss through to the development of modern statistics. 


At the beginning of the \xx\ century a mathematics professor at the Massachusetts Institute of Technology could proudly announce that ``scientific investigations of all kinds'' relied on a mature computational technology called ``The Adjustment of Observations'' or ``The Method of Least Squares'' \citep [1 and 17] {Bartlett1900-HM}.\footnote {\citet [v--vi] {Bartlett1900-HM} lists English, French, and German textbooks from the end of the \xix\ century. See \citet {Ghilani2006} for a treatment from the beginning of the \xxi\ century.} The subject divided into two cases. Note that matrices were not used in these formulations until the mid \xx\ century. 
\newcounter{case}
\begin {list} {case \arabic{case}.} {\usecounter {case} \setlength {\leftmargin} {3.5em} \setlength {\labelwidth} {3.1em} \setlength {\itemindent} {0em}}
\item\label {case:one} The ``adjustment of indirect observations'' had Legendre's original, overdetermined equations that today would be stated as $\min_x \| b - A x \|_2$, where $A$ is an $m \times n$ matrix, $m > n$, and $b$ and $x$ are compatibly sized column vectors. These problems were solved by reducing them to $A^t A x = c$ with $c = A^t b$ where ${}^t$ is transposition \citep [17, arts.\ 23--25] {Bartlett1900-HM} \citep [93 art.\ 74, 106 art.\ 84] {Wright1906-HM}.

\item\label {case:two} The ``adjustment of conditioned observations'' was formulated by \citet {Gauss1826-HM} to find minimum 2-norm solutions of underdetermined equations, $\min_{Ax = b} \| x \|_2$, where $m < n$. He reduced these problems to $A A^t u = b$ with $x = A^t u$. The latter equations for $x$ were called the \textit {correlate equations\/} \citep [111, arts.\ 116--118] {Bartlett1900-HM} \citep [152 art.\ 119, 243 art.\ 179] {Wright1906-HM}.
\end {list} 
Each row of the over-determined (case \ref {case:one}) or under-determined (case \ref {case:two}) equations, $A x = b$, was called a \textit {condition\/}. The reduced forms of both problems were called \textit {normal equations\/} and were solved by various forms of elimination.\footnote {\citet [415--420] {Stigler1999-HM} reports that \citet [84] {Gauss1822a-HM} appears to have used \textit {Normalgleichungen\/} first but also just once and offhand. Gauss did not explain what he meant by the name. He could have simply meant ``ordinary'' equations.}

\subsection {Gauss's Formulation of the Method of Least Squares}
\label {sec:Gauss}

Gauss performed meticulous calculations almost as a leisure activity throughout his life \citep [227-228] {Dunnington2004-HM}. He thought calculations so important that he included them in his papers striving to make them brief and intuitively clear. Accordingly, his publication that mentioned ``common elimination'' was followed with a detailed explanation. Rather than solve the simultaneous linear equations ``of the minimum'' as Legendre had done, instead \citet {Gauss1810-HM} formulated his calculations in terms of quadratic forms.

A canonical expression for quadratic forms had already been introduced by \citet {Lagrange1759-HM}. In modern matrix notation, a quadratic form is given by $x^t A x$ for a column vector of variables $x$ and a symmetric matrix $A$. Lagrange substituted new variables for linear combinations of the original variables, which amounts to constructing a triangular matrix $U$ of substitution coefficients so that $A = U^t D U$ for a diagonal matrix $D$. The quadratic form thus became a weighted sum of squares, $(Ux)^t D (Ux)$, which Lagrange used to ascertain local extrema. Such was the approach that Gauss took for least squares.\footnote {Gauss may have learned the quadratic construction from Lagrange. At the time when Gauss first considered least squares in 1795, he borrowed from the G\"ottingen library the journal volume that contains Lagrange's paper on extrema \citep [398] {Dunnington2004-HM}. Since Gauss did not need this volume for his major work from this period which was on number theory (in which he did cite papers from the journal but none from this particular volume \citep [art.\ 202, n.\ 9] {Gauss1801-HM}), it seems likely that he borrowed the volume to consult Lagrange's paper to learn about minimizing quadratic forms.} Since neither \citet {Lagrange1759-HM} nor \citet {Gauss1801-HM} dealt with quadratic forms of more than three variables,\footnote {\citet [83--85] {Hawkins1977a} summarizes the use of quadratic forms in \citet {Gauss1801-HM}.} Gauss had yet to systematically extend the construction of the canonical expression to arbitrarily many unknowns, which he finally did in his papers about least squares.\footnote {A difficulty in extending the construction was the subscriptless notation of the time. \citet [161, n.] {Farebrother1999}\ attributes double subscript notation to Cauchy in 1815 \citep [114] {Cauchy1905-HM}. Gauss first used numeral indices in 1826, but as superscripts \citep [83] {Gauss1826-HM}.}

\citet [22] {Gauss1810-HM} wrote the overdetermined equations as
\setlength {\arraycolsep} {0.0625em}
\newcommand {\p} {^{\,\prime}}
\newcommand {\pp} {^{\,\prime\prime}}
\newcommand {\ppp} {^{\,\prime\prime\prime}}
\begin {displaymath}
\setlength {\arraycolsep} {0.0625em}
\begin {array} {l c l c c l c c l c c l c c l}
n& \plus& a& p& \plus& b& q& \plus& c& r& \plus& d& s& \plus& \dots
\equals w
\\
n\p& \plus& a\p& p& \plus& b\p& q& \plus& c\p& r& \plus& d\p& s& \plus&
\dots \equals w\p
\\
n\pp& \plus& a\pp& p& \plus& b\pp& q& \plus& c\pp& r& \plus& d\pp& s&
\plus& \dots \equals w\pp
\\
\multicolumn {2} {l} {\dots}
\end {array}
\end {displaymath}
(original notation, but with ``\dots''\ replacing ``etc.''). The symbols $n$, $a$, $b$, $c$, \dots, with or without primes, are numbers. The purpose is to find values for the variables $p$, $q$, $r$, \dots\ to minimize
\begin {displaymath}
\Omega \equals w w \plus w\p w\p \plus w\pp w\pp \plus \dots \, .
\end {displaymath}
Gauss introduced a bracket notation
\begin {displaymath}
[xy] \equals x y \plus x\p y\p \plus x\pp y\pp \plus \dots \, ,
\end {displaymath}
where the letter $x$ either is $y$ or lexicographically precedes $y$.\footnote {That is, ``$[xy]$'' is a composite symbol that Gauss meant to be used only when the letter ``$x$'' precedes or matches the letter ``$y$'' in alphabetic order.} Gauss did not have a name for his bracket symbols; they would be known as ``auxiliaries'' because they represented intermediate values that were needed to solve the least squares problem by Gauss's approach. This notation expressed the normal equations (name not yet introduced) as follows.
\begin {equation}
\label {eqn:normal}
\begin {array} {c c c c c c c c c c c c c}
[an]& \plus& [aa]& p& \plus& [ab]& q& \plus& [ac]& r& \plus& [ad]& s
\plus \dots \equals 0
\\{}
[bn]& \plus& [ab]& p& \plus& [bb]& q& \plus& [bc]& r& \plus& [bd]& s
\plus \dots \equals 0
\\{}
[cn]& \plus& [ac]& p& \plus& [bc]& q& \plus& [cc]& r& \plus& [cd]& s
\plus \dots \equals 0
\\
\dots
\end {array}
\end {equation}
As \citet [73] {Legendre1805-HM} and he himself had done before \citep [214] {Gauss1809-HM}, Gauss again remarked that these equations could be solved by elimination \citep [22] {Gauss1810-HM}, but he did not explicitly perform the calculation. Instead, he noted that the brackets give the coefficients of the variables in the sum of squares, a quadratic form. 
\begin {displaymath}
\begin {array} {c c c c c c c c c c c c c c c c c c}
\Omega& =& \multicolumn {16} {l} {[nn] \plus 2 \, [an] \, p \plus 2 \,
[bn] \, q \plus 2 \, [cn] \, r \plus 2 \, [dn] \, s \plus \dots}
\\
& \plus& [aa]& pp& \plus& 2& [ab]& pq& \plus& 2& [ac]& pr& \plus& 2&
[ad]& ps& \plus \dots
\\
& \plus& [bb]& qq& \plus& 2& [bc]& qr& \plus& 2& [bd]& qs& \plus&
\multicolumn {2} {l} {\dots}
\\
& \plus& [cc]& rr& \plus& 2& [cd]& rs& \plus& \multicolumn {2} {l}
{\dots}
\\
& \multicolumn {2} {l} {\dots}
\end {array}
\end {displaymath}
Gauss extended his bracket notation to
\begin {equation}
\label {eqn:brackets}
\begin {array} {r c l}
[xy, 1]& \equals& \displaystyle [xy] - {[ax] [ay] \over [aa]}
\\ \noalign {\medskip}
[xy, 2]& \equals& \displaystyle [xy,1] - {[bx,1] [by,1] \over [bb,1]}
\\ \noalign {\medskip}
[xy, 3]& \equals& \displaystyle [xy,2] - {[cx,2] [cy,2] \over [cc,2]}
\\ \noalign {\medskip}
\multicolumn {2} {c} {\dots}
\end {array}
\end {equation}
and so on. The values of these formulas are the coefficients remaining after successive combinations of variables have been grouped into perfect squares. The first of these combinations of variables, $A$,
\begin {equation}
\label {eqn:Gauss}
\begin {array} {r c l}
A& \equals& [an] \plus [aa] p \plus [ab] q \plus [ac] r \plus [ad] s
\plus \dots
\\ \noalign {\smallskip}
B& \equals& [bn,1] \plus [bb,1] q \plus [bc,1] r \plus [bd,1] s \plus
\dots
\\ \noalign {\smallskip}
C& \equals& [cn,2] \plus [cc,2] r \plus [cd,2] s \plus \dots
\\ \noalign {\smallskip}
\multicolumn {3} {l} {\dots,}
\end {array}
\end {equation}
simplifies the quadratic form by grouping all occurrences of the variable $p$ into a perfect square, thus removing $p$ from the balance of $\Omega$:
\begin {displaymath}
\begin {array} {c c c c c c c c c c c c c c c c c c}
\displaystyle \Omega - \vrule depth0pt height3.25ex width0pt\smash {A^2 \over [aa]}& =& \multicolumn {16} {l} {[nn,1] \plus 2 \,
[bn,1] \, q \plus 2 \, [cn,1] \, r \plus 2 \, [dn,1] \, s \plus \dots}
\\
& \plus& [bb,1]& qq& \plus& 2& [bc,1]& qr& \plus& 2& [bd,1]& qs& \plus&
\multicolumn {2} {l} {\dots}
\\
& \plus& [cc,1]& rr& \plus& 2& [cd,1]& rs& \plus& \multicolumn {2} {l}
{\dots}
\\
& \multicolumn {2} {l} {\dots.}
\end {array}
\end {displaymath}
If this process is repeated with $B$, $C$, \dots, then eventually,
\begin {displaymath}
\Omega \minus {A^2 \over [aa]} \minus {B^2 \over [bb,1]} \minus {C^2
\over [cc,2]} \minus \, \dots \, \equals [nn,\mu] \, ,
\end {displaymath}
where $\mu$ is the quantity of variables. Each combination $A$, $B$, $C$, \dots\ has one fewer unknown than the preceding combination. Thus $A = 0$, $B = 0$, $C = 0$, \dots\ can be solved in reverse order to obtain the values for $p$, $q$, $r$, \dots, also in reverse order, at which $\Omega$ attains its minimum, $[nn,\mu]$. In later theoretical discussions of least squares methods, e.g.\ \citet [art.\ 13] {Gauss1826-HM}, Gauss always referred back to his 1810 paper for details of the calculations for transforming quadratic forms. 

Gauss's contributions to the method of least squares became known almost immediately. By 1819 even a gymnasium prospectus, \citet {Paucker1819}, cited \citet {Legendre1805-HM} and \citet {Gauss1809-HM} for formulating the statistical inference problem (modern terminology), although Gauss was not cited for any particular method of calculating the inferences. Gauss's process that neatly tied together (what we now call) linear algebra, optimization theory, and his probabilistic justification for the inferences, seems to have been adopted as an algorithm slowly and by geodesists not by mathematicians. 

\subsection {Geodesy for Cartography}
\label {sec:cartography}

The original motivation for the computational developments in least squares was to cope with observational errors in two major scientific activities of the \xix\ century, and Gauss was involved in both. One was astronomy which was pursued for navigation and timekeeping as well as for its intrinsic interest. \citet {Gauss1809-HM} described the calculations that he had invented to derive orbital formulas from a few observations. He reported the results of his own calculations before and after 1809 in many papers that constituted the bulk of his early publications.\footnote {\citet {Nievergelt2001} explains the orbital calculations that were used in \xix\ century astronomy, while \citet {Grier2005} describes the institutional history of computing groups in national observatories where the calculations were made.} 

The other scientific field was geodetic research for cartography, an activity continually and directly sponsored by governments. Indeed, the first scientific agency of the United States was the Coast Survey Office founded in 1807.\footnote {\citet {Cajori1929} traces the colorful early years of the Coast Survey. The enthusiasm for scientific research at the Coast Survey and the overlapping responsibilities of the younger Geological Survey were the subject of vehement debate in the government \citep {Congress1886-HM}.} Gauss became prominent in geodesy through the many papers he wrote during his protracted survey of Hanover. This small German state, roughly coincident with modern Lower Saxony, was Gauss's home for most of his life. 

Friedrich Bessel warned Gauss that the toil of the survey would detract from his research \citep [120] {Dunnington2004-HM}. However, the time was well spent. Although Legendre invented the method of least squares to solve a geodetic problem, an English survey officer explained that it was Gauss from whom geodesists adopted the approach. 
\begin {quotation}
\noindent If in effecting a [survey] triangulation one observed only just so many angles as were absolutely necessary to fix all the points, there would be no difficulty in calculating [the locations]; only one result would be arrived at. But it is the invariable custom to observe more angles than are absolutely needed, and it is these supernumerary angles which give rise to complex calculations. Until the time of Gauss and Bessel computers had simply used their judgement as they best could as to how to employ and utilize the supernumerary angles; the principal of least squares showed that a system of corrections ought to be applied, one to each observed bearing or angle, such that subject to the condition of harmonizing the whole work, the sum of their squares should be an absolute minimum. The first grand development of this principle is contained in this work of Bessel's.
\same --- \citet [26--27] {Clarke1880-HM}
\end {quotation} 
Clarke refers to a Prussian triangulation to which Bessel applied the panoply of Gauss's techniques.\footnote {\citet [21] {Clarke1880-HM} refers to \citet [130] {Bessel1838-HM} which is discussed by \citet [11] {Torge2001}.} Bessel's endorsement of Gauss in a survey for a major government drew the attention of other geodesists.

\begin {figure} 
\centering
\includegraphics [scale=1] {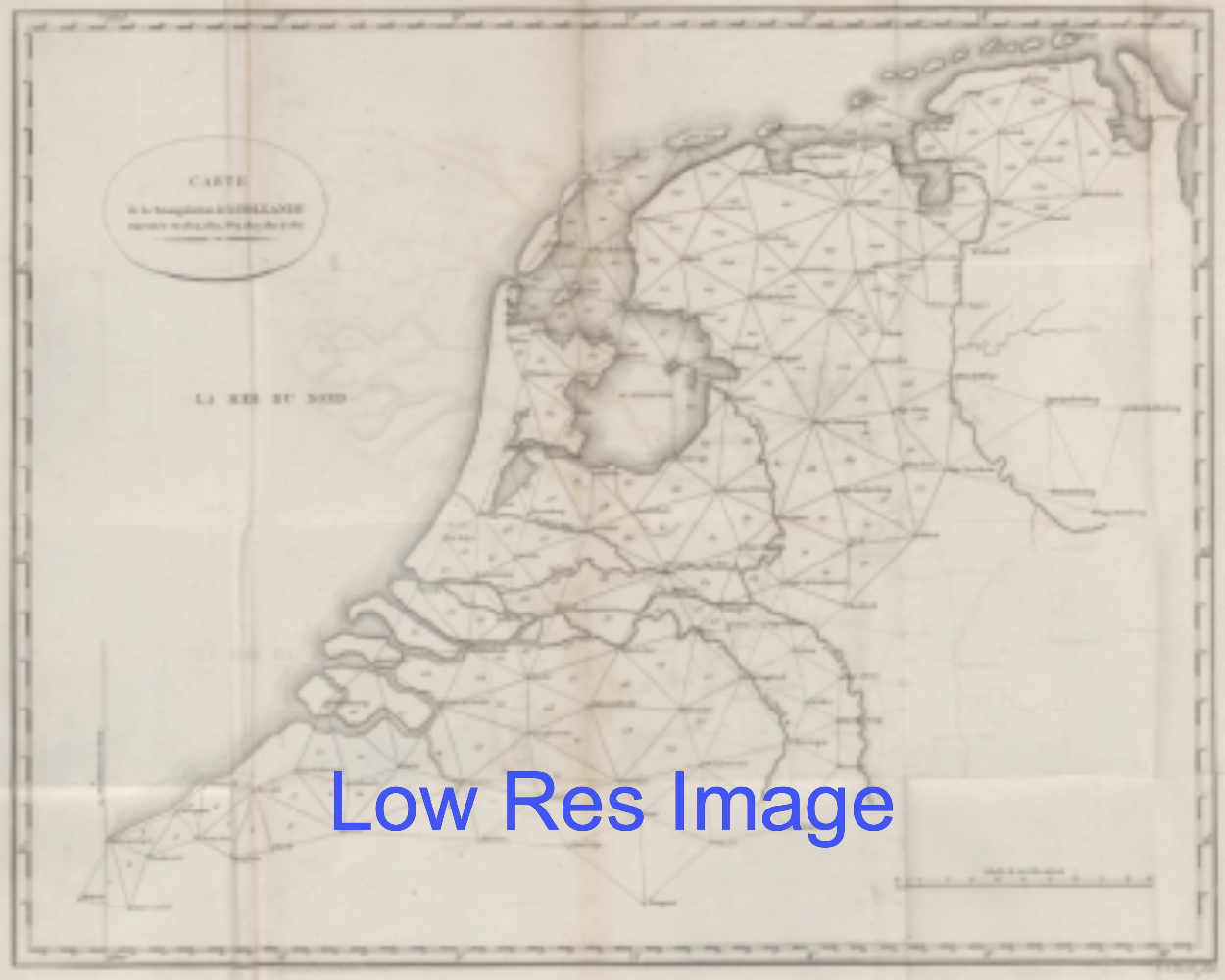}%
\caption {\citet [86, art.\ 23] {Gauss1826-HM} used a small portion of a triangulation of Holland to illustrate survey adjustments. He took the data from \citep {Krayenhoff1827-HM}, a book in his personal library. The map seen here is from the earlier edition \citet {Krayenhoff1813}. Courtesy of the Bancroft Library, University of California, Berkeley.}
\label {fig:holland}
\end {figure}
\begin {figure} 
\vspace*{1ex}
\centering
\includegraphics [scale=1] {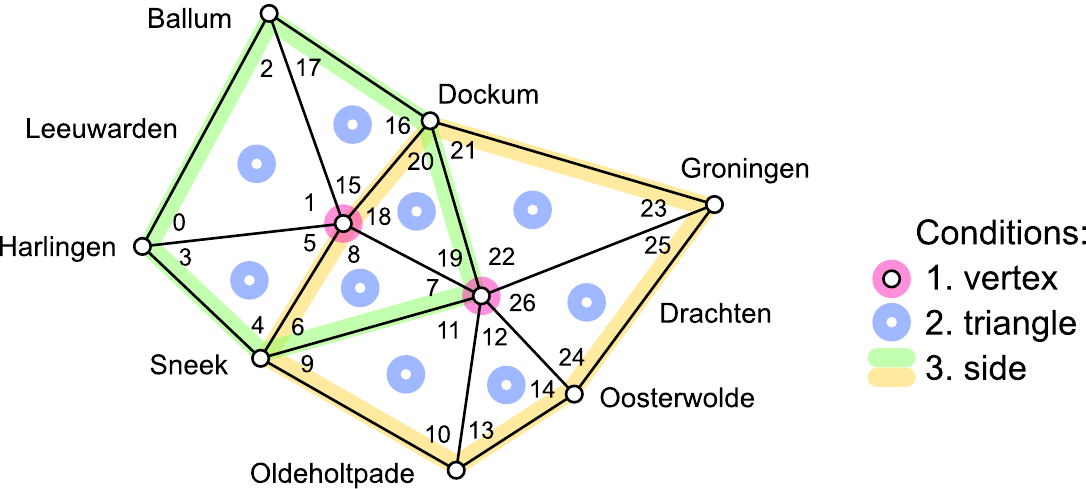}%
\caption {The portion of Figure \ref {fig:holland} that Gauss readjusted. He had to find adjustments for 27 angles, but he had only 13 conditions for them to satisfy, consisting of 2 vertex conditions, 9 triangle conditions, and 2 side conditions. \citet [489] {Jordan1895-HM} reproduces a similar figure.}
\label {fig:krayenhoff}
\end {figure}

As Clarke explained, a critical step in making maps was to reconcile the measurement errors in the raw data gathered by surveyors. At each vertex surveyors measured either the included angles of triangles or the azimuthal directions to nearby vertices; in the latter case the angles were differences of the directions.\footnote {Each survey prescribed strict uniformity in taking measurements. The measuring instrument was and is called a \textit {theodolite\/}. The book by \citet {Clarke1880-HM}, the research paper by \citet {Nell1881-HM}, the encyclopedia by \citet {Jordan1895-HM}, and the textbook by \citet {Wright1906-HM} have further explanations and more examples from the past.} \citet [art.\ 22] {Gauss1826-HM} stipulated that the angles satisfy three kinds of conditions: (1) the sum of angles around an interior vertex equals $2 \pi$ radians, (2) the sum of angles in a triangle equals $\pi$ plus the spherical excess, and (3) \textit {side conditions\/} that chain together linearized sine laws for circuits of triangles with common edges. The unknowns are perturbations to the measurements intended to make the angles satisfy the conditions. Gauss formed a side condition from the triangles around each interior vertex, in which case (it is easy to see) even an ideal net consisting of $f$ nonoverlapping triangles and $v$ vertices has $3f$ angles, but only $3f-2v+4$ conditions. Thus, in general, the adjustment problems were under- not overdetermined. In his last major theoretical work on least squares, \citet {Gauss1826-HM} introduced the solution method given in section \ref {sec:squares} as case \ref {case:two}. He illustrated the method by readjusting a small part of the Dutch triangulation (Figures \ref {fig:holland} and \ref {fig:krayenhoff}). In comparison the British Isles triangulation was quite irregular \citep [II, plate xviii] {OrdnanceSurvey1858-HM}. Most surveys could not afford to follow a highly regular triangulation and also had missing measurements from inaccessible vertices (e.g.\ mountaintops) or blocked sight lines. 

As the \xix\ century progressed, the growing use of least squares methods created a recurring need to solve dauntingly large problems. Of necessity these large problems were broken into smaller subproblems. 
\begin {quotation}
\noindent In the principal triangulation of Great Britain and Ireland there are 218 stations, at 16 of which there are no observations, the number of observed bearings is 1554, and the number of equations of condition 920.\footnote {In \citet [158] {Stigler1986} the numbers 1554 and 920 are mistakenly switched.} The reduction of so large a number of observations in the manner we have been describing [i.e.\ least squares case \ref {case:two}] would have been quite impossible, and it was necessary to have recourse to methods of approximation\dots.

\dots\ the network covering the kingdom was divided into a number of blocks, each presenting a not unmanageable number of equations of condition. One of these being corrected or computed independently of the others, the corrections so obtained were substituted (as far as they entered) in the equations of condition of the next block, and the sum of the squares of the remaining equations in that figure made a minimum. The corrections thus obtained for the second block were substituted in the third and so on. \dots\ The number of blocks is 21, in 9 of them the number of equations of condition is not less than 50, and in one case the number is 77. These calculations --- all in duplicate --- were completed in two years and a half --- an average of eight computers being employed\dots.

In connection with so great a work successfully accomplished, it is but right to remark how much it was facilitated by the energy and talents of the chief computer, Mr.\ James O'Farrell.
\same --- \citet [237, 243] {Clarke1880-HM}
\end {quotation}
The systems of linear equations for the subproblems were still very large for hand computing. In these situations, professional computers found it useful to emulate Gauss's calculations of fifty years earlier. 

\section {Perfecting Gaussian Elimination for Professional Computers}
\label {sec:part4}

\subsection {Overview of Perfecting Calculations}

``Numerical mathematics stayed as Gauss left it until World War II,'' concluded \citet [287] {Goldstine1972} in his history of the subject. In the 120 years after \citet {Gauss1826-HM} there were at best a dozen noteworthy publications on solving simultaneous linear equations. The two main thrusts were simplifications for professional computers, which began with Gauss himself, and matrix interpretations. These developments overlap only slightly chronologically, and they did not influence each other until the middle of the \xx\ century. The former will be addressed here and the latter in the following section. Four authors, Gauss, Doolittle, Cholesky, and Crout, demonstrably influenced professional hand-computing practice for Gaussian elimination. 

\subsection {Gauss's Convenient Notation}
\label {sec:notation}

The triangulation of the great \citet {OrdnanceSurvey1858-HM} and the formulation of its adjustment \citep {Clarke1880-HM} are preserved in detail, but exactly how the calculations were performed is missing. Since the details of major calculations were not archived, whether and why Gauss's solution method was used --- and what it was considered to be --- must be inferred from sources such as William \citet {Chauvenet1868-HM}.\footnote {Chauvenet is remembered as a founder of the United States Naval Academy and by the Mathematical Association of America's Chauvenet Prize.} There appear to have been three advantages to Gauss's method of solving the normal equations (described in section \ref {sec:Gauss}). 

\begin {figure} [t]
\centering
\includegraphics [scale=1] {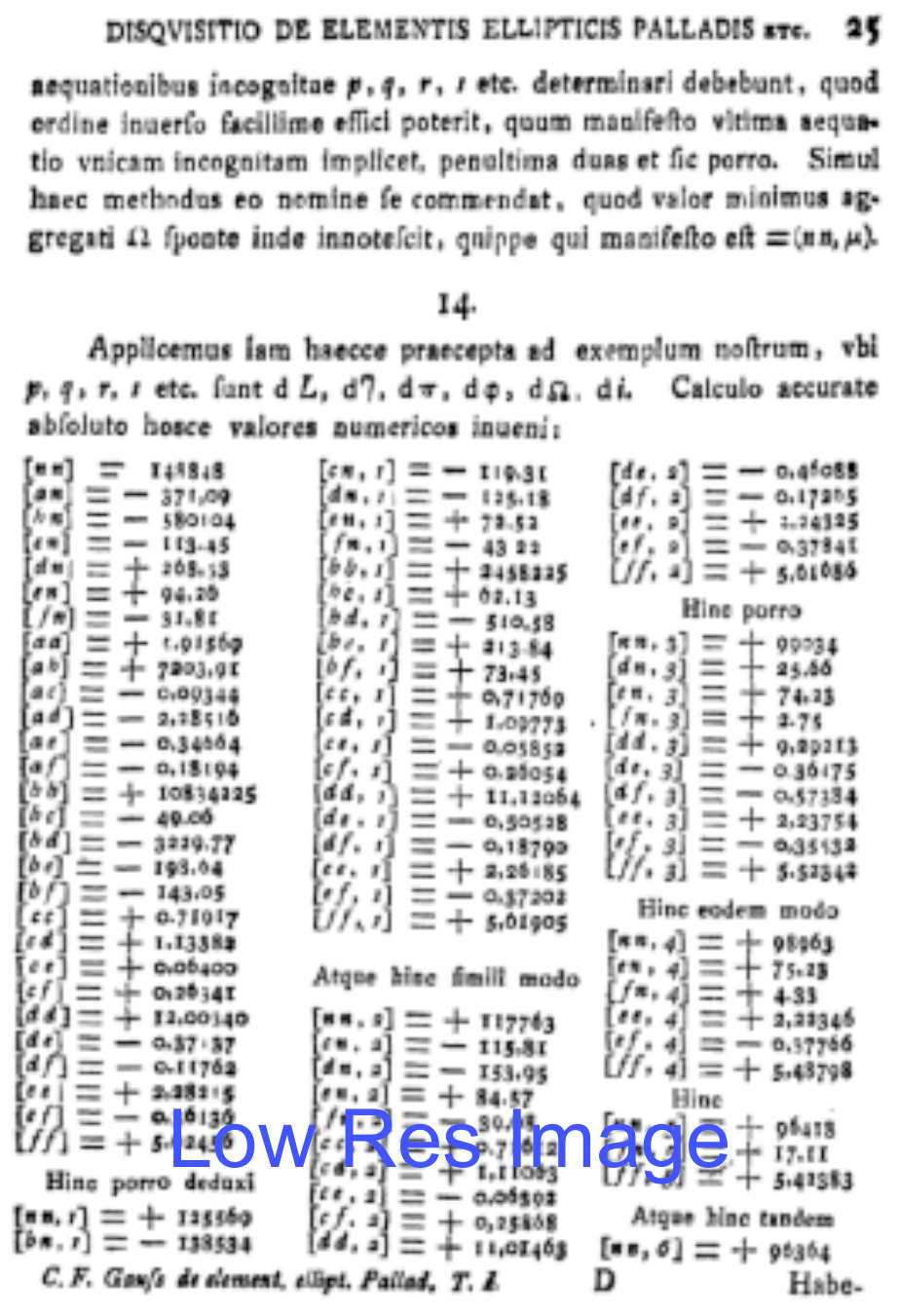}%
\caption {How Gauss performed Gaussian elimination. He wrote down the numbers in the order he computed them, using his bracket notation to identify the values \citep [25] {Gauss1810-HM}.}
\label {fig:pagefromgauss}
\end {figure}

First, the bracket notation conceptually separated the algebra from the arithmetic so the workflow could be addressed. Gauss relieved computers of the tedium of having to rewrite equations, and in so doing, he enabled them to consider how best to organize their work.
\begin {quotation}
\noindent By whatever method of elimination is performed, we shall necessarily arrive at the same final values of the unknown quantities; but when the number of equations is considerable, the method of substitution, with Gauss's convenient notation, is universally followed.
\same \mbox {--- \citet [514] {Chauvenet1868-HM}}
\end {quotation}
It may come as some surprise to learn that when Gauss performed Gaussian elimination, he simply listed all the numbers in the order he computed them, using his bracket notation to identify the values: $[cd] = 1.13382$, $[cd, 1] = 1.09773$, $[cd, 2] = 1.11063$, etc.\ (Figure \ref {fig:pagefromgauss}). \citet [138] {Dunnington2004-HM} echoed Bessel in regretting how much time Gauss spent calculating for his interminable survey projects: Gauss estimated that he used one million numbers! 

Second, in contrast to the method depicted in Figure \ref {fig:today} and prescribed in algebra textbooks, Gauss realized economies by avoiding duplicate calculations for symmetric equations.\footnote {The normal equations are symmetric in the sense that the same coefficient is attached to the $\jth$ variable in the $\kth$ equation, and vice versa.}

\begin {quotation}
\noindent By means of a peculiar notation proposed by Gauss, the elimination by substitution is carried on so as to preserve throughout the symmetry which exists in the normal equations. \same --- \citet [530] {Chauvenet1868-HM}
\end {quotation}
Chauvenet even foreshadowed modern practice for measuring the work of a calculation by counting how many bracket or auxiliary values had to be computed: 156 for 8 normal equations, etc. 

Third, Gauss included refinements such as estimates of precision, variances, and weights for the unknowns, all expressed in his bracket notation. The difficulty of making changes to the formulas, whose complexity was compounded by their relationship to Gauss's statistical theories, and the demonstrable benefits of Gauss's ideas, engendered a reluctance to alter his computational prescriptions. His efficient method for overdetermined problems by itself was conceptually difficult because the solution of the normal equations, $A A^t u = b$, was not the solution of the problem, $x = A^t u$. Gauss's bracket notation was still being taught a hundred years after his 1811 paper.\footnote {\citet {Bartlett1900-HM}, \citet {Johnson1905-HM}, and \citet {Wright1906-HM} offer instruction on the use of bracket notation. For continued emphasis on the advantages of preserving symmetry, see \citet [84--85] {Palmer1912-HM} and \citet [133, chap.\ 4] {Golub1996} .} At the beginning of the \xx\ century, \citet {Wright1906-HM} saw professional computers using only two methods to solve normal equations: either the brackets of Gauss or the tables of Doolittle.

\subsection {Doolittle, Legendary Computer}
\label {sec:Doolittle}

Myrick Hascall Doolittle was a computer who solved Gauss's normal equations to adjust triangulations at the United States Coast and Geodetic Survey for a period of 38 years from 1873 to 1911.\footnote {For further information on Doolittle, see the biographical sketch by \citet [78--79] {Grier2005}, and the primary biographical sources quoted at length by \citet {Farebrother1987}.} He is still remembered for his professional acumen as a human computer. Demonstrating the efficacy of his methods, \citet [117] {Doolittle1878-HM} notes that he solved $41$ normal equations in a week using paper and pencil. For perspective, \citet [676] {Fox1987b} reports that four mathematicians, Alan Turing among them, needed about two weeks to solve $18$ equations with ``desk computing equipment'' in 1946.\footnote {Using Gaussian elimination to find $n$ unknowns in $n$ equations requires about $2 n^3 / 3$ arithmetic operations. The approximately 3,888 operations for $n = 18$ compares well with the 4,800 minutes in ten working days. The question is thus not why were Fox et al.\ slower, but rather how was Doolittle so quick.} This does not, however, imply that Doolittle was a numerical savant. A colleague, Mr.\ J.\ G.\ Porter, duplicated the calculation in the same time to check for errors.

\newcommand {\disp} [1] {$\displaystyle {#1}$}
\newcommand {\dsd} {\vrule depth1.5ex height2.5ex width0pt}
\newcommand {\xss} {\vrule depth1.3ex height0ex width0pt}
\newcommand {\sds} {\vrule depth3.0ex height4.0ex width0pt}
\newcommand {\dds} {\vrule depth3.0ex height4.0ex width0pt}
\newcommand {\ddd} {\vrule depth3.0ex height4.0ex width0pt}
\newcommand {\sdd} {\vrule depth3.0ex height4.0ex width0pt}
\newcommand {\ssd} {\vrule depth1.5ex height2.5ex width0pt}
\newcommand {\sss} {\vrule depth1.5ex height2.5ex width0pt}
\renewcommand {\arraystretch} {1.0}

\begin {figure}[h!]
\begin {center}
\footnotesize
\scriptsize
{\small Table A}\\[1.0ex]
\begin {tabular} {r c c c c c c}
\xss& \parbox[b]{2.5em}{\centering recip\-rocal}& $w$& $x$& $y$& $z$& \parbox[b]{3.5em}{\centering absolute term}\\ \cline{2-7}
step 1.\ssd&& $[aa]$& $[ab]$& $[ac]$& $[ad]$& $[an]$\\
2.\sds& \disp {-1\over[aa]}& $w={}$& \disp {-{[ab]\over[aa]}}& \disp {-{[ac]\over[aa]}}& \disp {-{[ad]\over[aa]}}& \disp {-{[an]\over[aa]}}\\ \cline{2-7}
5.\dsd&&& $[bb,1]$& $[bc,1]$& $[bd,1]$& $[bn,1]$\\
6.\sds& \disp {-1\over[bb,1]}&& $x={}$& \disp {-{[bc,1]\over[bb,1]}}& \disp {-{[bd,1]\over[bb,1]}}& \disp {-{[bn,1]\over[bb,1]}}\\ \cline{2-7}
10.\dsd&&&& $[cc,2]$& $[cd,2]$& $[cn,2]$\\
11.\sds& \disp {-1\over[cc,2]}&&& $y={}$& \disp {-{[cd,2]\over[cc,2]}}& \disp {-{[cn,2]\over[cc,2]}}\\ \cline{2-7}
16.\dsd&&&&& $[dd,3]$& $[dn,3]$\\
17.\sds& \disp {-1\over[dd,3]}&&&& $z={}$& \disp {-{[dn,3]\over[dd,3]}}\\ \cline{2-7}
\end {tabular}\\[8ex]
{\small Table B}\\[1.0ex]
\begin {tabular} {r c c c c}
\xss& $x$& $y$& $z$& \parbox[b]{3.5em}{\centering absolute term}\\ \cline{2-5}
step 3.\ssd& $[bb]$& $[bc]$& $[bd]$& $[bn]$\\
4.\sds& \disp {-{[ab]\over[aa]}\times [ab]}& \disp {-{[ab]\over[aa]}\times[ac]}& \disp {-{[ab]\over[aa]}\times[ad]}& \disp {-{[ab]\over[aa]}\times [an]}\\ \cline{2-5}
7.\dsd&& $[cc]$& $[cd]$& $[cn]$\\
8.\sdd&& \disp {-{[ac]\over[aa]} \times [ac]}& \disp {-{[ac]\over[aa]} \times [ad]}& \disp {-{[ac]\over[aa]} \times [an]}\\
9.\dds&& \disp {-{[bc,1]\over[bb,1]} \times [bc,1]}& \disp {-{[bc,1]\over[bb,1]} \times [bd,1]}& \disp {-{[bc,1]\over[bb,1]} \times [bn,1]}\\ \cline{2-5}
12.\dsd&&& $[dd]$& $[dn]$\\
13.\sdd&&& \disp {-{[ad]\over[aa]} \times [ad]}& \disp {-{[ad]\over[aa]} \times [an]}\\
14.\ddd&&& \disp {-{[bd,1]\over[bb,1]} \times [bd,1]}& \disp {-{[bd,1]\over[bb,1]} \times [bn,1]}\\
15.\dds&&& \disp {-{[cd,2]\over[cc,2]} \times [cd,2]}& \disp {-{[cd,2]\over[cc,2]} \times [cn,2]}\\ \cline{2-5}
\end {tabular}
\end {center}
\caption {How \citet [115] {Doolittle1878-HM} performed Gaussian elimination, transcribed from his numerical example using Gauss's bracket notation to identify the quantities. The step numbers indicate the order of forming the rows. Each row in table B is a multiple, by a single number, of a partial row in table A. The rows 5, 10, and 16 in table A are sums of rows 3--4, 7--9. and 12--15 in table B, respectively.}
\label {fig:tablesAandB}
\end {figure}

Doolittle's speed stemmed from streamlining the work for hand computing. His one paper on the subject \citep {Doolittle1878-HM} contained a small numerical example that is restated here in Gauss's bracket notation to clarify the calculation. Corresponding to Gauss's equation (\ref{eqn:normal}), Doolittle's normal equations were as follows.
\begin {equation}
\label {eqn:DoolittleNormalEquations}
\begin {array} {c c c c c c c c c c c c c c}
0 \equals& [aa]& w& \plus& [ab]& x& \plus& [ac]& y& \plus& [ad]& z& \plus& [an]
\\{}
0 \equals& [ab]& w& \plus& [bb]& x& \plus& [bc]& y& \plus& [bd]& z& \plus& [bn]
\\{}
0 \equals& [ac]& w& \plus& [bc]& x& \plus& [cc]& y& \plus& [cd]& z& \plus& [cn]
\\{}
0 \equals& [ad]& w& \plus& [bd]& x& \plus& [cd]& y& \plus& [dd]& z& \plus& [dn]
\end {array}
\end {equation}
Doolittle expressly planned the calculation to derive Newton's and Rolle's substitution formulas for each variable, which correspond to rearrangements of Gauss's equations $A = 0$, $B = 0$, \dots. Doolittle called these formulas the ``explicit functions'' for the variables. He was able to colocate many of the numbers by completing the formula for a given variable before undertaking any calculations for the next. Among the practices that he introduced, is placing the numbers of the calculation in tables to identify them. Doolittle kept the coefficients of the substitution formulas in table A (Figure \ref {fig:tablesAandB}), while he used a second table, B, to record in columns the sums that give the values in table A. Exactly how Doolittle conducted the work may be lost. ``For the sake of perspicuity,'' he noted, ``I have here made some slight departures from actual practice'' \citep [117] {Doolittle1878-HM}.

The salient feature of Doolittle's tables is a reduction in the labor of division and multiplication. All divisions reduce to multiplications through reciprocals formed (one per variable) in the first column of table A. All multiplications have a single multiplier repeatedly applied to several multiplicands in a row of table A and recorded in another row of table B. For example, row 8 in table B results from a single multiplier in table A (row 2 column $y$) applied to several multiplicands in row 1, beginning at the same column ($y$) and moving rightward. The reduction of work occurs because \citet [117] {Doolittle1878-HM} performed multiplication using the 3-digit tables of \citet {Crelle1864-HM}. Since he reused the multiplier for an entire row of calculations, Doolittle could open Crelle's tables to the one page for that multiplier and all the multiplicands. \citet [93] {Schott1879-HM} emphasized that using multiplication tables was innovative, noting that ``logarithms are altogether dispensed with.''\footnote {This comment suggests that computers generally did use logarithms to evaluate Gauss's brackets. Although the logarithm of a sum is not usually expressed in terms of logarithms of the summands, the subtraction in equation (\ref {eqn:brackets}) could be done via special tables of ``Gauss's logarithms.''}

\begin {figure} 
\newcommand {\leveling} {\vrule depth0pt height4.2ex width0pt}
\begin {center}
\footnotesize
\setlength {\tabcolsep} {0.375em}
\begin {minipage}[t]{16.0em}
\centering
{\small Table C}\\[1.0ex]
\begin {tabular} {c c c c}
\xss \parbox[b]{2.5em}{\centering recip\-rocal}& \leveling $x$& $y$& $z$\\ \hline
\sdd\disp {-1\over[aa]}& \disp {-{[ab]\over[aa]}}& \disp {-{[ac]\over[aa]}}& \disp {-{[ad]\over[aa]}}\\ 
\ddd\disp {-1\over[bb,1]}&& \disp {-{[bc,1]\over[bb,1]}}& \disp {-{[bd,1]\over[bb,1]}}\\
\ddd\disp {-1\over[cc,2]}&&& \disp {-{[cd,2]\over[cc,2]}}\\ 
\ddd\disp {-1\over[dd,3]}\\ \hline
\end {tabular}
\end {minipage}
\hspace {2em}
\begin {minipage}[t]{21,5em}
\centering
{\small Table D}\\[1.0ex]
\begin {tabular} {c c c c c}
\xss \leveling $w$& $x$& $y$& $z$\\ \hline
\sdd\disp {-{[an]\over[aa]}}& \disp {-{[bn,1]\over[bb,1]}}& \disp {-{[cn,2]\over[cc,2]}}& \disp {-{[dn,3]\over[dd,3]}}\\ \cline{4-4}
\ddd\disp {z \times -{[ad]\over[aa]}}& \disp {z \times -{[bd,1]\over[bb,1]}}& \disp {z \times -{[cd,2]\over[cc,2]}}& \parbox{4.0em}{\centering column sum${}=z$}\\ \cline {3-3}
\ddd\disp {y \times -{[ac]\over[aa]}}& \disp {y \times -{[bc,1]\over[bb,1]}}& \parbox{4.0em}{\centering column sum${}=y$}\\ \cline {2-2}
\ddd\disp {x \times -{[ad]\over[aa]}}& \parbox{4.0em}{\centering column sum${}=x$}\\ \cline {1-1}
\ddd\parbox{4.0em}{\centering column sum${}=w$}
\end {tabular}
\end {minipage}
\end {center}
\caption {How Doolittle performed back-substitution. Table C he copied from table A. The sum of each column in table D evaluates the ``explicit function'' for a different variable.}
\label {fig:tablesCandD}
\end {figure}

The back-substitutions was performed with similar economy. \citet [116] {Doolittle1878-HM} distinguished between numbers used once and those used many times, so he prepared for the substitutions by copying the reciprocals and ``explicit function'' coefficients from table A to table C (Figure \ref{fig:tablesCandD}). The value of $z$ was available in the final row of table A. The remaining variables were evaluated in table D where each row consists of one multiplier applied to one column, this time, of table C. The sums of the columns in table D give the other variables.

Doolittle's method included several contributions for which he is not now credited. He owed his speed in part to using just 3-digit arithmetic in the multiplication tables, and hence everywhere in the calculation. This could, as Doolittle appreciated, introduce severe rounding errors.\footnote {Rounding errors are relatively much smaller for modern computing equipment that adheres to international standards which require using roughly 7 or 14 decimal digits in calculations \citep {ieee1985}.} Thus, an important aspect of Doolittle's method was the ability to correct the 3-digit approximate solution with comparatively little extra work. Since the angle adjustment problem itself corrected numbers that were approximately known --- the measured angles --- it may have seemed natural to further correct the adjustments. Both \citet {Doolittle1878-HM} and \citet [93] {Schott1879-HM} described the correction process without giving it a name; today it is called \textit {iterative improvement\/} or \textit {iterative refinement\/}.\footnote {This process has both practical and theoretical importance in numerical analysis \citep [231, chap.\ 12] {Higham2002}. In matrix notation that is consistent with the equations (\ref {eqn:DoolittleNormalEquations}), suppose $\bar x$ is an approximate solution to the simultaneous linear equations $Ax + b = 0$, where $A$ and $b$ are a given matrix and vector, and where $x$ is the true vector solution. The correction to $\bar x$ can be obtained by solving $A e + r = 0$ where $r = A \bar x + b$ is called the \textit {residual\/} vector. It is straightforward to see that $x = \bar x + e$, but the calculation of $e$ will incur some errors so in principle the correction process may need to be repeated, hence \textit {iterative\/} improvement.}

\begin {figure}
\begin {center}
\footnotesize
\setlength {\tabcolsep} {0.333em}
\begin {minipage}[t]{22.0em}
\centering
{\small Table E}\\[1.0ex]
\begin {tabular} {c c c c}
\xss 1.& 2.& 3.& 4.\\ \hline
\sdd $r_1$& $r_2$& $r_3$& $r_4$\\ \cline{1-1}
\ddd\parbox{4.0em}{\centering column sum${}=s_1$}& \disp {s_1 \times -{[ab]\over[aa]}}& \disp {s_1 \times -{[ac]\over[aa]}}& \disp {s_1 \times -{[ad]\over[aa]}}\\ \cline{2-2}
\ddd& \parbox{4.0em}{\centering column sum${}=s_2$}& \disp {s_2 \times -{[bc,1]\over[bb,1]}}& \disp {s_2 \times -{[bd,1]\over[bb,1]}}\\ \cline{3-3}
\ddd&& \parbox{4.0em}{\centering column sum${}=s_3$}& \disp {s_3 \times -{[cd,2]\over[cc,2]}}\\ \cline{4-4}
\ddd&&& \parbox{4.0em}{\centering column sum${}=s_4$}\\ 
\end {tabular}
\end {minipage}
\hspace {1em}
\begin {minipage}[t]{23em}
\centering
{\small Table F}\\[1.0ex]
\begin {tabular} {c c c c c}
\xss $w$& $x$& $y$& $z$\\ \hline
\sdd\disp {s_1 \times {-1\over[aa]}}& \disp {s_2 \times {-1\over[bb,1]}}& \disp {s_3 \times {-1\over[cc,2]}}& \disp {s_4 \times {-1\over[dd,3]}}\\ \cline{4-4}
\ddd\disp {z_2 \times -{[ad]\over[aa]}}& \disp {z_2 \times -{[bd,1]\over[bb,1]}}& \disp {z_2 \times -{[cd,2]\over[cc,2]}}& \parbox{4.0em}{\centering column sum${}=z_2$}\\ \cline {3-3}
\ddd\disp {y_2 \times -{[ac]\over[aa]}}& \disp {y_2 \times -{[bc,1]\over[bb,1]}}& \parbox{4.0em}{\centering column sum${}=y_2$}\\ \cline {2-2}
\ddd\disp {x_2 \times -{[ad]\over[aa]}}& \parbox{4.0em}{\centering column sum${}=x_2$}\\ \cline {1-1}
\ddd\parbox{4.0em}{\centering column sum${}=w_2$}
\end {tabular}
\end {minipage}
\end {center}
\caption {How Doolittle performed iterative refinement. Table E corresponds to the final column of table B. The calculations in table B apply to the constant terms, while in table E they apply to the residuals. Table F is similar in construction to table D.}
\label {fig:tablesEandF}
\end {figure}

In describing the refinement process, Doolittle wrote $w_1$, $x_1$, $y_1$, $z_1$ for the values obtained from tables A--D. When these approximations are substituted into the normal equations they give residual values.
\begin {displaymath}
\begin {array} {c c c c c c c c c c c c c c}
r_1 \equals& [aa]& w_1& \plus& [ab]& x_1& \plus& [ac]& y_1& \plus& [ad]& z_1& \plus& [an]
\\{}
r_2 \equals& [ab]& w_1& \plus& [bb]& x_1& \plus& [bc]& y_1& \plus& [bd]& z_1& \plus& [bn]
\\{}
r_3 \equals& [ac]& w_1& \plus& [bc]& x_1& \plus& [cc]& y_1& \plus& [cd]& z_1& \plus& [cn]
\\{}
r_4 \equals& [ad]& w_1& \plus& [bd]& x_1& \plus& [cd]& y_1& \plus& [dd]& z_1& \plus& [dn] \end {array}
\end {displaymath}
Doolittle emphasized that this particular calculation had to be carried out to a high degree of accuracy after which the residuals $r_1$, $r_2$, $r_3$, $r_4$ (my notation) could be rounded to three digits for the remaining steps. In table E (Figure \ref {fig:tablesEandF}) he repeated the calculations on the residuals that he had performed on the constant terms of the normal equations in the final column of table B. In Table F (Figure \ref {fig:tablesEandF}) he duplicated the back-substitution of table D, but in this case for the corrections, which Doolittle named $w_2$, $x_2$, $y_2$, $z_2$. The corrected solutions, $w = w_1 + w_2$, $x = x_1 + x_2$, etc., were accurate to two or three digits in Doolittle's example. Note, these solutions still must be entered into the correlate equations to obtain the angle adjustments.

Another innovation, ``one of the principal advantages,'' was a provision to include new equations and variables. \citet [117--118] {Wright1906-HM} provide a clearer text than \citet [17--18] {Doolittle1878-HM}, whose description is somewhat brief. Doolittle suggested large problems could be solved by successively including conditions in the minimization problem, that is, by appending equations and variables to the normal equations. He recommended ordering the conditions so as to preserve zeroes in the elimination, and he suggested an ordering based on the geometric interpretation of the conditions. Doolittle thus anticipated the work on sparse matrix factorizations that would be done a hundred years later, for example \citet {George1981} and \citet {Duff1986}. 

\begin {figure} [t]
\centering
\begin {minipage} [t] {0.45\textwidth}
\centering
\includegraphics [scale=1] {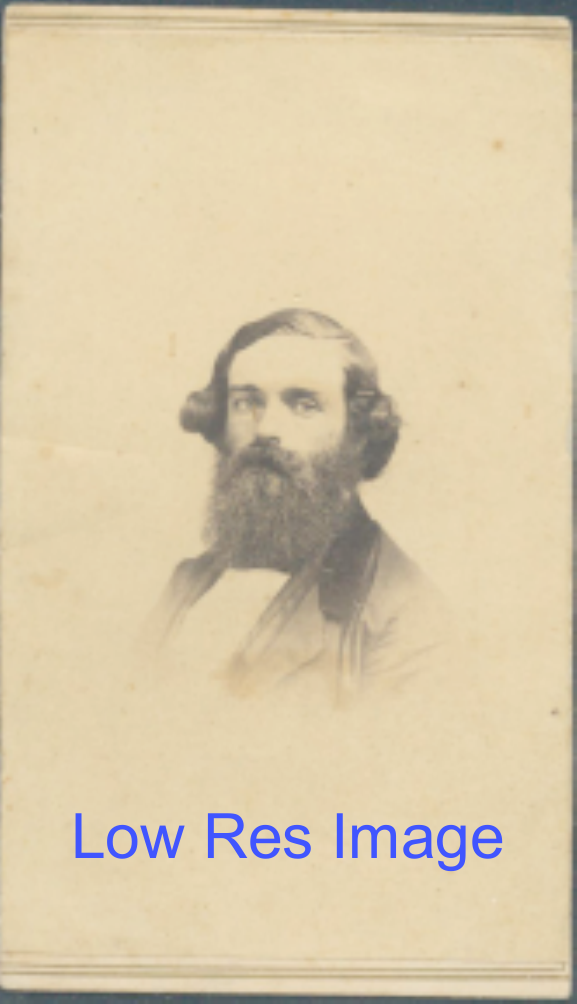}%
\caption {Myrick Hascall Doolittle, 1830--1911, circa 1862. Courtesy of the Antiochiana collection at the Olive Kettering Library of Antioch University.}
\label {fig:Doolittle}
\end {minipage}
\hfil
\begin {minipage} [t] {0.45\textwidth}
\centering
\includegraphics [scale=1] {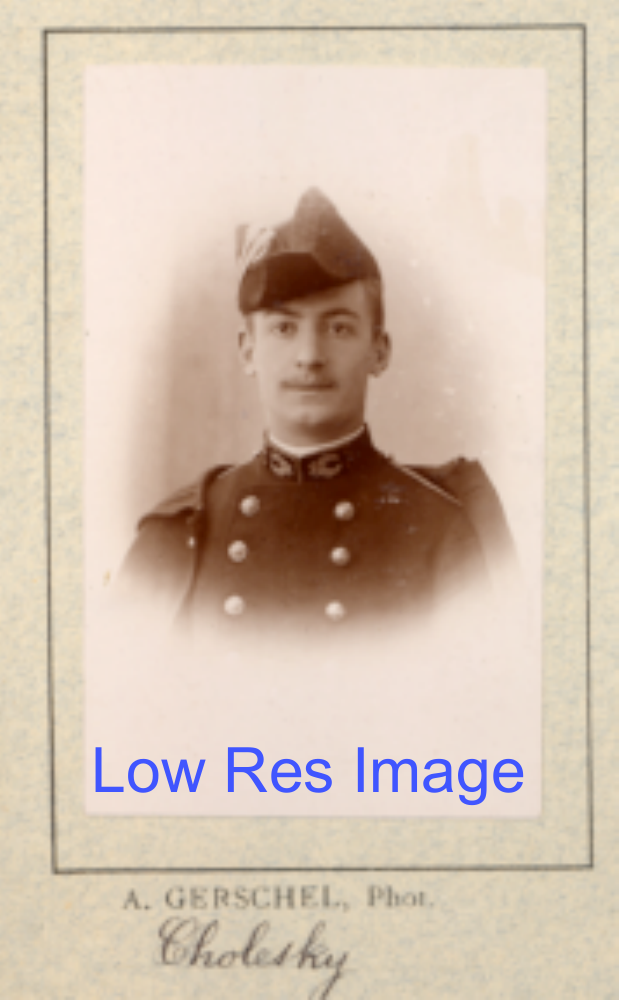}%
\caption {Andr\'e-Louis Cholesky, 1875--1918, at the start of his student days at the \'Ecole Polytechnique, 1895--1897. Courtesy of the Archives of the Biblioth\`eque Centrale of the \'Ecole Polytechnique.}
\label {fig:CholeskyPortrait}
\end {minipage}
\end {figure}

\citet [112] {Dwyer1941a} remarked that ``from Doolittle down to the present'' no formal proof had ever been offered that Doolittle's tables do solve the normal equations. Some justification is needed because Doolittle did not explicitly reduce to zero the coefficients of eliminated variables. For example, the combination of the $\first$ and $\second$ equations with the $\third$ involves a choice of multipliers to eliminate $w$ and $x$. In rows 8 and 9 of table B, Doolittle ignored the variables to be eliminated, yet he applied the proper multipliers to update coefficients for the retained variables $y$ and $z$. This saving is possible because, thanks to the underlying symmetry, Doolittle knew that the multipliers are also the coefficients of his ``explicit functions.'' \citet {Dwyer1941b} gave one proof, and other explanations could be given, such as expanding Gauss's bracket formulas for the coefficients in his equations $A = 0$, $B = 0$, \dots. The relationship between Gauss's brackets and Gaussian elimination seems to have been known, as evidenced by \citet [530] {Chauvenet1868-HM}. Although Doolittle did not express his method symbolically, he undoubtedly had the training and ability to do so.\footnote {After receiving a Bachelor's degree at Antioch College, Doolittle taught mathematics there, and then he studied under Benjamin Peirce at Harvard College.}

Doolittle's reticence in publishing his method may be because discussing computing methods per se was thought to be neither appropriate (as judged by the mathematical community) nor desirable (from the standpoint of the computer). \citet [156] {Grier2005} concluded from his study of hand computers that until the \xx\ century computing was a craft skill passed from masters to apprentices. A reluctance to disclose methods is consistent with Grier's picture of journeymen computers. Indeed, \citet [93] {Schott1879-HM} emphasized that Doolittle's paper of 1878 was written at the express request of the Coast Survey Director. The \citet {WashingtonStar1913-HM} reported Doolittle ``contributed numerous papers on his favourite subject'' to the Philosophical Society of Washington, yet of Doolittle's three publications listed in the bibliography compiled by \citet [365] {Gore1889}, calculating was discussed only in \citep {Doolittle1878-HM}.

Doolittle's 1878 paper made a strong impression on computers. The claims made for the ``Coast Survey method'' by Doolittle and by \citet {Schott1879-HM} caused \citet {Werner1883} to examine the paper immediately and skeptically. He calculated Doolittle's tables three different ways: with Crelle's multiplication tables, with logarithms, and with the Thomas Arithmometer, an early calculating machine.\footnote {The greater accuracy of the latter two methods obviated the need for iterative refinement which was neglected by later authors.} The Doolittle tables evidently passed muster because they were reprinted for many years. \citet [iv, 65] {Jordan1895-HM} remarked on the ``old, classic'' notation of Gauss and then exhibited without attribution a calculation using logarithms and Doolittle's table B.\footnote {\label {fn:Jordan}Willhelm Jordan's encyclopedic reference work on geodesy was constantly updated. He credited the material to scarcely anyone by name except for Gauss. \citet {Althoen1987} explain that it is this Jordan, from Hanover, who was responsible for the Gauss-Jordan algorithm for inverting matrices. \citet {Katz1988} has a biographical sketch of Jordan.} \citet {Wright1906-HM} offered the Gauss brackets and the Doolittle tables as competing approaches. They did not understand the subtleties because they touted the use of Crelle's 3-digit tables (p.\ 120) yet they omitted the iterative refinement. 

Doolittle continued to be associated with computing in the \xx\ century even after computers adopted mechanical calculators and least squares became known as regression analysis. Howard Tolley, a computer at the Coast and Geodetic Office who became an official at the United States Department of Agriculture, politely scolded economists and statisticians when they neglected to credit geodesists \citep [159--164] {Grier2005}. \citet [497] {Tolley1927} cited several textbooks teaching Doolittle's method, and they illustrated what were still essentially Doolittle's tables by a calculation that, by then, was done with a Monroe four-function electric calculator.\footnote {\citet [500] {Tolley1927} reported that Miss Helen Lee, a computer, could use a form of Doolittle's method to solve 5 normal equations with 6-digit numbers (two whole and four fractional) in 50 minutes (or 40 minutes if rounded to two fractional digits).} Paul \citet {Dwyer1941a} (who is discussed in section \ref {sec:part5}) found enough similarities between Doolittle's calculations and a class of matrix computing methods to name the methods after Doolittle, in which form Doolittle is commemorated to this day. 




\subsection {Mechanical Calculators}
\label {sec:calculators}

Just as multiplication tables made possible Doolittle's method in the last quarter of the \xix\ century, calculating machines enabled other methods at the beginning of the next century.\footnote {The same of course can be said of electronic computers later in the \xx\ century. Innovation in mathematical practice can thus depend on manufacturing.} The most pressing need was help with multiplication, a tedious operation and the one most frequently encountered. The multiplying machines most common in scientific calculations were based first on Leibniz stepped drums and then on pinwheel gears \citep [ xiv--xv] {Williams1982-HM}.\footnote {\citet {Baxendall1929} and \citet {Murray1948, Murray1961} describe the mechanisms, and \citet {Marguin1994} shows pictures of the computing instruments.} 

The first commercial multiplier was the Thomas Arithmometer \citep {vanEyk1854-HM} based on Leibniz drums. Widespread acceptance came only in the late 1870s, and the mechanism was never sufficiently rugged to be entirely reliable in the \xix\ century \citep [17 and 19] {Johnson1997}. Indeed, the machines the Coast Survey first acquired in 1890 performed only addition \citep [93] {Grier2005}. Thus, \xix\ century technology offered little improvement over logarithm and multiplication tables. 

The first mass-produced multiplying machines were based on the version of the pinwheel mechanism that Willgodt Odhner patented in 1890 \citep {Apokin2001}. The explosive growth in the market for these machines soon supported many European manufacturers such as Brunsviga, Dactyle, and Original Odhner. \citet [84] {Whipple1914-HM} reports that by 1912 Brunsviga alone had manufactured 20,000 pinwheel calculators. 

The ability to accumulate sums of products was concomitant with machine multiplication. Multiplying calculators performed ``long multiplication'' in which the multiplicand shifts, to align with different digits of the multiplier, for repeated addition to an accumulator.\footnote {In fact, contrary to paper and pencil multiplication, the accumulator and multiplier moved rather than the multiplicand.} The accumulator need not contain zero at the start, so a product could increment whatever value was already there; a sum of products could be formed in this way. 

Calculations that were arranged to accumulate products were both faster and more accurate. Doolittle's table B could be discarded because it recorded only products that were to be summed. Further, the accumulator had more fractional digits than the multiplicand register, which was also the addend for addition. Any number removed from the accumulator had to be rounded before it could be reentered later. Therefore, as \citet [112] {Dwyer1941a} explained, summing products without reentry was ``theoretically more accurate'' because it gave ``the approximation resulting from a number of operations rather than the combination of approximations resulting from the operations.''

Electrified versions of calculators began to appear after 1920 \citep {Swartzlander1995}. Availability was limited by cost during the Great Depression and then by rationing during World War II. The Mathematical Tables Project gave most computers only paper and pencils because the cost of an electric machine rivaled the annual salary of a computer \citep [220] {Grier2005}.\footnote {The celebrated but short lived American Mathematical Tables Project project commenced in 1938 to aid the unemployed and closed in 1948 having aided the war effort.} In the United States, just three manufacturers of multiplying machines survived after the war: Friden from San Leandro, California, Monroe of Orange, New Jersey, and the ultra quick and quiet Marchant ``silent speed'' calculators from Oakland, California \citep {Chase1980}. 

Another multiplying device was the multiplying punch that was introduced in the 1930s. It was a variation of commercial punch card equipment, and access was again limited by high cost. These machines could solve $10$ simultaneous linear equations in about $4$ hours \citep [462] {Verzuh1949}. A literature developed for using them in scientific and statistical calculations, but nothing suggests it influenced Gaussian elimination.\footnote {For examples of punch card computing methods besides Verzuh see \citet {Brandt1935} and \citet {Eckert1940-HM}. The latter was the famous ``Orange Book'' of calculating methods for punch card equipment.} 

\subsection {Cholesky: Machine Algorithm} 
\label {sec:Cholesky}

The first algorithm intended for a machine appears to be that developed by the military geodesist and World War I casualty Andr\'e-Louis Cholesky.\footnote {For a biography and a discussion of his work see \citet {BrezinskiGross2005}; see also \citet {Brezinski2006}.} Like Doolittle, Cholesky calculated Gauss's angle adjustments that were formulated as an underdetermined least squares problem (case \ref {case:two} of the least squares method described in section \ref {sec:squares}). Commandant \citet {Benoit1924-HM} published his colleague's method posthumously.\footnote {Benoit served with Cholesky as a military geodesist and also published a eulogy for him.} A similar but not identical manuscript was recently found among Cholesky's military papers and has now been published \citep {Cholesky2005-HM}. 

Although Cholesky's method is most simply expressed using matrices which were known when he developed it, he did not use them, so his invention is better understood in algebraic notation. Cholesky wrote the condition equations as \citep [eqn.\ 1] {Benoit1924-HM}, 
\begin {equation}
\label {eqn:original}
\renewcommand {\arraystretch} {1.25}
\setlength {\arraycolsep} {0.0625em}
\begin {array} {l c c l c c l c c c c c c c c c}
a_1& x_1& \plus& a_2& x_2& \plus& a_3& x_3& \plus& \dots& \plus& a_n& x_n& \plus& K_1& \equals 0
\\
b_1& x_1& \plus& b_2& x_2& \plus& b_3& x_3& \plus& \dots& \plus& b_n& x_n& \plus& K_2& \equals 0
\\
\multicolumn {2} {c} {\vdots}&& \multicolumn {2} {c} {\vdots}&& \multicolumn {2} {c} {\vdots}&&&& \multicolumn {2} {c} {\vdots}&& \vdots
\\
\ell_1& x_1& \plus& \ell_2& x_2& \plus& \ell_3& x_3& \plus& \dots& \plus& \ell_n& x_n& \plus& K_p& \equals 0
\end {array}
\end {equation}
where $n > p$. He wrote the normal equations as \citep [eqn.\ 5] {Benoit1924-HM}, 
\begin {equation}
\label {eqn:normal-Cholesky}
\renewcommand {\arraystretch} {1.25}
\setlength {\arraycolsep} {0.0625em}
\begin {array} {c c c c c c c c c c c c c}
a^1_1& \lambda_1& \plus& a^2_1& \lambda_2& \plus& \dots& \plus& a^p_1& \lambda_p& \plus& K_1& \equals 0
\\
a^1_2& \lambda_1& \plus& a^2_2& \lambda_2& \plus& \dots& \plus& a^p_2& \lambda_p& \plus& K_2& \equals 0
\\
\multicolumn {2} {c} {\vdots}&& \multicolumn {2} {c} {\vdots}&&&& \multicolumn {2} {c} {\vdots}&& \vdots
\\
a^1_p& \lambda_1& \plus& a^2_p& \lambda_2& \plus& \dots& \plus& a^p_p&
\lambda_p& \plus& K_p& \equals 0
\end {array}
\end {equation}
where $a_j^k = a_k^j$ is the sum of products of coefficients in the $\jth$ and $\kth$ conditions, for example, $\smash {a_2^p} = b_1 \ell_1 + \dots + b_n \ell_n$. 

\newcommand {\benoit} {\beta}

Cholesky's remarkable insight was that, since many underdetermined systems share the same normal equations, for any normal equations there may be some condition equations that can be directly solved more easily \citep [70] {Benoit1924-HM}. He found his alternate equations in the convenient, triangular form (introducing new unknowns, $y$ in place of $x$, and new coefficients, $\benoit$ in place of $a$).\footnote {\citet {Benoit1924-HM} unfortunately used $\alpha$ for the new coefficients which here has been replaced by $\benoit$ to more easily distinguish them from the original coefficients $a$.}
\begin {equation}
\label {eqn:lower}
\renewcommand {\arraystretch} {1.25}
\setlength {\arraycolsep} {0.0625em}
\begin {array} {c c c c c c c c c c c c c c c c}
\benoit^1_1& y_1&&&&&&&&&&&& \plus& K_1& \equals 0
\\
\benoit^1_2& y_1& \plus& \benoit^2_2& y_2&&&&&&&&& \plus& K_2& \equals 0
\\
\benoit^1_3& y_1& \plus& \benoit^2_3& y_2& \plus& \benoit^3_3&
y_3&&&&&& \plus& K_3& \equals 0
\\
\multicolumn {2} {c} {\vdots}&& \multicolumn {2} {c} {\vdots}&& \multicolumn {2} {c} {\vdots}&&{\ddots}&&&&& \vdots
\\
\benoit^1_p& y_1& \plus& \benoit^2_p& y_2& \plus& \benoit^3_p& y_3& \plus& \dots& \plus& \benoit^p_p&
y_p& \plus& K_p& \equals 0
\end {array}
\end {equation}
Cholesky discovered that the coefficients in equation (\ref {eqn:lower}) are given by straightforward formulas \citep [72] {Benoit1924-HM}.
\begin {equation}
\label {eqn:coefficients}
\begin {array} {r c l}
\benoit^i_i& \equals& \sqrt { \, \strut \smash {a^i_i - (\benoit^1_i)^{\,2}_{\strut} - (\benoit^2_i)^{\,2}_{\strut} - \dots - (\benoit^{i-1}_i)^{\,2}_{\strut}}}\\ \noalign {\bigskip}
\benoit^i_{i+r}& =& \displaystyle {a^i_{i+r} - \benoit^1_i \, \benoit^1_{i+r} - \benoit^2_i \, \benoit^2_{i+r} - \dots - \benoit^{i-1}_{i+r} \, \benoit^{i-1}_{i+r} \over \benoit^i_i}
\end {array}
\end {equation}
The solution $\lambda$ of the normal equations expresses the solution of the condition equations as a linear combination of the (transposed) coefficients in the condition equations. For Cholesky, these combinations become a system of equations to be solved for $\lambda$ \citep [eqn.\ 7] {Benoit1924-HM}.
\begin {equation}
\label {eqn:upper}
\renewcommand {\arraystretch} {1.25}
\setlength {\arraycolsep} {0.0625em}
\begin {array} {c c c c c c c c c c c c c c c c}
y_1& \equals& \benoit_1^1& \lambda_1& \plus& \benoit_2^1& \lambda_2& \plus& \dots& \plus& \benoit_p^1& \lambda_p 
\\
y_2& \equals&&&& \benoit_2^2& \lambda_2& \plus& \dots& \plus& \benoit_p^2& \lambda_p 
\\
\vdots&&&&&&&&{\ddots}&& \multicolumn {2} {c} {\vdots}
\\
y_p& \equals&&&&&&&&& \benoit_p^p& \lambda_p 
\end {array}
\end {equation}
Since the original condition equations (\ref {eqn:original}) and Cholesky's equations (\ref {eqn:lower}) have the same normal equations (\ref {eqn:normal-Cholesky}), the quantities $\lambda$ are the same for both problems. Once obtained, $\lambda$ can be used to evaluate $x$ as usual (from the transpose of the original coefficients), thereby solving the original condition equations. Thus Cholesky's method was first to form his new coefficients using equation (\ref {eqn:coefficients}), then to solve (\ref {eqn:lower}) by forward substitution for $y$, next to solve (\ref {eqn:upper}) by backward substitution for $\lambda$, and finally to evaluate $x$.
\begin {equation}
\label {eqn:transpose}
\renewcommand {\arraystretch} {1.25}
\setlength {\arraycolsep} {0.0625em}
\begin {array} {c c c c c c c c c c c c c c c c}
x_1& \equals& a_1& \lambda_1& \plus& b_1& \lambda_2& \plus& \dots& \plus& \ell_1& \lambda_p 
\\
x_2& \equals& a_2& \lambda_1& \plus& b_2& \lambda_2& \plus& \dots& \plus& \ell_2& \lambda_p 
\\
x_3& \equals& a_3& \lambda_1& \plus& b_3& \lambda_2& \plus& \dots& \plus& \ell_3& \lambda_p 
\\
\vdots&& \multicolumn {2} {c} {\vdots}&& \multicolumn {2} {c} {\vdots}&&&& \multicolumn {2} {c} {\vdots}
\\
x_n& \equals& a_n& \lambda_1& \plus& b_n& \lambda_2& \plus& \dots& \plus& \ell_n& \lambda_p 
\end {array}
\end {equation}

Benoit published Cholesky's method in his own format with revised notation and with a table showing how the calculation was conducted (Figure \ref {fig:Cholesky}). He indicated the $\ith$ column of new coefficients (from $\smash {\benoit^i_i}$ to $\smash {\benoit^i_{p}}$) should be calculated before going on to the next column. The table is more compact than Doolittle's tables because much of the intermediate work is not recorded: each $\benoit^i_{i+r}$ is a sum of products that can be accumulated with a machine. Both Benoit and Cholesky mentioned using calculators to form these sums, and Cholesky specifically referred to the Dactyle brand of the Odhner design. \citet {Cholesky2005-HM} reported solving 10 equations with 5-digit numbers in 4 to 5 hours: he was as quick as a multiplying punch although he worked with fewer digits.

\begin {figure}[t]
\renewcommand {\baselinestretch} {1.0}
\renewcommand {\arraystretch} {1.0}
\begin {center}
\footnotesize
\begin {tabular} {c | c c c c c c c c c | c c |}
\multicolumn {2} {c |} {\xss}& $\lambda_1$& $\lambda_2$& $\dots$& $\lambda_{i-1}$& $\lambda_{i}$& $\dots$& $\lambda_{p-1}$& $\lambda_p$& $K$& $\lambda$\\ \hline
1\sss& \multicolumn {1} {c |} {$\benoit_1^1$}& $a_1^1$& $a_1^2$&& $a^{i-1}_1$& $a^i_1$&& $a^{p-1}_1$& $a^p_1$& $K_1$& $\lambda_1$\\ \cline {3-3}
2\sss& $\benoit_2^1$& \multicolumn {1} {c |} {$\benoit_2^2$}& $a_2^2$&& $a^{i-1}_2$& $a^i_2$&& $a^{p-1}_2$& $a^p_2$& $K_2$& $\lambda_2$\\ \cline {4-4}
\sss& \multicolumn {3} {c |} {}&&&&&&&&\\ 
\sss& \multicolumn {4} {c |} {}&&&&&&&\\ \cline {6-6}
$i$\sss& $\benoit_i^1$& $\benoit_i^2$& $\benoit_i^3$&& \multicolumn {1} {c |} {$\benoit_i^i$}& $a_i^i$&& $a_i^{p-1}$& $a_i^p$& $K_i$& $\lambda_i$\\ \cline {7-7}
\sss& \multicolumn {6} {c |} {}&&&&&\\ 
\sss& \multicolumn {7} {c |} {}& $a_{p-1}^{p-1}$& $a_{p-1}^{p}$&&\\ \cline {9-9}
$p$\sss& $\benoit_p^1$& $\benoit_p^2$& $\benoit_p^3$&& $\benoit_p^i$&&& \multicolumn {1} {c |} {$\benoit_p^p$}& $a_p^p$& $K_p$& $\lambda_p$\\ \hline
\sss $y$& $y_1$& $y_2$& $y_3$&& $y_i$&&& \multicolumn {1} {c |} {$y_p$}\\ \cline {1-9}
\end {tabular}
\end {center}
\caption {How Cholesky may have organized the calculation to solve the normal equations by what is now known as the Cholesky factorization, from the fold-out table of \citet {Benoit1924-HM}. The left column and top row are labels; the other positions would be occupied by numbers. The coefficients, $a$, and the constants terms, $K$, of the normal equations are placed for reference above the diagonal. The coefficients, $\benoit$, of Cholesky's manufactured condition equations are written below the diagonal. The solution, $y$, of his condition equations stands in the bottom row. The solution, $\lambda$, of the normal equations is recorded in the right column. Some additional rows and columns for arithmetic checks and for accuracy estimates have been omitted.}
\label {fig:Cholesky}
\end {figure}

For 20 years after the publication by \citet {Benoit1924-HM}, Cholesky's method was seldom discussed. It was, however, used in Sweden during this time by the geodesist Tryggve \citet {Rubin1926} who augmented it to produce other statistically relevant quantities.\footnote {\citet [286] {Brezinski2006} claims incorrectly that after \citet {Benoit1924-HM} ``a period of 20 years followed without any mention of the work.''}
\begin {quotation}
\noindent The normal equations have been solved by the Cholesky-Rubin method, which offers the advantage that the solution is easily effected on a calculating machine (Cholesky) and that the most probable values and their mean errors are derived simultaneously (Rubin). 
\same --- \citet [30] {AhlmannRosenbaum1933}
\end {quotation}
\citet [22] {Jensen1944}, writing in a Danish geodetic publication, remarked that Cholesky's method of solving the normal equations ``ought to be more generally used.'' Soon thereafter it was independently discovered in matrix form by Dwyer and was put to use in America (see section \ref {sec:part5}). 

\subsection {Crout: ``each element is determined by one continuous machine operation''}
\label {sec:Crout}




Prescott Crout was a professor of mathematics at the Massachusetts Institute of Technology with an interest in mathematics for electrical engineering. In a paper that is model of brevity, \citet {Crout1941a} listed the gamut of uses for linear equations that had developed after least squares, and described a new method to solve general systems of equations rather than only normal equations from the method of least squares. Like Cholesky, Crout formulated his method to emphasize sums of products.

\citet [1239] {Crout1941a} explained ``the method was originally obtained by combining the various processes which comprise Gauss's method, and adapting them for use with a computing machine.'' He wrote the coefficients in a rectangular matrix. Crout meant ``matrix'' in the sense of ``table'' with the constant terms in the final column; he did not employ matrix algebra. His method consisted of a few terse rules for transforming the numbers. So much had changed in the application of Newton's and Rolle's elimination rules that all semblance of symbolic algebra had vanished. Crout's pithy instructions are here restated more expansively.
\begin {enumerate}
\item The first column is left unchanged. The first row to the right of the diagonal is divided by the diagonal entry.
\item An entry on or below the diagonal is reduced by the sum of products between entries to the left in its row and the corresponding entries above in its column. This calculation is permitted only after all those entries themselves have been transformed.
\item As in 2, similarly for an entry above the diagonal except it is lastly divided by the diagonal entry in its row.
\item On completing the above steps, (condensing further instructions) conduct a back-substitution using the coefficients above the diagonal and the final column of transformed constants.
\end {enumerate}
These steps were obviously intended for use by computers who had calculators that could accumulate sums of products. Like Doolittle before him, Crout included instructions to improve the accuracy of the solution by an unnamed iterative refinement process.\footnote {\citet [49] {Wilkinson1960} commented that the sum of product methods as originally stated could give poor results because they omitted provisions to reorder the equations and variables such as stated in section \ref {sec:today}. This contingency is unnecessary for normal equations but it can be important for others.} He concluded the paper with a rigorous proof that the method necessarily solves the intended equations. 

\begin {figure} [t]
\centering
\begin {minipage} [t] {0.45\textwidth}
\centering
\includegraphics [scale=1] {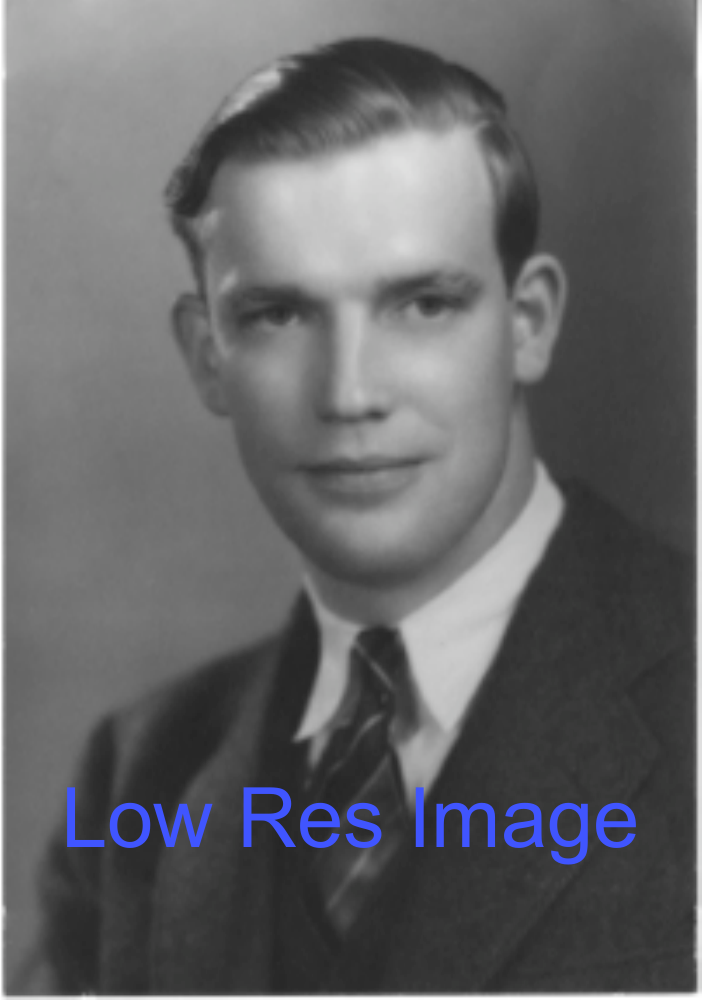}%
\caption {Prescott Durand Crout, 1907--1984, circa 1936. Courtesy MIT Museum.}
\label {fig:Crout}
\end {minipage}
\hfil
\begin {minipage} [t] {0.45\textwidth}
\centering
\includegraphics [scale=1] {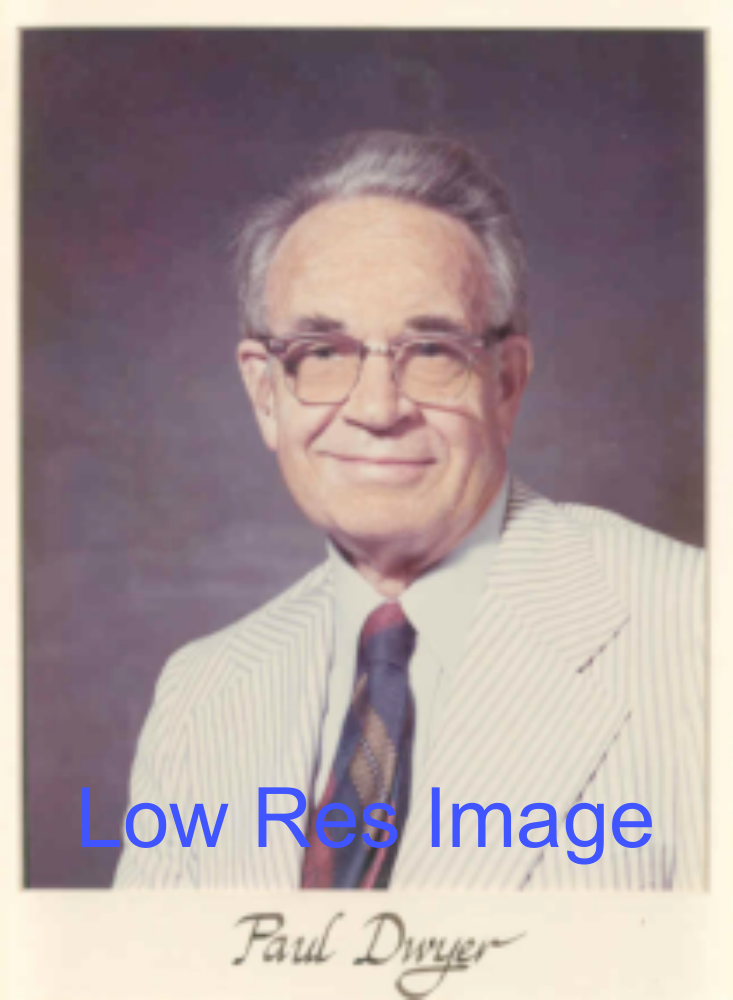}%
\caption {Paul Sumner Dwyer, 1901--1982, circa 1960s. Courtesy of Prof.\ Ingram Olkin of Stanford University.}
\label {fig:Dwyer}
\end {minipage}
\end {figure}

In the same year Paul \citet [108] {Dwyer1941a} published a survey of computing methods which included a method similar to Crout's, the ``abbreviated method of single division.'' Dwyer did not find the method in the literature but rather synthesized it to complete his classification scheme: for normal equations or not, with ``single division'' or not, and ``abbreviated'' or not. The latter was Dwyer's name for methods that accumulated sums. \citet [259] {Waugh1945} reported ``great recent interest'' in such ``compact'' methods, but they cited no recent literature other than their own papers and those of Prescott Crout.

Crout's method was rapidly adopted by those with access to calculators, and in the following decade his name appeared dozens of times in the literature. Crout himself reprised his paper \citep {Crout1941b-HM} for a series of manuals issued by the Marchant Company, a leading calculator manufacturer, and the Company prepared very detailed instructions for Crout's method in a subsequent manual  \citep {Marchant1941-HM}. \citet [352] {Wilson1952} particularly recommended the Marchant manuals to academic researchers. \citet {Black1949} cited both \citet {Crout1941a} and Doolittle as interpreted by \citet {Dwyer1941b} in describing how to use calculators for engineering work. 

\section {Describing Gaussian Elimination with Matrices}
\label {sec:part5}

\subsection {Overview of Matrix Descriptions}

The milieux of Gauss, Doolittle, Cholesky, and Crout, who worked on Gaussian elimination through symbolic algebra, was eventually brought to an end by those who worked with matrix algebra. The interpretation of Gaussian elimination through matrices led to a consolidation by showing how all the computing variations were trivially related through what are now called matrix decompositions. Henceforth one formalism would suffice for all methods.

\begin {figure}
\centering
\includegraphics [scale=1] {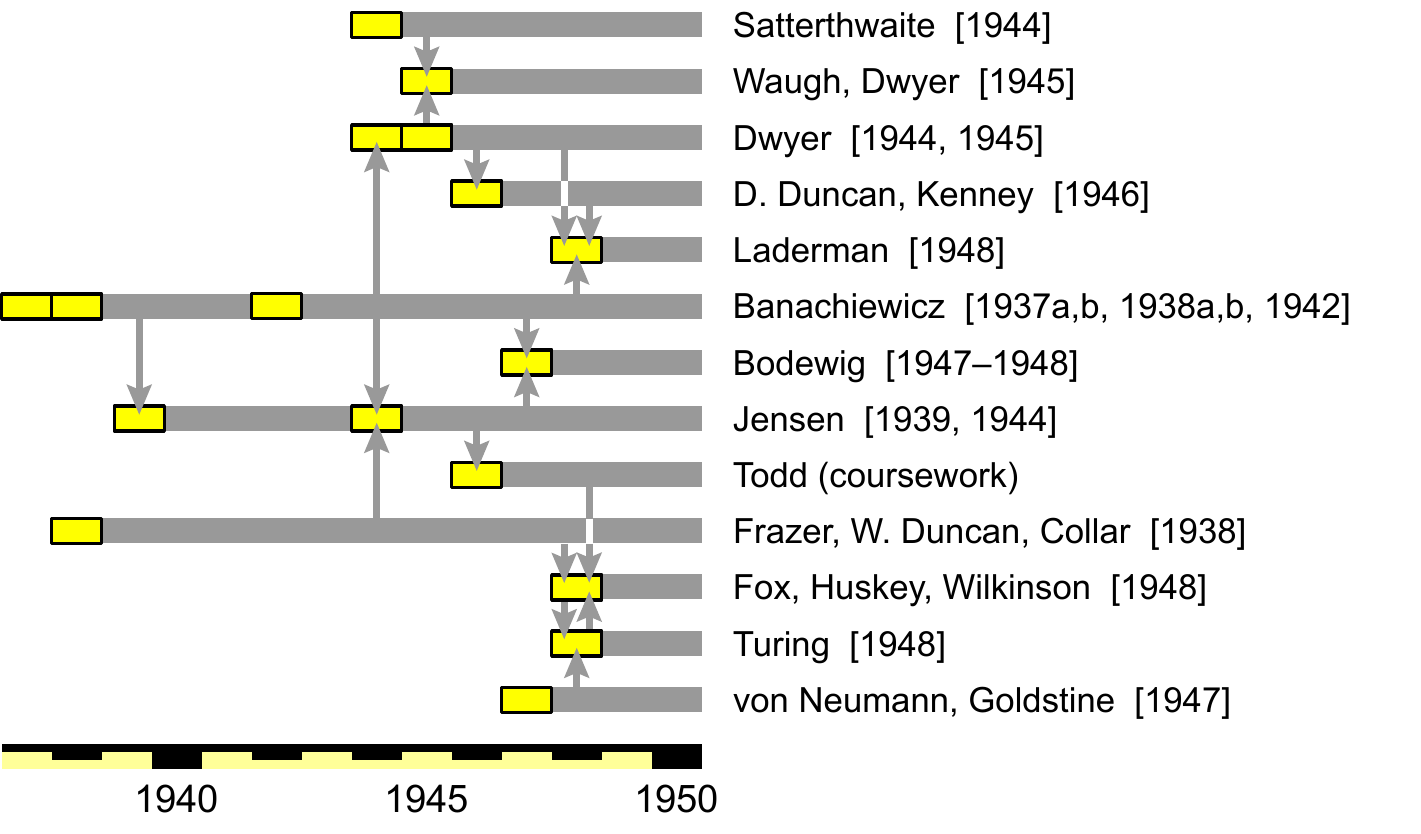}%
\caption {Citation patterns among the earliest authors who mention Gaussian elimination in the form of a matrix decomposition. Arrows show the flow of information to the citing papers at the head from the cited oeuvre at the tail. For example Dwyer in 1944 cited one or more works of Banachiewicz. Boxes indicate published papers, except in the case of Todd who taught a course in 1946 that is described by \citet {TausskyTodd2006}.}
\label {fig:factoring}
\end {figure}

\Citet [488] {vonNeumann1954-HM} cited matrices as an example of the delay between a mathematical discovery and its use outside mathematics. \citet {Hawkins1975, Hawkins1977a, Hawkins1977b} found matrix algebra was independently invented by Eisenstein (1852), Cayley (1858), Laguerre (1867), Frobenius (1878) and Sylvester (1881) to clarify subjects such as determinants, quadratic forms, and the elementary divisors that Weierstrass used to study ordinary differential equations. The first significant application of matrices outside pure mathematics was Heisenberg's matrix (quantum) mechanics in 1925. Although accounts of Gaussian elimination without matrix notation continued well into mid century,\footnote {For example \citet {Aitken1937, Aitken1951}, \citet {Crout1941b-HM, Crout1941a}, \citet {Dwyer1941b}, \citet {Hotelling1943}, \citet [422--423] {Bargmann1946-HM}, and \citet {MacDuffee1949-HM}.} from the end of the 1930s several authors offered matrix descriptions within a short span of time. The citation analysis in Figure \ref {fig:factoring} identifies six groups or individuals who made unprecedented contributions: Banachiewicz, Frazer et al., Jensen, Dwyer, Satterthwaite, and von Neumann and Goldstine.

\nocite {Banachiewicz1937b}%
\nocite {Banachiewicz1937a}%
\nocite {Banachiewicz1938c}%
\nocite {Banachiewicz1938a}%
\nocite {Banachiewicz1942}%
\nocite {Bodewig1947-HM}%
\nocite {Duncan1946}%
\nocite {Dwyer1944}%
\nocite {Dwyer1945}%
\nocite {Fox1948}%
\nocite {Jensen1939}%
\nocite {Jensen1944}%
\nocite {Laderman1948}%
\nocite {Satterthwaite1944}%
\nocite {Turing1948}%
\nocite {vonNeumann1947-HM}%
\nocite {Waugh1945}

\subsection {Toeplitz: First Matrix Decomposition} 
\label {sec:Toeplitz}

David Hilbert's study of integral equations inspired his students Erhard Schmidt and Otto Toeplitz to promote infinite matrices as a representation for linear operators on function spaces.\footnote {Another Hilbert prot\'eg\'e, John von Neumann, showed the infinite matrix approach to be fundamentally inadequate \citep [330] {Bernkopf1968}.} As part of that research, Toeplitz used determinants to examine the invertibility of infinite matrices, ``If one uses the symbolism of matrix calculus (see for example Frobenius),'' then any finite symmetric matrix $S$ with all leading principal determinants not zero, has a matrix $U$ with 
\begin {equation}
\label {eqn:toeplitz}
U^\prime U = S^{-1} \quad \mbox {equivalently} \quad U^{-1}
U^{\prime-1} = S
\end {equation}
(original notation) where $^\prime$ is transposition \citep [102] {Toeplitz1907-HM}.\footnote {Toeplitz remarked that Lagrange and Gauss knew such a decomposition for quadratic forms, and he cited \citet {Jacobi1857-HM} for a similar decomposition of bilinear forms.} In this way Toeplitz first exhibited what now is called the Cholesky factor, $U^{-1}$, though he left the entries in a computationally useless form given by determinants. It was clear from the determinantal formulas that $U$ was a lower triangular matrix. \citet [198] {TausskyTodd2006} suggested that equation (\ref {eqn:toeplitz}) was the first expression of any such formula in matrix notation. Moreover, the equation may be unique to Toeplitz because his $U$ and $S$ are related by inversion. Taussky and Todd thought a proof might be based on formulas of \citet [39] {Gantmacher1959-HM}. 

\subsection {Banachiewicz: Cracovian Algebra}
\label {sec:Banachiewicz}

One path to matrix algebra which has not been recognized is a motivation deriving from computing that prompted the Cracovian algebra of the polymath Tadeusz Banachiewicz.\footnote {The first of several biographical sketches of Banachiewicz, a concentration camp survivor, was given by \citet {Witkowski1955}.} Cracovians began as a notational scheme for astronomical calculations and became an algebra different from Cayley's for matrices of arbitrary dimension. The distinguishing feature of the Cracovian algebra is a column-by-column product which is more natural when calculating with columns of figures by hand.
\begin {quotation}
\noindent It must, however, be conceded that in practice it is easier to multiply column by column than to multiply row by column \dots. It may, in fact, be said, that the computations are made by cracovians and the theory by matrices 
\same --- \citet [5] {Jensen1944}
\end {quotation}
Banachiewicz posed least squares problems in terms of Cracovians for the purpose of improving computations.\footnote {For descriptions of the algorithms and additional references to the original papers, see \citet {Kocinski2004}.} \citet [3, 19] {Jensen1939} reports hearing Banachiewicz advocate this approach at meetings of the Baltic Geodetic Commission as early as 1933. \citet {Banachiewicz1938a} independently discovered Cholesky's method, and although later than Cholesky, he initially had greater impact. Banachiewicz inspired the work of Jensen, and he was widely cited: by \citet [45] {Jensen1944}, \citet [89] {Dwyer1944}, \citet [78] {Cassinis1946-HM}, \citet [part V, 90] {Bodewig1947-HM}, \citet {Laderman1948}, again by \citet [103] {Dwyer1951}, and by \citet [301] {Forsythe1953}. Nevertheless, he was neglected by the influential computational mathematician \citet {Householder1956, Householder1964}. 

\subsection {Frazer, Duncan, and Collar: Elementary Matrices}
\label {sec:Frazer}

The work of R.\ A.\ Frazer, W.\ J.\ Duncan, and A.\ R.\ Collar exemplifies the practical aspects of the source for matrix algebra that \citet {Hawkins1977b} finds in Weierstrass's study of linear differential equations.\footnote {See \citet {Pugsley1961}, \citet {Relf1961}, and \citet {Bishop1987} for biographies of Frazer, Duncan, and Collar, respectively.} In this case the stimulus for Frazer's interest in matrices may have been Henry Baker with whom Frazer studied as a Cambridge undergraduate.\footnote {\citet {Pugsley1961} remarks on the influence of Cayley's student Baker whose work included matrix formulations of differential equations such as \citet {Baker1903}.}  Collar and Duncan used matrix notation in their approximate formulations of structural problems in the precursor to the finite element method \citep {Felippa2001}.  The three collaborators thus translated into questions about matrices the problems they studied for Frazer's aerodynamics section at the National Physical Laboratory, notably problems of airframe vibration or flutter. Their textbook about matrices in engineering analysis became an international standard \citep [80] {Pugsley1961}.

\citet [96--99] {Frazer1938} viewed elimination as ``building up the reciprocal matrix in stages by elementary operations'' which could produce a triangular matrix ``such that its reciprocal can be found easily.'' They demonstrated column elimination of a $4 \times 4$ matrix $a$ (their notation). They had $a M_1 M_2 M_3 = \tau$ where the $M_i$ are ``post multipliers'' that effect the elementary operations and $\tau$ is ``the final triangular matrix,'' and they remarked that $M_1 M_2 M_3$ was itself a triangular matrix but ``opposite-handed'' from $\tau$. However, they did not continue their presentation to the modern conclusion, neither commenting on matrix factoring nor writing out a factorization such as $a = \tau \, (M_1 M_2 M_3)^{-1}$. 

The modern literature retains the emphasis on elementary operations that are accomplished by multiplication with \textit {elementary matrices\/}. \citet [13--15] {Jensen1944} borrowed the approach to establish the connection between Gaussian elimination and triangular factoring. He restated it in the more conventional form of row operators acting on the left but using the same notation $M_i$. Textbooks now use this row exposition in terms of ``premultiplication'' to establish the relationship between Gaussian elimination and matrix factoring, although they attribute it to neither Frazer et al.\ nor Jensen. For example \citet [94-95] {Golub1996} use, remarkably, the same notation $M_i$ but ahistorically call it a ``Gauss transformation.''

\subsection {Jensen: Synthesis and Conduit}
\label {sec:Jensen}

\begin {figure} [t]
\centering
\begin {minipage} [t] {0.45\textwidth}
\centering
\includegraphics [scale=1] {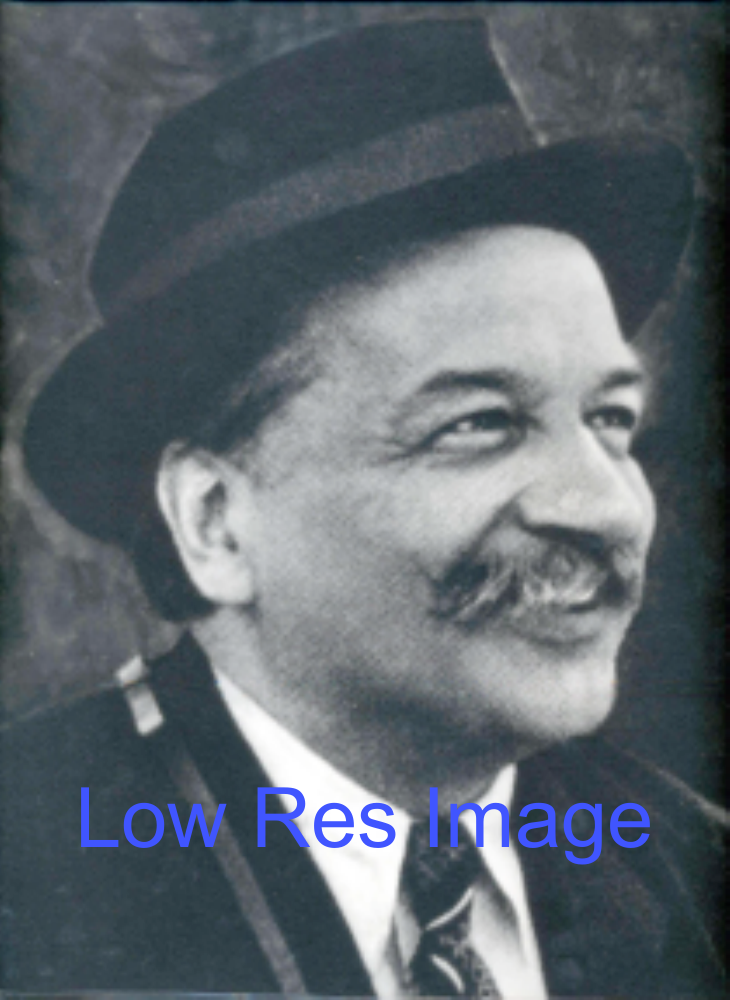}%
\caption {Tadeusz Banachiewicz, 1882--1954, at the 1946 meeting of the International Astronomical Union in Copenhagen. Courtesy of Prof.\ Adam Strza{\l}kowski of the Jagiellonian University.}
\end {minipage}
\hfil
\begin {minipage} [t] {0.5\textwidth}
\centering
\includegraphics [scale=1] {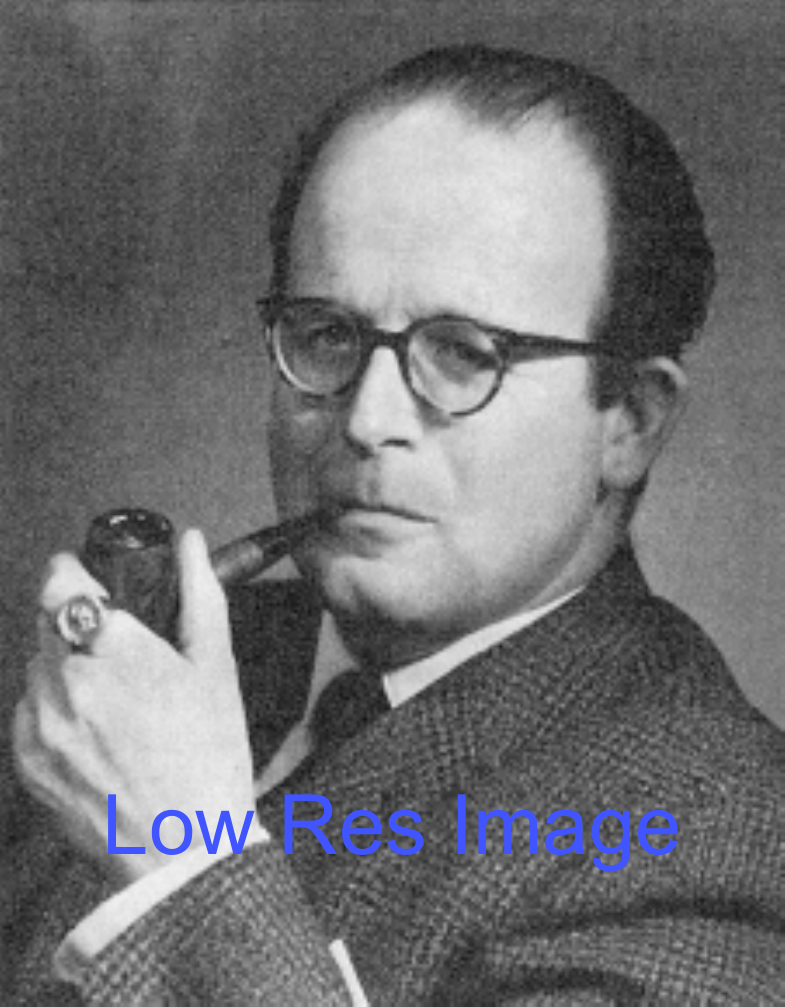}%
\caption {Henry Jensen, 1915--1974. Courtesy University of Copenhagen.}
\label {fig:Jensen}
\end {minipage}
\end {figure}

\newcommand {\squarematrix} {{\setlength {\unitlength} {1.5ex} \begin {picture} (1, 1) 
\put (0, 0) {\line (1, 0) {1}} 
\put (1, 0) {\line (0, 1) {1}} 
\put (1, 1) {\line (-1, 0) {1}}
\put (0, 1) {\line (0, -1) {1}}
\end {picture}}}

\newcommand {\lowertriangle} {{\setlength {\unitlength} {1.5ex} \begin {picture} (1, 1) 
\put (0, 0) {\line (1, 0) {1}} 
\put (0, 0) {\line (0, 1) {1}} 
\qbezier (1, 0) (0.5, 0.5) (0, 1)
\end {picture}}}

\newcommand {\uppertriangle} {{\setlength {\unitlength} {1.5ex} \begin {picture} (1, 1) 
\put (1, 0) {\line (0, 1) {1}} 
\put (0, 1) {\line (1, 0) {1}} 
\qbezier (1, 0) (0.5, 0.5) (0, 1)
\end {picture}}}

Two survey papers, by Jensen\footnote {Henry Jensen did his work on elimination at the Danish Geod{\ae}tisk Institut, and he later joined the faculty of the University of Copenhagen.} and by Bodewig, helped to establish the matrix interpretation of Gaussian elimination by describing many algorithms independent of origin in a common notation \citep {Jensen1944} \citep {Bodewig1947-HM}. \citet [3, 19] {Jensen1939} characterized his contribution as extending to ``matrix symbolism'' the results of Banachiewicz for least squares problems and for the normal equations. \citet [11] {Jensen1944} found that three algorithms for solving the normal equations were similar in that they could be interpreted as ``reducing the matrix in question to a triangular matrix:'' the ``Gauss'ian algorithm,'' the Cracovian method, and Cholesky's method.\footnote {The other, non-triangular methods that Jensen discussed were: solution by determinants, the method of equal coefficients, Boltz's method, and Kruger's method. He attributed the method of equal coefficients to B.\ I.\ Clasen in 1888; it is now more commonly attributed to Jordan as the Gauss-Jordan method, c.f.\ footnote \ref {fn:Jordan}.} To emphasize the similarities among the triangular factoring algorithms, Jensen used pictograms for triangular matrices with zeroes either \textit {under\/} \uppertriangle\ or \textit {over\/} \lowertriangle\ the main diagonal (original terminology and pictures). His primary interest was the normal equations with coefficient matrix $A^* A = N$ for a rectangular matrix $A$, where $^*$ is transposition. \citet [15, eqn.\ 15; 22, eqn.\ 3] {Jensen1944} explained that the ``Gauss'ian algorithm'' amounted to $N = \lowertriangle \, \uppertriangle \,$, where \lowertriangle\ and \uppertriangle\ are not related by transposition (he refers to Gauss's calculation for the symmetric normal equations), while in Cholesky's method the factors are so related. 

What \citet [13--16] {Jensen1944} called the ``Gauss'ian algorithm'' was his original synthesis of three different approaches. He began with his own row-oriented version of Frazer et al.'s transformation of a matrix to triangular form. He conducted the transformation symbolically by using Gauss's brackets as the matrix entries, thereby connecting the transformations to the work of Gauss. Although Jensen did not mention Doolittle, he recommended that the calculated numbers should ``be conveniently tabulated'' not in matrices but rather in a form that resembled Doolittle's table B. 

Jensen's survey was evidently read. The use in modern textbooks of the row formulation of Frazer et al.'s elementary transformations has already been noted. Similarly, Jensen's recommendation that Cholesky's method ``ought to be more generally used than is the case'' \citep [22] {Jensen1944} was heard by \citet {Fox1948} and by \citet {Turing1948} through John Todd who taught the method from Jensen's paper in 1946 \citep [197] {TausskyTodd2006}.

\subsection {Dwyer: Abbreviated Doolittle Method and Square Root Method}
\label {sec:Dwyer}

Paul Sumner Dwyer was a professor at the University of Michigan and a president of the Institute of Mathematical Statistics.\footnote {For biographical information of Dwyer, see the article by his longtime colleague Cecil \citet {Craig1972}.} Dwyer collaborated with an official at the United States Department of Agriculture named Frederick Waugh.\footnote {\citet [68--69] {Fox1989} summarizes Waugh's career at the Department.} The Department participated in creating a field of economics, known as econometrics \citep {Fox1989}, whose research methodology then consisted of data analysis by, essentially, the method of least squares. The partnership with Waugh began in the 1940s and exposed Dwyer to computers financed by government resources. This meant he could expect calculations to be mechanized, which placed him at the forefront of computing practice. 

\newcommand {\col} [1] {\parbox {3.5em} {\raggedleft #1}}
\newcommand {\xdd} {\vrule depth1.0ex height0ex width0pt}
\newcommand {\uud} {\vrule depth0.5ex height2.25ex width0pt}
\newcommand {\ud} {\vrule depth0.5ex height2.0ex width0pt}
\newcommand {\uudd} {\vrule depth1.0ex height2.25ex width0pt}
\newcommand {\udd} {\vrule depth1.0ex height2.0ex width0pt}
\newcommand {\ux} {\vrule depth0ex height1.5ex width0pt}

\begin {figure} 
\renewcommand {\baselinestretch} {1.0}
\renewcommand {\arraystretch} {1.0}
\begin {center}
\footnotesize
\setlength {\tabcolsep} {0.25em}
\begin {tabular} {r r c | r r r r c | r c | c l}
\cline {4-10}
\uudd&&&\multicolumn {1} { c} {\col {$x_1$}}&
\multicolumn {1} {c} {\col {$x_2$}}& 
\multicolumn {1} {c} {\col {$x_3$}}&
\multicolumn {1} {c} {\col {$x_4$}}&&
\multicolumn {1} {c} {\col {r.\ h.\ s.}}&\\ \cline {4-10}
&1&& \uud 1.0000& .4000& .5000& .6000&& .2000&&& original equations\\ \cline {4-10}
&2&& \ud .4000& 1.0000& .3000& .4000&& .4000&\\
&3&& \ud .5000& .3000& 1.0000& .2000&& .6000&\\
&4&& \ud .6000& .4000& .2000& 1.000&& .8000&\\ \cline {4-10}
$1 \Rightarrow$& 5&&\uud 1.0000& .4000& .5000& .6000&& .2000&&& elimination, part 1\\ \cline {4-10}
$5, 2 \Rightarrow$& 6&&\uud& .8400& .1000& .1600&& .3200&\\
$5, 3 \Rightarrow$& 7&&\ud& .1000& .7500& $-$.1000&& .5000&\\
$5, 4 \Rightarrow$& 8&&\ud& .1600& $-$.1000& .6400&& .6800&\\ \cline {4-10}
$6 \Rightarrow$& 9&&\uud& 1.0000& .1190& .1905&& .3810&&& elimination, part 2\\ \cline {4-10}
$9, 7 \Rightarrow$& 10&&\uud&& .7381& $-$.1190&& .4619&\\
$9, 8 \Rightarrow$& 11&&\ud&& $-$.1190& .6095&& .6190&\\ \cline {4-10}
$10 \Rightarrow$& 12&&\uud&& 1.0000& $-$.1612&& .6258&&& elimination, part 3\\ \cline {4-10}
$12, 11 \Rightarrow$& 13&&\uud&&& .5903&& .6935&\\ \cline {4-10}
$13 \Rightarrow$& 14&&\uud&&& 1.0000&& 1.1748&&& back-substitution\\
$11, 14 \Rightarrow$& 15&&\ud&& 1.0000&&& .8152&\\
$9, 14, 15 \Rightarrow$& 16&&\ud& 1.0000&&&& .0602&\\
$5, 14, 15, 16 \Rightarrow$& 17&&\ud1.0000&&&&& $-$.9366&\\ \cline {4-10}
\end {tabular}
\end {center}
\caption {A form of Gaussian elimination that Paul \citet {Dwyer1941a} called ``the method of single division'' and which he found equivalent, except for cosmetic changes, to the ``method of pivotal condensation'' of \citet {Aitken1937}, and to an earlier method of \citet {Deming1928}. The nonunitary numbers in rows 5, 9, and 12 are the upper diagonal entries in Crout's table.}
\label {fig:singledivision}
\end {figure}

\citet {Dwyer1941a} began from a comparison of solution methods that emphasized sources in American and English statistics. He considered 22 papers from 1927 to 1939 and suggested their bibliographies should be consulted for an even more thorough picture of the subject. Dwyer's review painstakingly uncovered the similarities between proliferating methods distinguished by minor changes to the placement of numbers in tables. For example, he noted a method of \citet {Deming1928} was equivalent, except for superficial changes, to a ``method of pivotal condensation'' of \citet {Aitken1937}, both of which he included under the rubric ``method of single division'' (Figure \ref {fig:singledivision}). The differences among these methods seem pedestrian until one has to choose the most effective way to calculate and record all the numbers in the tables by hand. Nevertheless, in another paper, \citet [260] {Waugh1945} summarized the field more succinctly, observing that the methods are more or less the same, except ``Crout divides the elements of each row by the leading element while we divide the elements of columns.'' \citet [111-112] {Dwyer1941a} credited to others, including Waugh, the observation that accumulating calculators made it unnecessary to record the series terms in Doolittle's table B. Dwyer called the streamlined procedure the ``abbreviated method of single division -- symmetric'' or the ``abbreviated Doolittle method.''

\citet {Dwyer1944} independently interpreted Gaussian elimination as matrix factoring. His primary interest was the case \ref {case:one} least squares problem described in section \ref {sec:squares}. Beginning from the coefficient matrix $A$ of the normal equations, he showed the abbreviated Doolittle method was an ``efficient way of building up'' some ``so called triangular'' matrices $S$ and $T$ with $A - S^t T = 0$ \citep [86] {Dwyer1944}. Rather than as product of the elementary matrices of \citet {Frazer1938} and \citet {Jensen1944}, Dwyer obtained the decomposition from a sum of vector outer products \citep [86, eqns.\ 21, 23] {Dwyer1944}. He remarked that this formula would be the ``key'' to ``a more general theory'' \citeibid [88] {Dwyer1944} of decompositions which \citet {Waugh1945} subsequently developed to invert nonsymmetric matrices. For the normal equations it was possible to choose $S = T$ for a ``square root'' method \citep [88] {Dwyer1944} (modern Cholesky method) which \citet {Dwyer1945} developed in a later paper. He added that he found no other matrix algebra interpretation of solving equations except from \citet {Banachiewicz1942} who, he noted, also had a square root method \citep [89] {Dwyer1944}. 

Dwyer particularly influenced computers in the United States. \citet {Laderman1948} reported that \citet {Duncan1946} popularized Dwyer's square root method, and that it was even used at the Mathematical Tables Project. European mathematicians such as \citet {Fox1954}, however, preferred to apply Cholesky's name. Since Dwyer's many papers and his book on linear equations (1951) \nocite {Dwyer1951}always invoked the memory of Doolittle, Dwyer's own name was never attached to either of the computing methods that he championed. 

\subsection {Satterthwaite: Three-factor Decomposition}

\citet {Dwyer1941a} had concluded his survey of methods with a process that he attributed to R.\ A.\ Fisher for solving equations with multiple right sides. The presentation inspired the statistician Franklin \citet {Satterthwaite1944} to calculate inverse matrices.\footnote {\citet {Satterthwaite1944} published between \citet {Dwyer1944} and \citet {Waugh1945} and in the same journal, but he cited only Dwyer's papers from 1941.} He linked the tables of the calculation to inverses through the triple ``factorization'' $A = (R_1 + I) \, S_1 (I + T_1)$ where $S_1$ is diagonal, and $R_1$ and $T_1$ are respectively ``pre-'' and ``postdiagonal'' (strictly triangular) \citep [374, eqn.\ 3.1; 376, sec.\ 5] {Satterthwaite1944}. In consideration of having to round the numbers that were written into the tables, Satterthwaite showed the process was accurate provided $A$ was close to being an identity matrix. He therefore suggested using the process to improve inverses: if $F \approx A^{-1}$, then  calculate $(FA)^{-1} F$. His paper appears to have had little influence. It received only ten citations of which four were by Dwyer and Waugh.

\subsection {Bodewig: Literature Surveys}

\begin {figure} [t]
\centering
\begin {minipage} [t] {0.45\textwidth}
\centering
\includegraphics [scale=1] {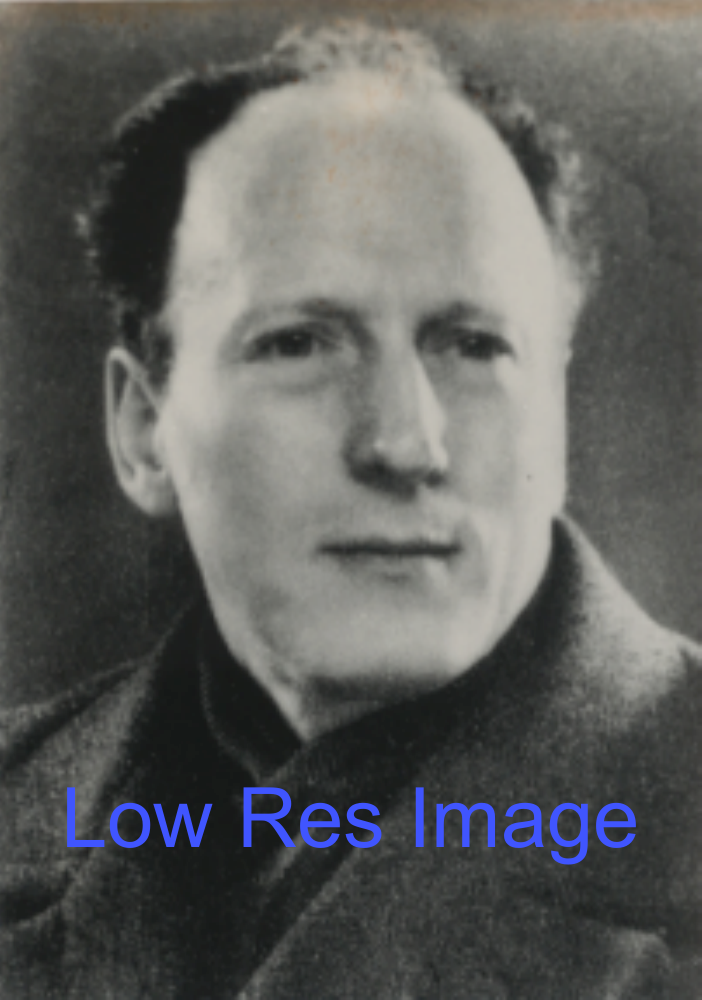}%
\caption {Ewald Konrad Bodewig. Courtesy Archives of the Mathematisches Forschungsinstitut Oberwolfach.}
\label {fig:Bodewig}
\end {minipage}
\hfil
\begin {minipage} [t] {0.50\textwidth}
\centering
\includegraphics [scale=1] {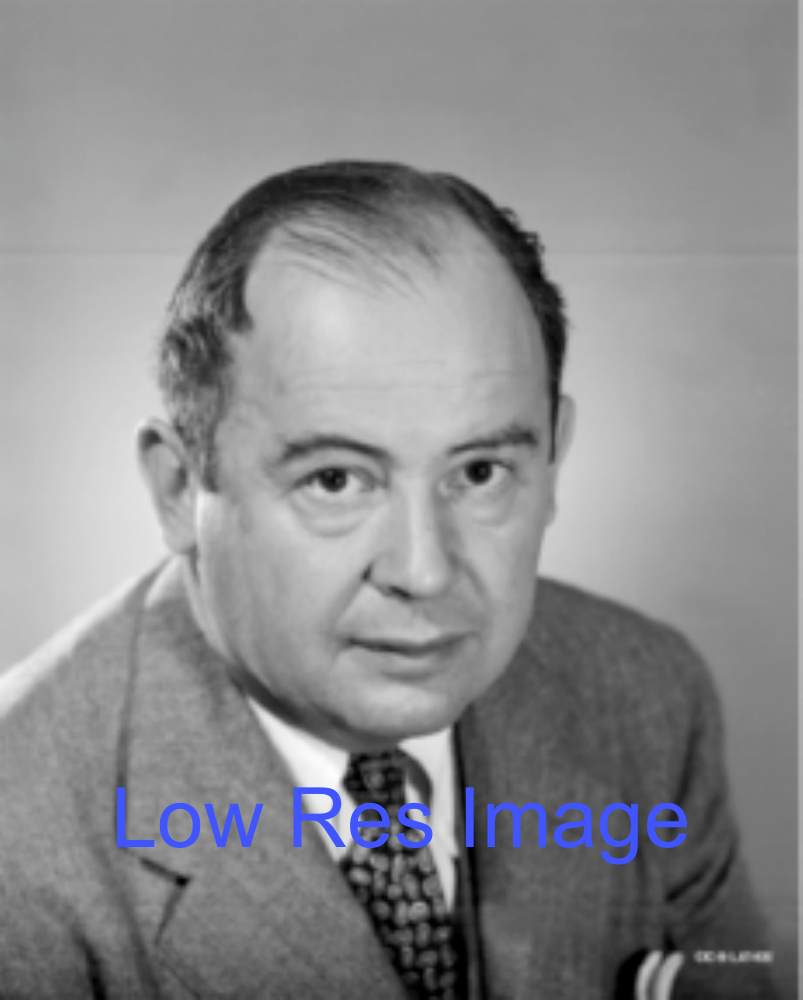}%
\caption {John von Neumann, 1903--1957, circa March, 1947. Courtesy of the Los Alamos National Laboratory Archives.}
\label {fig:vonNeumann}
\end {minipage}
\end {figure}

The mathematician Ewald Bodewig followed Jensen's approach on a grand scale in a five-part paper that summarized in matrix notation all methods for solving linear equations that were known up until 1947 \citep {Bodewig1947-HM}.\footnote {During the same period \citet {Bodewig1944-HM, Bodewig1948-HM, Bodewig1949-HM} wrote evaluations for \textit {Mathematical Reviews\/} of all three of the papers by \citet {Jensen1944}, \citet {vonNeumann1947-HM}, and \citet {Turing1948}. Shortly thereafter Bodewig became distraught over the paltry wages paid to mathematicians and announced his semi-retirement from mathematics in a letter which received much comment \citep {Bodewig1950, Betts1951}.} His interesting bibliography lay at the end of the paper but, like Jensen's, it was not comprehensive because he neglected authors such as Crout, Doolittle, and Dwyer. He took Jensen too literally because he labelled triangular matrices ``left'' or ``right'' in order to indicate the location of the \textit {zeroes\/} instead of the location of the \textit {non\/}zeroes: Bodewig's ``right'' was $\lowertriangle = {\mathfrak D}_r$ and ``left'' was $\uppertriangle = {\mathfrak D}_l$ where $\mathfrak D$ stood for Dreiecksmatrix \citeibid [part I, 444] {Bodewig1947-HM}. He repeated Jensen's row version of Frazer et al.'s presentation of Gaussian elimination, and he emphasized their summary formula using Jensen's pictograms, $\, \lowertriangle \; \squarematrix = \uppertriangle \,$ \citeibid [part I, 449] {Bodewig1947-HM}. He also followed \citet {Jensen1944} in describing Cholesky's method as $\mathfrak S = \mathfrak D_l^\prime \mathfrak D_l^{}$ for a symmetric $\mathfrak S$, where $^\prime$ is transpose \citeibid [part I, 450] {Bodewig1947-HM}.

\subsection {Von Neumann and Goldstine: The Combination of Two Tricks}
\label {sec:vonNeumann}

Among the first authors to describe Gaussian elimination in terms of matrix algebra, John von Neumann and his collaborator Herman Goldstine were alone in making a nontrivial use of the matrix decomposition. They devised a decomposition method related to Gaussian elimination that was appropriate for use with the first electronic computers which were then being built. They established bounds on the rounding errors of the computation in terms of the ratio of the largest to the smallest singular values of the coefficient matrix. Their result, which is beyond the scope of the present discussion, marks the beginning of modern research interests in computational mathematics. The ratio that they discovered is now called the matrix condition number. 

Von Neumann and Goldstine were also the only authors to show exactly how the ``traditional schoolbook method'' \citep [3] {Farebrother1988} of equation (\ref {eqn:Cohen}) calculates triangular factors of the coefficient matrix. They described the elimination algorithm as creating a sequence of ever-smaller reduced matrices and vectors from $A x = y$ \citep [1051] {vonNeumann1947-HM}, 
\begin {displaymath}
A = A^{(1)}, \; A^{(2)}, \; A^{(3)}, \; \ldots , \; A^{(n)} ,
\qquad
y = y^{(1)}, \; y^{(2)}, \; y^{(3)}, \; \ldots , \; y^{(n)} , 
\end {displaymath}
where the rows and columns of $A^{(i)}$ and $y^{(i)}$ are numbered from $i$ to $n$.\footnote {Is customary to neglect the possibility of reordering the equations and variables in presentations such as this to avoid complicating the notation.} For $i = 1, \ldots, n-1$ the computation is,
\begin {eqnarray}
\label {eqn:LU}
A^{(i+1)}_{j,k}& =& \makebox [9em] [l] {$\displaystyle A^{(i)}_{j,k} -
A^{(i)}_{j,i} A^{(i)}_{i,k} / A^{(i)}_{i,i} $}
\qquad 
\mbox {for $j$, $k > i$}
\\
\noalign {\smallskip}
\label {eqn:y}
y^{(i+1)}_{j}& =& \makebox [9em] [l] {$\displaystyle y^{(i)}_{j} -
(A^{(i)}_{j,i} / A^{(i)}_{i,i} ) \, y^{(i)}_{i}$}
\qquad 
\mbox {for $j > i$} \, .
\end {eqnarray}
Next, the algorithm solves by substitution the equations $B^\prime x = z$ where the entries of $B^\prime$ and $z$ are chosen from the reduced matrices and vectors (the first row of each). For a matrix $C$ of the multipliers $\smash {A^{(i)}_{j,i} / A^{(i)}_{i,i}}$ with $j \ge i$ (note the unit diagonal), von Neumann and Goldstine summed equation (\ref {eqn:y}) over $i$ and rearranged to give $C z = y$. From this equation and $B^\prime x = z$ they concluded $C B^\prime = A$. 
\begin {quotation}
\noindent We may therefore interpret the elimination method as \dots\ the combination of two tricks: First, it decomposes $A$ into a product of two semi-diagonal matrices \dots\ [and second] it forms their inverses by a simple, explicit, inductive process. 
\next --- \citet [1053] {vonNeumann1947-HM}\footnote {Von Neumann and Goldstine wrote ``semi-diagonal'' for ``triangular.''}
\end {quotation}
With this passage von Neumann and Goldstine finally expressed in a new idiom what Newton and Rolle had first explained over two hundred years earlier.
\begin {quotation}
\noindent 
And you are to know, that by each {\AE}quation one unknown Quantity may be taken away, and consequently, when there are as many {\AE}quations and unknown Quantities, all at length may be reduc'd into one, in which there shall be only one Quantity unknown. 
\same \mbox {--- \citet [60--61] {Newton1720-HM} and prior years}

\end {quotation}

Von Neumann and Goldstine found a lack of symmetry in what is now called the ``classic'' algorithm because the first factor always had $1$'s on its main diagonal. They divided the second factor by its diagonal to obtain $B^\prime = DB$, hence $A = CDB$ which they said was ``a new variant'' of Gaussian elimination \citep [1031] {vonNeumann1947-HM} that is now written $A=LDU$,
\begin {equation}
\label {eqn:LandDandU}
L_{j,i} = A^{(i)}_{j,i} / A^{(i)}_{i,i}
\qquad
D_{i,i} = A^{(i)}_{i,i}
\qquad
U_{i,k} = A^{(i)}_{i,k} / A^{(i)}_{i,i}
\qquad
j,k \ge i \, ,
\end {equation}
where $\smash {A^{(i)}_{j,k}}$ are the entries of the reduced matrices given by equation (\ref {eqn:LU}). 

\subsection {Decompositions and Diagonals}

In summary, four groups of authors took conceptually different routes to interpret Gaussian elimination as constructing triangular decompositions in Cayleyan matrix algebra. Inspired by Banachiewicz, \citet {Jensen1944} obtained a triangular decomposition from products of the elementary matrices of \citet {Frazer1938}. \citet {Dwyer1944} expressed a matrix as a sum of outer vector products from which a triangular decomposition followed. \citet {Satterthwaite1944} found entries of the factors in tables used to solve equations with multiple right hand sides. \Citet {vonNeumann1947-HM} chose entries of the triangular factors directly from numbers that occurred in the ``schoolbook'' or symbolic algebra form of Gaussian elimination. 

What accounts in part for the existence of multiple inventors of elimination algorithms is the fact that the triangular matrix decompositions are not unique. The triple term decomposition, $A=LDU$, is unique only from the requirement that all the diagonal entries of $L$ and $U$ equal $1$. The decompositions $A = L(DU)$ and $A = (LD)U$ that differently apportion the diagonal are known today by the names Doolittle and Crout, respectively. If $A$ is a symmetric and positive definite matrix, then $L = U$ and the decomposition $A = (L D^{1/2}) (L D^{1/2})^t$ is known by the name Cholesky.\footnote {Similarly, what Gauss computed might be characterized as $A = (LD) D^{-1} (LD)^t$ because he only retained his bracket values which give the entries of $(LD)$ and those of $D$ on the diagonal.} The Doolittle naming convention is misleading because the calculations of Doolittle --- and Gauss also, of course --- applied only to symmetric, positive definite matrices. All these names are anachronisms because the inventors did not use matrix notation. 

\section {Conclusion}
\label {sec:part6}

\subsection {Summary}

An algorithm equivalent to Gaussian elimination appeared in mathematical texts from ancient China. A few instances of the method emerged during the late Renaissance in European textbooks, but since they were geographically scattered, it seems unlikely they were influenced by contact with the orient. \citet {Newton1720-HM}, writing circa 1670, and \citet {Rolle1690-HM} separately gave rules for successively eliminating variables to solve simultaneous equations, and they chose words other than \textit {elimination\/} to describe their work. The features of Newton's rule can be traced through successive algebra texts from many authors resulting in a fairly standardized presentation by 1800. \citet {Lacroix5th1804-HM} emphatically called this schoolbook process ``elimination.''

At the turn of the century both Legendre and Gauss remarked the textbook rule could be used to solve least squares problems that were formulated as normal equations. Not content with that, \citet {Gauss1810-HM} devised a special solution process. His approach was adopted by professional hand computers working in geodesy, and was then replaced (Doolittle, Cholesky) to make better use of computing aids including manual calculators. Astronomers, geodesists, and statisticians (Banachiewicz, Jensen, Dwyer) endeavored to improve least squares calculations through studies based on matrix algebra. Nothing entirely new resulted, although Dwyer's independently discovered ``square root'' (Cholesky's) method was widely adopted in the United States.

These methods, including that of Gauss, cannot substitute for schoolbook elimination if other than normal equations are of interest. A need to solve more general equations arose in engineering in the \xx\ century. \citet {Frazer1938} proposed a formalism for calculations in terms of ``elementary'' matrix operations, while \citet {Crout1941a}, without matrices, devised a version of the schoolbook method that was publicized by a calculator manufacturer. \citet {Dwyer1944} and \citet {Satterthwaite1944} observed that any matrix could be decomposed into triangular factors, which the latter and \citet {Waugh1945} exploited to calculate inverse matrices. \Citet {vonNeumann1947-HM} finally interpreted schoolbook elimination as a matrix decomposition for the purpose of studying the accuracy of calculations made by electronic computers. Gaussian elimination subsequently became the subject of a large literature patterned after their study.

\subsection {Remembering Gauss Imperfectly}

At the beginning of the \xix\ century, the triumph of calculating the location of Ceres had earned the youthful Gauss fame enough to realize his wish for a life free of teaching mathematics. Nevertheless, by the end of that century, histories of mathematics would neglect Gauss's applied work.\footnote {\citet [46] {Buhler1981} and \citet [405--410] {Dunnington2004-HM} explain Gauss's career choices; \citet {Matthiessen1878} and \citet {Cajori1894-HM} are examples of early histories of mathematics.} \citet [177] {GrattanGuinness2004} explains: the establishment of the research university as a social institution, which occurred roughly during Gauss's lifetime \citep {Clark2006}, increased employment for credentialed mathematicians, which both sped the professionalization of the subject and coincided with a preference for pure over applied mathematics. The same taste was manifest in historical scholarship, he adds, in that pure subjects became more thoroughly chronicled than applications. 

Grattan-Guinness argues for the existence of two recollections of the mathematical past: history recounts the development of ideas in the context of contemporary associations, while heritage remembers reinterpreted work that embodies the state of mathematical knowledge. The example of Gaussian elimination cautions that much of applicable mathematics may lack both: no mathematical history because applications may be deemed the purview of non-mathematical faculties, no heritage because applications might not enter or remain in the corpus. As a result, significant developments may come from what appears to be the periphery of the mathematical community. The interpretation of Gaussian elimination through matrix algebra, which is the most important development in understanding the subject for present purposes, arose as much outside as inside mathematics. 

Gauss's methods to determine orbits and to solve inference problems indelibly linked calculations to his reputation among astronomers and geodesists. Whereas in the first half of the \xix\ century algebra textbooks referred only to elimination, from the second half of the \xix\ century reference books for astronomy and geodesy always cited Gauss in order to recommend their least squares calculations. His bracket notation for what came to be called auxiliaries was correctly regarded as pertinent to the solution of just normal equations. For example, \citet [530] {Chauvenet1868-HM} had the ``elimination of unknown quantities from the normal equations \dots\ according to Gauss,'' and \citet [557] {Liagre1879-HM} ``l'\textit {\'elimination\/} des inconnues entre les \'equations du minimum (\'equations normales)'' by ``les coefficients auxiliaires de Gauss.'' 

\begin {table} 
\renewcommand {\baselinestretch} {1.0}
\renewcommand {\arraystretch} {1.0}
\caption {Nomenclature for solving simultaneous linear equations by elimination indicating the gradual evolution to ``Gaussian elimination'' over many years. Some of these excerpts contain opinions about the history of the subject that are not necessarily accurate.} 
\label {tab:nomenclature}
\begin {center}
\scriptsize
\setlength {\tabcolsep} {0.5em}
\renewcommand {\arraystretch} {1.25}
\renewcommand {\arraystretch} {1.00}
\newcommand {\before} {\vrule depth0pt height2.25ex width0pt}
\newcommand {\after} {\vrule depth1.0ex height0ex width0pt}
\renewcommand {\arraystretch} {1.00}
\newcommand {\firstbefore} {\vrule depth0pt height2.10ex width0pt}
\renewcommand {\before} {\vrule depth0pt height1.90ex width0pt}
\renewcommand {\after} {\vrule depth0.80ex height0ex width0pt}
\newcommand {\lastafter} {\vrule depth0.90ex height0ex width0pt}
\newcommand {\leftside} [1] {\begin {minipage} [t] {13.5em} \raggedleft #1\after \end {minipage}}
\newcommand {\rightside} [1] { \begin {minipage} [t] {15.5em} \raggedright #1\after \end {minipage}}
\newcommand {\names} [1] {\begin {minipage} [t] {6em} \before \raggedright #1\after \end {minipage}}
\begin {tabular} {| c | l | r | c | r l |}
\hline
\sc year
& \multicolumn {1} {c |} {\before \sc author \after}
& \multicolumn {1} {c |} {\before \sc page \after}
& \multicolumn {1} {l } {{\sc uses gauss's brackets} \hspace*{-10em}}
& \multicolumn {2} {c |} {\sc nomenclature}
\\ \hline \hline
1690& \firstbefore Rolle \nocite {Rolle1690-HM}\after& 34&& Des Substitutions& \\
& \before \nocite {Rolle1690-HM}\after& 42&& De la Methode& \\
1720& \before Newton \nocite {Newton1720-HM}\after& 60&& Transformation& of two or more Equations into one
\\
\multicolumn {2} {| l |} {and prior years}&&& \before in order to exterminate\after& the unknown Quantities
\\
1742& \before Hammond \nocite {Hammond1742}\after& 142&& is called the \textit {exterminating\/}&
[of] \textit {an unknown Quantity\/}
\\
1755 & \before Simpson \nocite {Simpson1755-HM}\after& 63&& \leftside {Of the Extermination\\ the reduction}& \rightside {of unknown quantities, or\\ of two, or more equations, to a single}
\\
1771& \before Euler \nocite {Euler1771-HM}\after& \multicolumn {2} {l} {part 2 page 27 \hspace*{-7em}}& Der nat\"urlichste Weg& bestehet nun darinn
\\
1804& \before Lacroix \nocite {Lacroix5th1804-HM}\after& 114&& se nomme \textit {\'elimination\/}&
\\
1805& \before Legendre \nocite {Legendre1805-HM}\after& 73&& par les m\'ethodes&
ordinaires
\\
1809& \before Gauss \nocite {Gauss1809-HM}\after& 214&& per eliminationem& vulgarem
\\
1818& \before Lacroix \nocite {Lacroix1818-HM}\after& 84&& is called \textit {elimination\/}&
\\
1822& \before Euler \nocite {Euler1822-HM}\after& 216&& the most natural method&
of proceeding
\\
1828& \before Bourdon \nocite {Bourdon1828-HM}\after& 79&& La m\'ethode d'elimination&
\\
1835& \before Davies \nocite {Davies1835-HM}\after& 93&& quantities may be eliminated& by the following rule
\\
1846& \before Clark \nocite {Clark1846-HM}\lastafter& 102&& Of elimination& when there are three or more \dots
\\
& \before \nocite {Clark1846-HM}\lastafter& 103&& Reduce& the equations\\ \hline
1868& \firstbefore Chauvenet \nocite {Chauvenet1868-HM}\after& 530& \checkmark& \leftside {Elimination\\ the method of substitution,}& \rightside {of the unknown quantities by\\ according to Gauss}
\\
1879& \before Liagre \nocite {Liagre1879-HM}\after& 557& \checkmark& Pour
proc\'eder a l'\textit {\'elimination\/}& \dots\ coefficients auxiliaires de Gauss
\\
1884& \before Merriman \nocite {Merriman1884-HM}\after& 51&& the method of
substitution& \begin {minipage} [t] {15em} \raggedright due to Gauss by which is preserved the symmetry\after \end {minipage}

\\
1888& \names {Doolittle \nocite {Doolittle1888-HM} (not \textit {the\/} Doolittle)}& 43& \checkmark& \leftside {In the elimination\\ the method of substitution}& \rightside {it will be convenient to employ\\ using \dots\ notation
proposed by Gauss.}
\\
1895& \before Jordan \nocite {Jordan1895-HM}\after& 80& \checkmark& \textbf {Gauss schen Elimination}& 
\\
1900& \before Bartlett \nocite {Bartlett1900-HM}\after& 97& \checkmark& ``Method of
Substitution''& proposed by Gauss
\\
1905& \before Johnson \nocite {Johnson1905-HM}\after& 120& \checkmark& method of
substitution& \rightside {as developed by Gauss, which has the advantage of preserving, \dots\ in the elimination, the symmetry}
\\
1906& \names {Wright and Hayford \nocite {Wright1906-HM}}& 106& \checkmark& method of substitution& \rightside {introduced by Gauss\\ and the Doolittle method}\\
1907& \before Helmert \nocite {Helmert1907-HM}\after& 120& \checkmark& Algorithmus& von C.\ F.\ Gauss
\\
1912& \before Palmer \nocite {Palmer1912-HM}\after& 92& \checkmark& Gauss's
method& of solution with logarithms
\\
1924& \before Benoit \nocite {Benoit1924-HM}\after& 24&& les m\'ethodes
ordinaires& y compris celle de Gauss
\\
1927& \names {Tolley and Ezekiel \nocite {Tolley1927}}& 107&& \leftside {Gauss himself developed a\kern1em \\ direct process of elimination}& \rightside {\ \\[0.15ex] which was much shorter than\\ solution by determinants\lastafter}
\\ \hline
1941& \firstbefore Crout \nocite {Crout1941a}\after& 1235&& Gauss's method&
\\
1941& \before Dwyer \nocite {Dwyer1941b}\after& 450& \multicolumn {1} {l}
{mentions \hspace*{-7em}}& notation& \dots\ suggested by Gauss
\\
1943& \before Hotelling \nocite {Hotelling1943}\after& 3&& Doolittle method&
\\
1944& \before Jensen \nocite {Jensen1944}\after& 13& \checkmark& Gauss'ian algorithm&
\\
1947& \before Bodewig \nocite {Bodewig1947-HM}\after& 447&& Gaussschen Methode&
\\
1947& \names {von Neumann and Goldstine \nocite {vonNeumann1947-HM}}& 1049&& the conventional elimination& method
\\
1948& \names {Fox, Huskey, and Wilkinson\nocite {Fox1948}}& 149&& Gaussian algorithm&
\\
1948& \before Turing \nocite {Turing1948}\lastafter& 287&& ``Gauss elimination&
process''
\\ \hline
1953& \firstbefore Forsythe \nocite {Forsythe1953}\after& 301&& systematic elimination& \rightside {of unknowns in the fashion of high- school algebra, as described by Gauss}
\\
& \before \nocite {Forsythe1953}\after& 315&& \textbf {Gaussian elimination}& (pivotal condensation)
\\
1953& \before Householder \nocite {Householder1953}\after& 68&& method of elimination& by addition and subtraction
\\ 
1960& \before Wilkinson \nocite {Wilkinson1960}\after& 45&& \textbf {Gaussian elimination}&
\\
1964& \before Householder \nocite {Householder1964}\after& 253&& \textbf {Gaussian elimination}&
\\
1977& \before Goldstine \nocite {Goldstine1977}\lastafter& 287&& \leftside {he [Gauss] discovered the\kern1em\\ method of \textbf {Gaussian elimination}}& \rightside {\ \\ for solving systems of linear equations}
\\ \hline
\end {tabular}
\end {center}
\vspace*{-5ex}
\end {table}

The invention of electronic computers led to new university faculties for instruction in the new way of computing which drew people of diverse training to the field. They extrapolated opinions of its history from whatever heritage they knew. Geodesy at one time accounted for the bulk of the simultaneous linear equations that were solved professionally, so its terminology may have appeared to be authoritative. The citations that precisely described the contribution of Gauss had already been shortened to an unspecific ``Gauss's procedure.'' This usage was misinterpreted as attributing ordinary, common, schoolbook elimination to Gauss, but only after World War II (Table \ref {tab:nomenclature}). Von Neumann [1947] was apparently the last prominent mathematician to simply write ``elimination'' as \citet {Lacroix5th1804-HM} and \citet {Gauss1809-HM} had done. There is no harm in naming elimination after Gauss because he did initiate the development of professional methods which are of considerable importance today, but it should be understood that ``Gaussian elimination'' is an honorific and not an attribution for a simple concept to which surprisingly many people have contributed.

\section* {Acknowledgments} 



I am grateful to the referees and especially to the editor J.\ Barrow-Green for comments and corrections that much improved this paper. I also thank C.\ Brezinski,  A.\ B.\ Carr, I.\ Ciolli, F.\ J.\ Conahan, R.W.\ Farebrother, C.\ F{\"o}rstel, A.\ George, J.\ J.\ Harbster, J.\ H{\o}yrup, L.\ Knoblauch, J.\ Kocinski, J.\ Marzahn, V.\ Monnier, M.\ R.\ Mutter, I.\ Olkin, R.\ Otnes, N.-O.\ Pr{\ae}gel, F.\ Pattman, K.\ Plofker, E.\ Robson, S.\ Sanders, E.\ Slottved L.\ E.\ Sj\"oberg, S.\ Stigler, A.\ Strza{\l}kowski, E.\ Sullivan, E.\ Swartzlander, G.\ Tee, A.\ E.\ Theberg, M.-C.\ Thooris, J.\ Todd, and I.\ Vetter who generously helped with photographs or references. Finally, I commend many institutions and Google books for making most of the following older sources available online. 


{
\raggedright

\begin{thebibliography}{212}
\expandafter\ifx\csname natexlab\endcsname\relax\def\natexlab#1{#1}\fi
\expandafter\ifx\csname url\endcsname\relax
  \def\url#1{\texttt{#1}}\fi
\expandafter\ifx\csname urlprefix\endcsname\relax\def\urlprefix{URL }\fi

\bibitem[{Ahlmann and Rosenbaum(1933)}]{AhlmannRosenbaum1933}
Ahlmann, H.~W., Rosenbaum, L., 1933. {Scientific Results of the
  Swedish-Norwegian Arctic Expedition in the Summer of 1931, Part I-III}.
  Geografiska Annaler 15, 1--68.

\bibitem[{Aitken(1937)}]{Aitken1937}
Aitken, A.~C., 1937. {Studies in practical mathematics I. The evaluation, with
  application, of a certain triple product matrix}. Proc. Royal Soc., Edinburgh
  57, 172--181.

\bibitem[{Aitken(1951)}]{Aitken1951}
Aitken, A.~C., 1951. Determinants and Matrices, seventh Edition. Oliver and
  Boyd, Edinburgh.

\bibitem[{Althoen and McLaughlin(1987)}]{Althoen1987}
Althoen, S.~C., McLaughlin, R., 1987. {Gauss-Jordan Reduction: A Brief
  History}. American Mathematical Monthly 94, 130--142.

\bibitem[{Apokin(2001)}]{Apokin2001}
Apokin, I.~A., 2001. {Willgodt Theophil Odhner and his ``Arifmometr:'' The
  Beginning of the Mass-Production of Calculating Machines}. In: Trogemann, G.,
  Nitussov, A.~Y., Ernst, W. (Eds.), Computing in Russia. Vieweg,
  Braunschweig/Wiesbaden, pp. 39--46.

\bibitem[{Baker(1903)}]{Baker1903}
Baker, H.~F., 1903. On the integration of linear differential equations. Proc.
  London Math. Soc. 35, 333--378.

\bibitem[{Banachiewicz(1937{\natexlab{a}})}]{Banachiewicz1937b}
Banachiewicz, T., 1937{\natexlab{a}}. Sur la r\'esolution num\'erique d'un
  syst\`eme d'\'equations lin\'eaires. Bull. intern. de l'Acad. Polonaise,
  S\'erie A. Sc. Math., 350.

\bibitem[{Banachiewicz(1937{\natexlab{b}})}]{Banachiewicz1937a}
Banachiewicz, T., 1937{\natexlab{b}}. {Zur Berechnung der Determinanten wie
  auch der Inversen, und zu darauf basierten Aufl\"osung linearer Gleichungen}.
  Acta Astronomica 3, 41--67.

\bibitem[{Banachiewicz(1938{\natexlab{a}})}]{Banachiewicz1938c}
Banachiewicz, T., 1938{\natexlab{a}}. M\'ethode de r\'esolution num\'erique des
  \'equations lin\'eaires, du calcul des d\'eterminants et des inverses, et de
  r\'eduction des formes quadratiques. Bull. intern. de l'Acad. Polonaise,
  S\'erie A. Sc. Math., 393--404.

\bibitem[{Banachiewicz(1938{\natexlab{b}})}]{Banachiewicz1938a}
Banachiewicz, T., 1938{\natexlab{b}}. Principes d'une nouvelle technique de la
  m\'ethode des moindres carr\'es. Bull. intern. de l'Acad. Polonaise, S\'erie
  A. Sc. Math., 134--135.

\bibitem[{Banachiewicz(1942)}]{Banachiewicz1942}
Banachiewicz, T., 1942. {An outline of the Cracovian algorithm of the method of
  least squares}. Astr. Jour. 50, 38--41.

\bibitem[{Bargmann et~al.(1946)Bargmann, Montgomery, and von
  Neumann}]{Bargmann1946-HM}
Bargmann, V., Montgomery, D., von Neumann, J., 1946. Solution of linear systems
  of high order. In: \rm \cite {Taub1963}. Vol.~5. pp. 421--477.

\bibitem[{Bartlett(1900)}]{Bartlett1900-HM}
Bartlett, D.~P., 1900. General Principles of the Method of Least Squares with
  Applications, 2nd Edition. Private printing, Boston.

\bibitem[{Bashmakova and Smirnova(2000)}]{Bashmakova2000-HM}
Bashmakova, I.~G., Smirnova, G.~S., 2000. The Beginnings and Evolution of
  Algebra (Tr.\ A.\ Shenitzer). Vol.~23 of Dolciani Mathematical Expositions.
  Mathematical Association of America, Washington.

\bibitem[{Baxendall(1929)}]{Baxendall1929}
Baxendall, D., 1929. Calculating machines. In: Encyclopedia Britannica, 14th
  Edition. Vol.~4. pp. 548--553.

\bibitem[{Benoit(1924)}]{Benoit1924-HM}
Benoit, April--June 1924. {Note sur une m\'ethode de r\'esolution des
  \'equations normales provenant de l'application de la m\'ethode des moindres
  carr\'es a un syst\`eme d'\'equations lin\'eaires en nombre inf\'erieur a
  celui des inconnues. --- Application de la m\'ethode a la r\'esolution d'un
  syst\`eme defini d'\'equations lin\'eaires (Proc\'ed\'e du Commandant
  Cholesky)}. Bulletin g\'eod\'esique~(2), 67--77. {T}he author is identified
  only as {C}ommandant {B}enoit.

\bibitem[{Bernkopf(1968)}]{Bernkopf1968}
Bernkopf, M., Jan. 1968. A history of infinite matrices. Archive for History of
  Exact Sciences 4~(4), 308--358.

\bibitem[{Bessel and Baeyer(1838)}]{Bessel1838-HM}
Bessel, F.~W., Baeyer, J.~J., 1838. Gradmessung in Ostpreussen und ihre
  Verbindung mit Preussischen und Russischen Dreiecksketten. K\"oniglichen
  Akademie der Wissenschaften, Berlin.

\bibitem[{Betts et~al.(1951)Betts, Bauer, Tiedjens, Landsberg, and
  Holter}]{Betts1951}
Betts, G.~L., Bauer, D., Tiedjens, V.~A., Landsberg, H.~E., Holter, N.~J.,
  1951. Scholars and the root of all evil. Science 113, 330--333.

\bibitem[{B\'ezout(1788)}]{Bezout1788-HM}
B\'ezout, E., 1788. Cours de math\'ematiques, \`a l'usage du corps royal de
  l'artillerie. Vol. 2 of 4. Ph.-D.\ Pierres, Paris.

\bibitem[{B\'ezout(2006)}]{Bezout2006-HM}
B\'ezout, E., 2006. General theory of algebraic equations (Tr.\ E.\ Feron).
  Princeton University Press, Princeton. First French edition 1779.

\bibitem[{Bishop(1987)}]{Bishop1987}
Bishop, R. E.~D., 1987. {Arthur Roderick Collar. 22 February 1908 -- 12
  February 1986}. Biographical Memoirs of Fellows of the Royal Society 33,
  164--185.

\bibitem[{Black(1949)}]{Black1949}
Black, A.~N., 1949. Further notes on the solution of algebraic linear
  simultaneous equations. The Quarterly Journal of Mechanics and Applied
  Mathematics 2~(3), 321--324.

\bibitem[{Bodewig(1944)}]{Bodewig1944-HM}
Bodewig, E., 1944. Mathematical Reviews 7, 488. Review of \citet {Jensen1944}.

\bibitem[{Bodewig(1947--1948)}]{Bodewig1947-HM}
Bodewig, E., 1947--1948. {Bericht \"uber die vershiedenen Methoden zur L\"osung
  eines Systems linear Gleichungen mit reellen Koffizienten. I through V}.
  Indagationes Mathematicae. Part I, 9(4):441-452 (1947). Part II, 9(5):
  518--530 (1947). Part III, 9(5):611--621 (1947). Part IV, 10(1):24--35
  (1948). Part V, 10(1):82--90 (1948). Also published in Koninklijke
  Nederlandsche Akademie van Wetenschappen, Proceedings. Part I,
  50(7--8):930--941 (1947). Part II, 50(9--10):1104--1116 (1947). Part III,
  50(9--10):1285--1295 (1947). Part IV, 51(1--2):53--64 (1948). Part V,
  51(1--2):211--219 (1948).

\bibitem[{Bodewig(1948)}]{Bodewig1948-HM}
Bodewig, E., 1948. {MR}0024235. Mathematical Reviews 9, 471b. Review of von
  Neumann and Goldstine (1947).

\bibitem[{Bodewig(1949)}]{Bodewig1949-HM}
Bodewig, E., 1949. {MR}0028100. Mathematical Reviews 10, 405c. Review of \citet
  {Turing1948}.

\bibitem[{Bodewig(1950)}]{Bodewig1950}
Bodewig, E., 1950. The payment of the learned man. Science 112, 538--539.

\bibitem[{Bourdon(1828)}]{Bourdon1828-HM}
Bourdon, M., 1828. \'El\'emens D'Alg\`ebre, 5th Edition. Bachelier, Paris.

\bibitem[{Brandt(1935)}]{Brandt1935}
Brandt, A.~E., 1935. Uses of the progressive digit method. In: Practical
  Application of the Punch Card Method in Colleges and Universities. Columbia
  University Press, New York, pp. 423--436.

\bibitem[{Brezinski(2006)}]{Brezinski2006}
Brezinski, C., 2006. The life and work of {A}ndr\'e {C}holesky. Numerical
  Algorithms 43~(3), 279--288.

\bibitem[{Brezinski and Gross-Cholesky(2005)}]{BrezinskiGross2005}
Brezinski, C., Gross-Cholesky, M., Dec. 2005. La vie et les travaux d'{A}ndr\'e
  {L}ouis {C}holesky. Bulletin de la Soci\'et\'e des amis de la biblioth\`eque
  de l'Ecole polytechnique (SABIX) 39, 7--32, available at \url
  {http://www.sabix.org/bulletin/sabixb39.htm}.

\bibitem[{Brezinski and Wuytack(2001)}]{Brezinski2001a}
Brezinski, C., Wuytack, L. (Eds.), 2001. Numerical Analysis: Historical
  Developments in the 20th Century. North-Holland, Amsterdam.

\bibitem[{Britton(1992)}]{Britton1992}
Britton, J.~L. (Ed.), 1992. Collected Works of A. M. Turing: Pure Mathematics.
  North-Holland, Amsterdam.

\bibitem[{B\"uhler(1981)}]{Buhler1981}
B\"uhler, W.~K., 1981. Gauss; A Biographical Study. Springer-Verlag, New York.

\bibitem[{Buteo(1560)}]{Buteo1560-HM}
Buteo, J., 1560. Logistica. Paris. The author is also known as Jean Borrel.

\bibitem[{Cajori(1890)}]{Cajori1890-HM}
Cajori, F., 1890. {The Teaching and History of Mathematics in the United
  States}. Government Printing Office, Washington. Bureau of Education Circular
  of Information No.\ 3.

\bibitem[{Cajori(1894)}]{Cajori1894-HM}
Cajori, F., 1894. A History of Mathematics. MacMillan and Co., New York.

\bibitem[{Cajori(1929)}]{Cajori1929}
Cajori, F., 1929. The chequered career of Ferdinand Rudolph Hassler, first
  superintendent of the United States Coast survey; a chapter in the history of
  science in America. The Christopher Publishing House, Boston.

\bibitem[{Cardano(1545)}]{Cardano1545}
Cardano, G., 1545. Artis Magnae Sive de Regulis Algebraicis. Translated by T.
  R. Witmer as \textit {The Great Art or The Rules of Algebra\/}, MIT Press,
  Cambridge, 1968.

\bibitem[{Cassinis(1946)}]{Cassinis1946-HM}
Cassinis, G., 1946. Risoluzione dei \`sistemi di equazioni algebriche lineari.
  Rendiconti del Seminario Matematico e Fisico di Milano 17, 62--78.

\bibitem[{Cauchy(1905)}]{Cauchy1905-HM}
Cauchy, A.~L., 1905. {\OE}uvres compl\`etes d'Augustin Cauchy. Vol.~1.
  Gauthier-Villars, Paris.

\bibitem[{Chase(1980)}]{Chase1980}
Chase, G.~C., Jul. 1980. History of mechanical computing machinery. Annals of
  the History of Computing 2~(3), 198--226, reprinted from the \textit
  {Proceedings of the Association of Computing Machinery\/}, 1952, with
  corrections and foreward by I.\ Bernard Cohen.

\bibitem[{Chauvenet(1868)}]{Chauvenet1868-HM}
Chauvenet, W., 1868. Treatise on the Method of Least Squares. J. B. Lippincott
  \& Co., Philadelphia.

\bibitem[{Cholesky(2005)}]{Cholesky2005-HM}
Cholesky, A.-L., 2005. La r\'esolution num\'erique des syst\`emes d'\'equations
  lin\'eaires. Bulletin de la Soci\'et\'e des amis de la biblioth\`eque de
  l'Ecole polytechnique (SABIX) 39. Available at \url
  {http://www.sabix.org/bulletin/sabixb39.htm}. From an original written 1910.

\bibitem[{Clark(1846)}]{Clark1846-HM}
Clark, D.~W., 1846. Elements of Algebra. Harper \& Brothers, New York.

\bibitem[{Clark(2006)}]{Clark2006}
Clark, W., 2006. Academic Charisma and the Origins of the Research University.
  University of Chicago Press, Chicago.

\bibitem[{Clark(1930)}]{Clark1930}
Clark, W.~E., 1930. The \={A}ryabha\d{t}\={\i}ya of \={A}ryabha\d{t}a: An
  Ancient Indian Work on Mathematics and Astronomy. University of Chicago
  Press, Chicago, reprinted by Kessinger Publishing, 2006.

\bibitem[{Clarke(1880)}]{Clarke1880-HM}
Clarke, A.~R., 1880. Geodesy. Clarendon Press, Oxford.

\bibitem[{Clarke(2008)}]{Clarke2008}
Clarke, B.~R., 2008. Linear Models: the Theory and Application of Analysis of
  Variance. Wiley series in probability and statistics. Wiley, Hoboken.

\bibitem[{{Coast and Geodetic Survey}(1881)}]{CoastandGeodeticSurvey1881-HM}
{Coast and Geodetic Survey}, 1881. Report of the Superintendent of the Coast
  and Geodetic Survey Showing the Progress of the Work During the Fiscal Year
  Ending with June, 1878. Government Printing Office, Washington.

\bibitem[{Cohen et~al.(2006)Cohen, Lee, and Sklar}]{Cohen2006}
Cohen, D., Lee, T., Sklar, D., 2006. Precalculus: With Unit-circle
  Trigonometry, 4th Edition. Thomson Brooks/Cole, Belmont.

\bibitem[{Craig(1972)}]{Craig1972}
Craig, C.~C., 1972. Remarks concerning {P}aul {S}.\ {D}wyer. In: Tracy, D.~S.
  (Ed.), Symmetric Functions in Statistics: Proceedings of a symposium in honor
  of Professor {P}aul {S}. {D}wyer. University of Windsor, Windsor, pp. 3--10.

\bibitem[{Crelle(1864)}]{Crelle1864-HM}
Crelle, A.~L., 1864. {Rechentafeln welche alles Multipliciren und Dividiren mit
  Zahlen unter tausend ganz ersparen: Bei gr\"osseren Zahlen aber die Rechnung
  erleichtern und sicherer machen}, 2nd Edition. G. Reimer, Berlin.

\bibitem[{Crout(1941{\natexlab{a}})}]{Crout1941a}
Crout, P.~D., 1941{\natexlab{a}}. A short method for evaluating determinants
  and solving systems of linear equations with real or complex coefficients.
  Transactions of the American Institute of Electrical Engineers 60,
  1235--1241.

\bibitem[{Crout(1941{\natexlab{b}})}]{Crout1941b-HM}
Crout, P.~D., 1941{\natexlab{b}}. A short method for evaluating determinants
  and solving systems of linear equations with real or complex coefficients.
  Marchant Methods MM-182, Marchant Calculating Machine Company, Oakland.

\bibitem[{Davies(1835)}]{Davies1835-HM}
Davies, C., 1835. {Elements of Algebra: Translated from the French of M.
  Bourdon}. Wiley \& Long, New York.

\bibitem[{Dedron and Itard(1959)}]{Dedron1959}
Dedron, J., Itard, J., 1959. Math\'ematiques et Math\'ematiciens. Magnard,
  Paris.

\bibitem[{Deming(1928)}]{Deming1928}
Deming, H.~G., 1928. A systematic method for the solution of simultaneous
  linear equations. Amer.\ Math.\ Monthly 35, 360--363.

\bibitem[{Dieudonn\'e(1981)}]{Dieudonne1981}
Dieudonn\'e, J., 1981. {Von Neumann}. In: Gillispie, C.~C. (Ed.), Dictionary of
  Scientific Biography. Vol.~14. Charles Scribner's Sons, New York, pp. 89--92.

\bibitem[{Domingues(2008)}]{Domingues2008}
Domingues, J.~C., 2008. Lacroix and the Calculus. Vol.~35 of Science networks
  historical studies. Birkh\"auser, Basel.

\bibitem[{Doolittle(1888)}]{Doolittle1888-HM}
Doolittle, C.~L., 1888. A Treatise on Practical Astronomy as Applied to Geodesy
  and Navigation, 2nd Edition. John Wiley \& Sons, New York.

\bibitem[{Doolittle(1878)}]{Doolittle1878-HM}
Doolittle, M.~H., 1878. Method employed in the solution of normal equations and
  in the adjustment of a triangularization. In: \rm \citet
  {CoastandGeodeticSurvey1881-HM}. pp. 115--120.

\bibitem[{Duff et~al.(1986)Duff, Erisman, and Reid}]{Duff1986}
Duff, I.~S., Erisman, A.~M., Reid, J.~K., 1986. Direct Methods for Sparse
  Matrices. Oxford University Press, Oxford.

\bibitem[{Duncan and Kenney(1946)}]{Duncan1946}
Duncan, D.~B., Kenney, J.~F., 1946. On the Solution of Normal Equations and
  Related Topics. Edwards Brothers, Ann Arbor, pamphlet.

\bibitem[{Dunnington(2004)}]{Dunnington2004-HM}
Dunnington, G.~W., 2004. Carl Friedrich Gauss: Titan of Science. The
  Mathematical Association of America, Washington. First edition, 1955.

\bibitem[{Dwyer(1941{\natexlab{a}})}]{Dwyer1941a}
Dwyer, P.~S., 1941{\natexlab{a}}. The solution of simultaneous equations.
  Psychometrika 6~(2), 101--129.

\bibitem[{Dwyer(1941{\natexlab{b}})}]{Dwyer1941b}
Dwyer, P.~S., 1941{\natexlab{b}}. {The Doolittle Technique}. The Annals of
  Mathematical Statistics 12~(4), 449--458.

\bibitem[{Dwyer(1944)}]{Dwyer1944}
Dwyer, P.~S., 1944. A matrix presentation of least squares and correlation
  theory with matrix justification of improved methods of solution. The Annals
  of Mathematical Statistics 15~(1), 82--89.

\bibitem[{Dwyer(1945)}]{Dwyer1945}
Dwyer, P.~S., 1945. The square root method and its use in correlation and
  regression. Journal of the American Statistical Association 40, 493--503.

\bibitem[{Dwyer(1951)}]{Dwyer1951}
Dwyer, P.~S., 1951. Linear Computations. John Wiley \& Sons, New York.

\bibitem[{Eckert(1940)}]{Eckert1940-HM}
Eckert, W.~J., 1940. Punched Card Methods in Scientific Computation. Thomas J.
  Watson Astronomical Computing Bureau, New York.

\bibitem[{Euler(1771)}]{Euler1771-HM}
Euler, L., 1771. Anleitung zur Algebra. Lund.

\bibitem[{Euler(1822)}]{Euler1822-HM}
Euler, L., 1822. Elements of Algebra (Tr.\ F. Horner), 3rd Edition. Longman and
  Co., London.

\bibitem[{Farebrother(1987)}]{Farebrother1987}
Farebrother, R.~W., 1987. A memoir of the life of {M}. {H}. {D}oolittle.
  Bulletin of the Institute of Mathematics and Its Application 23, 102.

\bibitem[{Farebrother(1988)}]{Farebrother1988}
Farebrother, R.~W., 1988. Linear Least Squares Computations. Marcel Dekker, New
  York.

\bibitem[{Farebrother(1999)}]{Farebrother1999}
Farebrother, R.~W., 1999. Fitting Linear Relationships: A History of the
  Calculus of Observations 1750--1900. Springer series in statistics.
  Springer-Verlag, New York.

\bibitem[{Felippa(2001)}]{Felippa2001}
Felippa, C.~A., 2001. A historical outline of matrix structural analysis: a
  play in three acts. Computers \& Structures 79~(14), 1313--1324.

\bibitem[{Forsythe(1953)}]{Forsythe1953}
Forsythe, G.~E., 1953. Solving linear algebraic equations can be interesting.
  Bull.\ Amer.\ Math.\ Soc. 59, 299--329.

\bibitem[{Fox(1989)}]{Fox1989}
Fox, K.~A., Jan. 1989. {Agricultural Economists in the Econometric Revolution:
  Institutional Background, Literature and Leading Figures}. Oxford Economic
  Papers, New Series 41~(1), 53--70.

\bibitem[{Fox(1954)}]{Fox1954}
Fox, L., 1954. Practical solution of linear equations and inversion of
  matrices. Applied Mathematics Series~39, U.\ S.\ Bureau of Standards.

\bibitem[{Fox(1987)}]{Fox1987b}
Fox, L., 1987. {James Hardy Wilkinson. 27 September 1919 -- 5 October 1986}.
  Biographical Memoirs of Fellows of the Royal Society 33, 670--708.

\bibitem[{Fox et~al.(1948)Fox, Huskey, and Wilkinson}]{Fox1948}
Fox, L., Huskey, H.~D., Wilkinson, J.~H., 1948. Notes on the solution of
  algebraic linear systems of equations. The Quarterly Journal of Mechanics and
  Applied Mathematics 1~(2), 150--173.

\bibitem[{Frazer et~al.(1938)Frazer, Duncan, and Collar}]{Frazer1938}
Frazer, R.~A., Duncan, W.~J., Collar, A.~R., 1938. Elementary Matrices and Some
  Applications to Dynamics and Differential Equations. Cambridge University
  Press, Cambridge.

\bibitem[{Friberg(2007{\natexlab{a}})}]{Friberg2007amazing}
Friberg, J., 2007{\natexlab{a}}. Amazing Traces of a Babylonian Origin in Greek
  Mathematics. World Scientific, Hackensack.

\bibitem[{Friberg(2007{\natexlab{b}})}]{Friberg2007remarkable}
Friberg, J., 2007{\natexlab{b}}. A Remarkable Collection of Babylonian
  Mathematical Texts. Springer Verlag, New York.

\bibitem[{Gantmacher(1959)}]{Gantmacher1959-HM}
Gantmacher, F.~R., 1959. The Theory of Matrices. Vol.~1. Chelsea Publishing
  Company, New York.

\bibitem[{Gauss(1801)}]{Gauss1801-HM}
Gauss, C.~F., 1801. Disquisitiones Arithmeticae. Leipzig.

\bibitem[{Gauss(1809)}]{Gauss1809-HM}
Gauss, C.~F., 1809. Theoria Motus Corporum Coelestium in Sectionibus Conicis
  Solum Ambientium. Perthes and Besser, Hamburg.

\bibitem[{Gauss(1810)}]{Gauss1810-HM}
Gauss, C.~F., 1810. Disquisitio de elementis ellipticis {P}alladis.
  Commentationes Societatis Regiae Scientiarum Gottingensis recentiores:
  Commentationes classis mathematicae 1 (1808--1811), 1--26.

\bibitem[{Gauss(1822)}]{Gauss1822a-HM}
Gauss, C.~F., Jan. 1822. {Anwendung der Wahrscheinlichkeitsrechnung auf eine
  Aufgabe der practischen Geometrie}. Astronomische Nachrichten 1~(6), 81--86.

\bibitem[{Gauss(1823)}]{Gauss1823-HM}
Gauss, C.~F., 1823. Theoria Combinatorionis Observationum Erroribus Minimis
  Obnoxiae. Heinrich Dieterich, G\"ottingen.

\bibitem[{Gauss(1826)}]{Gauss1826-HM}
Gauss, C.~F., 1826. Supplementum theoriae combinationis observationum erroribus
  minimis obnoxiae. Commentationes Societatis Regiae Scientiarum Gottingensis
  recentiores: Commentationes classis mathematicae 6 (1823--1827), 57--98.

\bibitem[{Gauss(1857)}]{Gauss1857-HM}
Gauss, C.~F., 1857. Theory of the Motion of the Heavenly Bodies Moving about
  the Sun in Conic Sections (Tr.\ C.\ H.\ Davis). Little, Brown, and Company,
  Boston.

\bibitem[{George and Liu(1981)}]{George1981}
George, A., Liu, J. W.~H., 1981. Computer Solution of Large Sparse Positive
  Definite Systems. Prentice-Hall, Englewood Cliffs.

\bibitem[{Ghilani and Wolf(2006)}]{Ghilani2006}
Ghilani, C.~D., Wolf, P.~R., 2006. Adjustment Computations: Spatial Data
  Analysis, 4th Edition. John Wiley \& Sons.

\bibitem[{Gillings(1972)}]{Gillings1972}
Gillings, R.~J., 1972. Mathematics in the Time of the Pharaohs. MIT Press,
  Cambridge.

\bibitem[{Gillispie(1997)}]{Gillispie1997-HM}
Gillispie, C.~C., 1997. {P}ierre-{S}imon {L}aplace, 1749--1827: a life in exact
  science. Princeton University Press, Princeton.

\bibitem[{Goldstine(1972)}]{Goldstine1972}
Goldstine, H.~H., 1972. {The Computer from Pascal to von Neumann}. Princeton
  University Press, Princeton.

\bibitem[{Goldstine(1977)}]{Goldstine1977}
Goldstine, H.~H., 1977. A History of Numerical Analysis from the 16th through
  the 19th Century. Vol.~2 of Studies in the History of Mathematics and
  Physical Sciences. Springer-Verlag, New York.

\bibitem[{Golub and Van~Loan(1996)}]{Golub1996}
Golub, G.~H., Van~Loan, C.~F., 1996. Matrix Computations, 3rd Edition. The
  Johns Hopkins University Press, Baltimore.

\bibitem[{Gore(1889)}]{Gore1889}
Gore, J.~H., 1889. A Bibliography of Geodesy. Government Printing Office,
  Washington, {U}nited States Coast and Geodetic Survey. Appendix No.\ 16 ---
  Report for 1887.

\bibitem[{Gosselin(1577)}]{Gosselin1577-HM}
Gosselin, G., 1577. De Arte Magna. Paris.

\bibitem[{Grattan-Guinness(1994)}]{Grattan-GuinnessEd1994}
Grattan-Guinness, I. (Ed.), 1994. Companion Encyclopedia of the History and
  Philosophy of the Mathematical Sciences. Routledge, London.

\bibitem[{Grattan-Guinness(2004)}]{GrattanGuinness2004}
Grattan-Guinness, I., 2004. The mathematics of the past: distinguishing its
  history from our heritage. Historia Mathematica 31, 163--185.

\bibitem[{Grier(2005)}]{Grier2005}
Grier, D.~A., 2005. When computers were humans. Princeton University Press,
  Princeton.

\bibitem[{Hagen(1867)}]{Hagen1867-HM}
Hagen, G., 1867. Grundz\"uge der Wahrscheinlichkeits-Rechnung. Ernst \& Korn,
  Berlin.

\bibitem[{Hald(2007)}]{Hald2007}
Hald, A., 2007. A History of Parametric Statistical Inference from {B}ernoulli
  to {F}isher, 1713--1935. Sources and studies in the history of mathematics
  and physical sciences. Springer, New York.

\bibitem[{Hammond(1742)}]{Hammond1742}
Hammond, N., 1742. The elements of algebra. London.

\bibitem[{Hawkins(1975)}]{Hawkins1975}
Hawkins, T., 1975. Cauchy and the spectral theory of matrices. Historia
  Mathematica 2, 1--29.

\bibitem[{Hawkins(1977{\natexlab{a}})}]{Hawkins1977a}
Hawkins, T., Jun. 1977{\natexlab{a}}. Another look at {C}ayley and the theory
  of matrices. Archives Internationales d'Histoire des Sciences 27~(100),
  82--112.

\bibitem[{Hawkins(1977{\natexlab{b}})}]{Hawkins1977b}
Hawkins, T., Jul. 1977{\natexlab{b}}. Weierstrass and the theory of matrices.
  Archive for History of Exact Sciences 17~(2), 119--164.

\bibitem[{Heath(1910)}]{Heath1910}
Heath, T.~L., 1910. Diophantus of Alexandria, A Study in the History of Greek
  Algebra, 2nd Edition. Cambridge University Press, Cambridge.

\bibitem[{Heefer(2005)}]{Heefer2005}
Heefer, A., 2005. The rhetoric of problems in algebra texbooks from {P}acioli
  to {E}uler. Thesis, Ghent University Centre for Logic and Philosophy of
  Science, Ghent, \url
  {http://logica.ugent.be/albrecht/thesis/AlgebraRhetoric.pdf}.

\bibitem[{Helmert(1907)}]{Helmert1907-HM}
Helmert, F.~R., 1907. Die Ausgleichungsrechnung nach der Methode der kleinsten
  Quadrate, 2nd Edition. B. G. Teubner, Leipzig.

\bibitem[{Higham(2002)}]{Higham2002}
Higham, N.~J., 2002. Accuracy and Stability of Numerical Algorithms, 2nd
  Edition. SIAM, Philadelphia.

\bibitem[{Hogendijk(1994)}]{Hogendijk1994}
Hogendijk, J.~P., 1994. Pure mathematics in {I}slamic civilization. In: \rm
  \cite {Grattan-GuinnessEd1994}. pp. 70--79.

\bibitem[{Horsburgh(1914)}]{Horsburgh1914-HM}
Horsburgh, E.~M. (Ed.), 1914. Handbook of the Napier tercentenary celebration,
  or, Modern instruments and methods of calculation. G.\ Bell, London.
  Reprinted with a new introduction by M.\ R.\ Williams, Charles Babbage
  Institute reprint series for the history of computing, vol,\ 3, Tomash
  Publishers, Los Angeles, 1982.

\bibitem[{Hotelling(1943)}]{Hotelling1943}
Hotelling, H., 1943. Some new methods in matrix calculation. The Annals of
  Mathematical Statistics 14~(1), 1--34.

\bibitem[{Householder(1953)}]{Householder1953}
Householder, A.~S., 1953. Principles of Numerical Analysis. McGraw-Hill, New
  York.

\bibitem[{Householder(1956)}]{Householder1956}
Householder, A.~S., 1956. Bibliography on numerical analysis. Journal of the
  ACM (JACM) 3~(2), 85--100.

\bibitem[{Householder(1964)}]{Householder1964}
Householder, A.~S., 1964. The Theory of Matrices in Numerical Analysis.
  Blaisdell, New York, reprinted by Dover, New York, 1975.

\bibitem[{H{\o}yrup(2002)}]{Hoyrup2002}
H{\o}yrup, J., 2002. Lengths, Widths, Surfaces: A Portrait of Old Babylonian
  Algebra and Its Kin. Springer-Verlag, New York.

\bibitem[{{Institute of Electrical and Electronics Engineers}(1985)}]{ieee1985}
{Institute of Electrical and Electronics Engineers}, 1985. {IEEE Standard for
  Binary Floating-Point Arithmetic}. ANSI/IEEE Std 754-1985, New York.

\bibitem[{Jacobi(1857)}]{Jacobi1857-HM}
Jacobi, C., 1857. {\"Uber eine elementare Transformation eines in Bezug auf
  jedes von zwei Variablen-Systemen linearen und homogenen Ausdrucks}. Journal
  f\"ur die reine und angewandte Mathematik 53, 265--270.

\bibitem[{Jensen(1939)}]{Jensen1939}
Jensen, H., 1939. {Herleitung einiger Ergebnisse der Ausgleichsrechnung mit
  Hilfe von Matrizen}. Meddelelse~13, Geod{\ae}tisk Institut, Copenhagen.

\bibitem[{Jensen(1944)}]{Jensen1944}
Jensen, H., 1944. An attempt at a systematic classification of some methods for
  the solution of normal equations. Meddelelse~18, Geod{\ae}tisk Institut,
  Copenhagen.

\bibitem[{Johnson(1997)}]{Johnson1997}
Johnson, S., 1997. Making the arithmometer count. Bulletin of the Scientific
  Instrument Society 52, 12--21.

\bibitem[{Johnson(1905)}]{Johnson1905-HM}
Johnson, W.~W., 1905. The Theory of Errors and Method of Least Squares, 1st
  Edition. John Wiley \& Sons, New York.

\bibitem[{Jordan(1895)}]{Jordan1895-HM}
Jordan, W., 1895. Handbuch der Vermessungskunde, 4th Edition. Vol.~1. J.\ B.\
  Metzler Verlag, Stuttgart.

\bibitem[{Kangshen et~al.(1999)Kangshen, Crossley, and Lun}]{Kangshen1999}
Kangshen, S., Crossley, J.~N., Lun, A. W.-C., 1999. The Nine Chapters of the
  Mathematical Art Companion and Commentary. Oxford University Press, New York.

\bibitem[{Katz(1988)}]{Katz1988}
Katz, V.~J., 1988. {Who is the Jordan of Gauss-Jordan}. Mathematics Magazine
  61~(2), 99--100.

\bibitem[{Katz(1997)}]{Katz1997b}
Katz, V.~J., 1997. Algebra and its teaching: {A}n historical survey. The
  Journal of Mathematical Behavior 16~(1), 25--38.

\bibitem[{Katz(1998)}]{Katz1998}
Katz, V.~J., 1998. A History of Mathematics: An Introduction, 2nd Edition.
  Addison-Wesley, Reading.

\bibitem[{Kinckhuysen(1661)}]{Kinckhuysen1661-HM}
Kinckhuysen, G., 1661. Algebra ofte Stel-konst Beschreven Tot dienst van de
  leerlinghen. Passchier Van Wesbusch, Haarlem.

\bibitem[{Kloyda(1938)}]{Kloyda1938}
Kloyda, M. T. a.~K., 1938. Linear and Quadratic Equations 1550--1660. Edwards
  Brothers, Ann Arbor, {U}niversity of Michigan lithoprint dissertation.

\bibitem[{Kocinski(2004)}]{Kocinski2004}
Kocinski, J., 2004. Cracovian Algebra. Nova Science Publishers, Hauppauge.

\bibitem[{Krayenhoff(1813)}]{Krayenhoff1813}
Krayenhoff, C. R.~T., 1813. Verzameling van Hydrographische en Topographische
  Waarnemigen in Holland. Doorman en Comp., Amsterdam.

\bibitem[{Krayenhoff(1827)}]{Krayenhoff1827-HM}
Krayenhoff, C. R.~T., 1827. Pr\'ecis historique des op\'erations g\'eod\'esique
  et Astronomiques, faites en Hollande. Imprimerie de l'\'etat., La Haye.

\bibitem[{Lacroix(1800)}]{Lacroix2nd1800-HM}
Lacroix, S.~F., 1800. Elemens d'alg\`ebre, \`a l'usage de l'Ecole centrale des
  Quatre-Nations, 2nd Edition. Impr.\ de Crapelet, chez Duprat, Paris.

\bibitem[{Lacroix(1804)}]{Lacroix5th1804-HM}
Lacroix, S.~F., 1804. Elemens d'alg\`ebre, \`a l'usage de l'Ecole centrale des
  Quatre-Nations, 5th Edition. Chez Courcier, Paris.

\bibitem[{Lacroix(1818)}]{Lacroix1818-HM}
Lacroix, S.~F., 1818. Elements of Algebra (Tr.\ J.\ Farrar). University Press,
  Cambridge, Massachusetts.

\bibitem[{Laderman(1948)}]{Laderman1948}
Laderman, J., 1948. The square root method for solving simultaneous linear
  equations. Mathematical Tables and other Aids to Computation 3~(21), 13--16.

\bibitem[{Lagrange(1759)}]{Lagrange1759-HM}
Lagrange, J.~L., 1759. Researches sur la m\'etode de maximis et minimis.
  Miscellanea Taurinensia 1. Journal title varies. Reprinted in Serret, J.-A.,
  1867, {\OE}uvres de Lagrange, vol.\ 1, Gauthier--Villars, Paris, 1--16.

\bibitem[{Lay-Yong and Kangshen(1989)}]{LayYong1989}
Lay-Yong, L., Kangshen, S., 1989. Methods of solving linear equations in
  traditional {C}hina. Historia Mathematica 16, 107--122.

\bibitem[{Legendre(1805)}]{Legendre1805-HM}
Legendre, A.~M., 1805. Nouvelle m\'ethodes pour la d\'etermination des orbites
  des com\`etes. Chez Didot, Paris.

\bibitem[{Liagre(1879)}]{Liagre1879-HM}
Liagre, J. B.~J., 1879. Calcul des probabilit\'es et th\'eorie des erreurs avec
  des applications aux sciences d'observation en g\'en\'eral et a la
  g\'eod\'esie en particulier, 2nd Edition. Muquardt, Bruxelles.

\bibitem[{Libbrecht(1973)}]{Libbrecht1973}
Libbrecht, U., 1973. Chinese Mathematics in the Thirteenth Century. MIT Press,
  Cambridge.

\bibitem[{MacDuffee(1949)}]{MacDuffee1949-HM}
MacDuffee, C.~C., 1949. Vectors and Matrices, 3rd Edition. Vol.~7 of Carus
  Mathematical Monographs. The Mathematical Association of America, Menasha.
  First edition 1943.

\bibitem[{Macomber(1923)}]{Macomber1923-HM}
Macomber, G.~L., 1923. The influence of the {E}nglish and {F}rench writers of
  the sixteenth, seventeenth, and eighteenth centuries on the teaching of
  algebra, {M}.A.\ thesis, University of California, Berkeley.

\bibitem[{{Marchant Calculating Machine Company}(1941)}]{Marchant1941-HM}
{Marchant Calculating Machine Company}, Sep. 1941. {Notes on the Use of the
  Marchant Calculator for Solution of Simultaneous Equations by the Method of
  Prescott D.\ Crout as Descrbied in Marchant Method MM-182}. Marchant Methods
  MM-183, Oakland.

\bibitem[{Marguin(1994)}]{Marguin1994}
Marguin, J., 1994. Histoire des instruments et machines \`a calculer: trois
  si\`ecles de m\'ecanique pensante, 1642--1942. Hermann, Paris.

\bibitem[{Martzloff(1997)}]{Martzloff1997-HM}
Martzloff, J.-C., 1997. {A History of Chinese Mathematics (Tr.\ S.\ S.\
  Wilson)}. Springer, Berlin.

\bibitem[{Matthiessen(1878)}]{Matthiessen1878}
Matthiessen, L., 1878. {Grundz\"uge der antiken und modernen Algebra der
  litteralen Gleichungen}. B. G. Teubner, Leipzig.

\bibitem[{Merriman(1884)}]{Merriman1884-HM}
Merriman, M., 1884. A Text Book on the Method of Least Squares. John Wiley \&
  Sons, New York.

\bibitem[{Murray(1948)}]{Murray1948}
Murray, F.~J., 1948. The Theory of Mathematical Machines, 2nd Edition. King's
  Crown Press, New York.

\bibitem[{Murray(1961)}]{Murray1961}
Murray, F.~J., 1961. Mathematical Machines. Columbia University Press, New
  York, in 2 volumes.

\bibitem[{Nell(1881)}]{Nell1881-HM}
Nell, P., 1881. {Schleiermacher's Methode der Winkelausgleichung in einen
  Dreiecksnetze}. Zeitschrift f\"ur Vermussungswesen 10~(1), 1--11, 109--121.

\bibitem[{Neugebauer(1969)}]{Neugebauer1969}
Neugebauer, O., 1969. The Exact Sciences in Antiquity, 2nd Edition. Dover, New
  York.

\bibitem[{Neugebauer(1973)}]{Neugebauer1973-HM}
Neugebauer, O., 1973. Mathematische Keilschrift-Texte, 2nd Edition. Vol. 2 and
  3 in one book. Springer-Verlag, Berlin.

\bibitem[{Newton(1707)}]{Newton1707}
Newton, I., 1707. Arithmetica Universalis. London.

\bibitem[{Newton(1720)}]{Newton1720-HM}
Newton, I., 1720. Universal Arithmetick: or, a Treatise of Arithmetical
  Composition and Resolution. To which is added, Dr.\ Halley's Method of
  finding the Roots of Equations Arithmetically. Senex, Taylor, et al., London.
  Translated from the Latin by the late Mr.\ Ralphson, and Revised and
  Corrected by Mr.\ Cunn. The 1728 edition is reproduced in \citet [v.\
  2]{Whiteside1964-1967}.

\bibitem[{Nievergelt(2001)}]{Nievergelt2001}
Nievergelt, Y., 2001. A tutorial history of least squares with applications to
  astronomy and geodesy. In: \rm \cite {Brezinski2001a}. pp. 77--112.

\bibitem[{{Ordnance Survey}(1858)}]{OrdnanceSurvey1858-HM}
{Ordnance Survey}, 1858. Ordnance Trigonometrical Survey of Great Britain and
  Ireland -- Account of the observations and calculations of the principal
  triangulation, 2 vols. Eyre and Spottiswoode, London.

\bibitem[{Palmer(1912)}]{Palmer1912-HM}
Palmer, A. d.~F., 1912. The Theory of Measurements. McGraw-Hill Book Company,
  New York.

\bibitem[{Paucker(1819)}]{Paucker1819}
Paucker, M.~G., 1819. Ueber die Anwendung der Methode der kleinsten
  Quadratsumme auf physikalische Beobachtungen. Johann Friedrich Steffenhagen
  und Sohn, Mitau, pamphlet. Program zur Eršffnung des Lehrkursus auf dem
  Gymnasium illustre zu Mitau.

\bibitem[{Peletier~du Mans(1554)}]{Peletier1554-HM}
Peletier~du Mans, J., 1554. L'Algebre. Lyon.

\bibitem[{Petersen and Arbenz(2004)}]{Petersen2004}
Petersen, W., Arbenz, P., 2004. Introduction to Parallel Computing: A Practical
  Guide with Examples in C. Oxford Texts in Applied and Engineering
  Mathematics. Oxford University Press.

\bibitem[{Plackett(1972)}]{Plackett1972}
Plackett, R.~L., Aug. 1972. {Studies in the History of Probability and
  Statistics, XXIX: The Discovery of the Method of Least Squares}. Biometrika
  59~(2), 239--251.

\bibitem[{Plofker(2009)}]{Plofker2009}
Plofker, K., 2009. Mathematics in India. Princeton University Press.

\bibitem[{Pugsley(1961)}]{Pugsley1961}
Pugsley, A.~G., 1961. {Robert Alexander Frazer. 1891--1959}. Biographical
  Memoirs of Fellows of the Royal Society 7, 75--84.

\bibitem[{Rahn(1659)}]{Rahn1659-HM}
Rahn, J.~H., 1659. Teutsche Algebra. Zurich.

\bibitem[{Relf(1961)}]{Relf1961}
Relf, E.~F., 1961. {William Jolly Duncan. 1894 -- 1960}. Biographical Memoirs
  of Fellows of the Royal Society 7, 37--51.

\bibitem[{Roberts and Kynaston(1995)}]{Roberts1995}
Roberts, R., Kynaston, D. (Eds.), 1995. The Bank of England: Money, Power and
  Influence 1694--1994. Oxford University Press, Oxford.

\bibitem[{Robson(1999)}]{Robson1999}
Robson, E., 1999. Mesopotamian mathematics, 2100-1600 BC: technical constants
  in bureaucracy and education. Vol.~14 of Oxford editions of cuneiform texts.
  Oxford University Press, Oxford.

\bibitem[{Robson(2008)}]{Robson2008}
Robson, E., 2008. Mathematics in Ancient Iraq: A Social History. Princeton
  University Press, Princeton.

\bibitem[{Rolle(1690)}]{Rolle1690-HM}
Rolle, M., 1690. Trait\'e d'alg\`ebre; ou principes generaux pour resoudre les
  questions de mathematique. E. Michallet, Paris.

\bibitem[{Ross(1831)}]{Ross1831-HM}
Ross, E.~C., 1831. Elements of Algebra Translated from the French of M. Bourdon
  for the Use of the Cadets of the U. S. Military Academy. E. B. Clayton, New
  York.

\bibitem[{Rubin(1926)}]{Rubin1926}
Rubin, T., 1926. Et nytt s\"att att l\"osa normalekvationer ({A} new method of
  solving normal equations). Svensk Lantm\"ateritidskrift~(1), 3--9.

\bibitem[{Satterthwaite(1944)}]{Satterthwaite1944}
Satterthwaite, F.~E., 1944. Error control in matrix calculation. The Annals of
  Mathematical Statistics 15~(4), 373--387.

\bibitem[{Saunderson(1761)}]{Saunderson1761-HM}
Saunderson, N., 1761. Selected Parts of Professor Saunderson's Elements of
  Algebra: For the Use of Students at the Universities, 2nd Edition. London.

\bibitem[{Schappacher(2005)}]{Schappacher2005}
Schappacher, N., 2005. Diophantus of alexandria: a text and its history, \url
  {http://www-irma.u-strasbg.fr/~schappa/NSch/Publications_files/Dioph.pdf}.

\bibitem[{Schott(1881)}]{Schott1879-HM}
Schott, C.~A., 1881. {Appendix No.\ 8}. In: \rm \cite
  {CoastandGeodeticSurvey1881-HM}. pp. 92--94.

\bibitem[{Sesiano(1982)}]{Sesiano1982}
Sesiano, J., 1982. {Books IV to VII of Diophantus' Arithmetica in the Arabic
  translation attributed to Qus\d{t}\={a} ibn L\={u}q\={a}}. Springer-Verlag,
  Heidelberg.

\bibitem[{Simpson(1755)}]{Simpson1755-HM}
Simpson, T., 1755. A Treatise of Algebra, 2nd Edition. John Nourse, London.

\bibitem[{Stewart(1998)}]{Stewart1998}
Stewart, G.~W., 1998. Matrix Algorithms 1:\ Basic Decompositions. SIAM,
  Philadelphia, Pennsylvania.

\bibitem[{Stigler(1986)}]{Stigler1986}
Stigler, S.~M., 1986. The History of Statistics: The Measurement of Uncertainty
  before 1900. Harvard University Press, Cambridge.

\bibitem[{Stigler(1999)}]{Stigler1999-HM}
Stigler, S.~M., 1999. Statistics on the Table: the History of Statistical
  Concepts and Methods. Harvard University Press, Cambridge.

\bibitem[{Swartzlander(1995)}]{Swartzlander1995}
Swartzlander, E., 1995. Calculators. IEEE Annals of the History of Computing
  17~(3), 75--77.

\bibitem[{Taub(1963)}]{Taub1963}
Taub, A.~H. (Ed.), 1963. {John von Neumann Collected Works}. Macmillan, New
  York.

\bibitem[{Taussky and Todd(2006)}]{TausskyTodd2006}
Taussky, O., Todd, J., 2006. {C}holesky, {T}oeplitz and the triangular
  factorization of symmetric matrices. Numerical Algorithms 41, 197--202.

\bibitem[{{The Monthly Review}(1801)}]{MonthlyReview1801-HM}
{The Monthly Review}, May--August 1801. {Review of \textit {El\'emens
  d'Alg\`ebra\/}} 35, 470--476.

\bibitem[{Toeplitz(1907)}]{Toeplitz1907-HM}
Toeplitz, O., 1907. {Die Jacobische Transformation der quadratischen Formen von
  unendlichvielen Ver\"anderlichen}. Nachrichten von der Gesellschaft der
  Wissenschaften zu G\"ottingen, Mathematisch-Physikalische Klasse, 101--109.

\bibitem[{Tolley and Ezekiel(1927)}]{Tolley1927}
Tolley, H.~R., Ezekiel, M., Dec. 1927. {The Doolittle Method for Solving
  Multiple Correlation Equations Versus the Kelley-Salisbury ``Iteration''
  Method}. Journal of the American Statistical Association 22~(160), 497--500.

\bibitem[{Torge(2001)}]{Torge2001}
Torge, W., 2001. Geodesy, 3rd Edition. Walter de Gruyter.

\bibitem[{Turing(1948)}]{Turing1948}
Turing, A.~M., 1948. Rounding-off errors in matrix processes. The Quarterly
  Journal of Mechanics and Applied Mathematics 1~(3), 287--308, reprinted in
  \cite {Britton1992}.

\bibitem[{{United States Congress}(1886)}]{Congress1886-HM}
{United States Congress}, 1886. Testimony Before the Joint Commission \dots.
  Government Printing Office, Washington.

\bibitem[{van Eyk(1854)}]{vanEyk1854-HM}
van Eyk, J.~A., 1854. Reken-werktuigen. De Volksvlijt, 387--396.

\bibitem[{Verzuh(1949)}]{Verzuh1949}
Verzuh, F., 1949. The solution of simultaneous linear equations with the aid of
  the 602 calculating punch. Mathematical Tables and other Aids to Computation
  3~(27), 453--462.

\bibitem[{von Neumann(1954)}]{vonNeumann1954-HM}
von Neumann, J., 1954. The role of mathematics in the sciences and in society.
  In: \citet [v. 6, pp. 477--490] {Taub1963}.

\bibitem[{von Neumann and Goldstine(1947)}]{vonNeumann1947-HM}
von Neumann, J., Goldstine, H.~H., 1947. Numerical inverting of matrices of
  high order. Bulletin of the American Mathematical Society 53~(11),
  1021--1099. Reprinted in \citet [v. 5, pp. 479--557] {Taub1963}.

\bibitem[{{Washington Star}(1913)}]{WashingtonStar1913-HM}
{Washington Star}, June 28 1913. {M}yrick {H}ascall {D}oolittle obituary, Part
  1 page 7.

\bibitem[{Waugh and Dwyer(1945)}]{Waugh1945}
Waugh, F.~V., Dwyer, P.~S., 1945. Compact computation of the inverse of a
  matrix. The Annals of Mathematical Statistics 16~(3), 259--271.

\bibitem[{Werner(1883)}]{Werner1883}
Werner, W., 1883. {Ueber die Methode der ``Coast and Geodetic Survey'' zur
  Aufl\"osung von Normalgleichungen}. Der Civilingenieur 29, 116--126.

\bibitem[{Whipple(1914)}]{Whipple1914-HM}
Whipple, F. J.~W., 1914. Calculating machines. In: \citet {Horsburgh1914-HM}.
  pp. 69--123.

\bibitem[{Whiteside(1964--1967)}]{Whiteside1964-1967}
Whiteside, D.~T. (Ed.), 1964--1967. The Mathematical Works of {I}saac {N}ewton,
  2 vols. Johnson Reprint Corporation, New York and London.

\bibitem[{Whiteside(1968--1982)}]{Whiteside1968-1982}
Whiteside, D.~T. (Ed.), 1968--1982. The Mathematical Papers of {I}saac
  {N}ewton, 8 vols. Cambridge University Press, Cambridge.

\bibitem[{Wilkinson(1960)}]{Wilkinson1960}
Wilkinson, J.~H., 1960. Rounding errors in algebraic processes. In: Information
  Processing. Published by R. Oldenbourg, Munich and Butterworths, London, pp.
  44--53, proceedings of the International Conference on Information
  Processing, UNESCO, Paris, 15--20 June 1959.

\bibitem[{Williams(1982)}]{Williams1982-HM}
Williams, M.~R., 1982. Introduction. In: \citet {Horsburgh1914-HM}. pp.
  ix--xxi.

\bibitem[{Wilson(1952)}]{Wilson1952}
Wilson, E.~B., 1952. An Introduction to Scientific Research. McGraw-Hill, New
  York.

\bibitem[{Witkowski(1955)}]{Witkowski1955}
Witkowski, J., Oct. 1955. {The Life and Work of Professor Dr.\ Thaddeus
  Banachiewicz}. Acta Astronomica Series C 5, 85--94.

\bibitem[{Wright and Hayford(1906)}]{Wright1906-HM}
Wright, T.~W., Hayford, J.~F., 1906. The Adjustment of Observations by the
  Method of Least Squares with Applications to Geodetic Work. D. Van Nostrand,
  New York.

\end{thebibliography}

}


\end{document}